\documentclass[12pt,reqno]{amsart}

\usepackage{amssymb,amsmath,amsthm,amscd,ifthen,xr}

\author{M. McKee}
\address{Department of Mathematics, University of Oklahoma, Norman, OK, 73019, USA.}
\email{mark.mckee.zoso@gmail.com}
\author{A. Pasquale}
\address{Institut Elie Cartan de Lorraine (IECL UMR CNRS 7502),
Universit\'e de Lorraine, F-57045 Metz, France.}
\email{angela.pasquale@univ-lorraine.fr}

\author{T. Przebinda}
\address{Department of Mathematics, University of Oklahoma, Norman, OK, 73019, USA.}
\email{tprzebinda@gmail.com}
\title[Semisimple orbital integrals on the symplectic space]{Semisimple orbital integrals on the symplectic space\\ for a real reductive dual pair}


\def\chc{chc}
\def\n{\mathfrak n}
\def\g{\mathfrak g}
\def\z{\mathfrak z}

\def\h{\mathfrak h}

\def\s{\mathfrak s}
\def\u{\mathfrak u}
\def\k{\mathfrak k}
\def\p{\mathfrak p}
\def\q{\mathfrak q}

\def\R{\mathbb{R}}
\def\C{\mathbb{C}}

\def\Ha{\mathbb{H}}
\def\DD{\mathbb{D}}

\def\a{\mathfrak a}
\def\b{\mathfrak b}
\def\l{\mathfrak l}
\def\z{\mathfrak z}
\def\c{\mathfrak c}

\def\so{\mathfrak s_{\overline 0}}
\def\ss1{\mathfrak s_{\overline 1}}

\def\hs1{\mathfrak h_{\overline 1}}

\def\G{\mathrm{G}}
\def\Zg{\mathrm{Z}}

\def\K{\mathrm{K}}
\def\H{\mathrm{H}}
\def\M{\mathrm{M}}
\def\Zg{\mathrm{Z}}

\def\Sg{\mathrm{S}}

\def\Id{\mathrm{I}}
\def\L{\mathrm{L}}
\def\Bbb{\mathbb}

\def\N{\mathrm{N}}
\def\A{\mathrm{A}}
\def\H{\mathrm{H}}

\def\GL{\mathrm{GL}}
\def\SL{\mathrm{SL}}
\def\SO{\mathrm{SO}}
\def\Sp{\mathrm{Sp}}

\def\Og{\mathrm{O}}
\def\Ug{\mathrm{U}}

\def\Im{\mathrm{Im}}
\def\Re{\mathrm{Re}}

\def \t{\tilde}
\def \wt{\widetilde}
\newcommand{\reg}[1]{ {#1}^{reg}}

\def\W{\mathsf{W}}
\def\V{\mathsf{V}}

\def\X{\mathsf{X}}
\def\Y{\mathsf{Y}}

\def\Dc{\mathbb {D}}
\def\Zb{\mathbb {Z}}
\def\Nb{\mathbb {N}}

\def\J{\mathcal{J}}

\def\End{\mathop{\hbox{\rm End}}\nolimits}
\def\det{\mathop{\hbox{\rm det}}\nolimits}
\def\ad{\mathop{\hbox{\rm ad}}\nolimits}
\def\Ad{\mathop{\hbox{\rm Ad}}\nolimits}
\def\Hom{\mathop{\hbox{\rm Hom}}\nolimits}

\def\sgn{\mathop{\hbox{\rm sgn}}\nolimits}

\def\lim{\mathop{\hbox{\rm lim}}\nolimits}

\def\ker{\mathop{\hbox{\rm ker}}\nolimits}

\newcommand{\anticomm}[2]{\null^{#1}#2}
\newcommand{\danticomm}[2]{\null^{\anticomm{#1}{#2}}#2}

\def\U{\mathcal{U}}

\def\P{\mathcal{P}}

\def\fontindex{\arabic}

\def\fonttitre{\textsf}
\newcounter{thh}

\newtheorem{thm}[thh]{\fonttitre{Theorem}}

\newtheorem{pro}[thh]{\fonttitre{Proposition}}
\newtheorem*{pro*}{\fonttitre{Proposition}}
\newtheorem{cor}[thh]{\fonttitre{Corollary}}
\newtheorem*{coro*}{\fonttitre{Corollary}}
\newtheorem{lem}[thh]{\fonttitre{Lemma}}
\newtheorem*{rem}{\fonttitre{Remark}}

\newtheorem*{defi*}{\fonttitre{Définition}}

\newtheorem*{nota*}{\fonttitre{Notation}}
\newenvironment{prf}{\begin{proof}}{\end{proof}}
\def\muet{ \ifthenelse{\equal{a}{b}}}
\def\nn{\nonumber}

\newcommand{\thmlist}{
\renewcommand{\theenumi}{\alph{enumi}}
\renewcommand{\labelenumi}{(\theenumi)}}

\begin{document}
\thanks{The second author is grateful to the University of Oklahoma for hospitality and financial support. The third author gratefully acknowledges hospitality and financial support from the Universit\'e de Lorraine and partial support from the NSA grant H98230-13-1-0205. }

\begin{abstract}
We prove a Weyl Harish-Chandra integration formula for the action of a reductive dual pair on the corresponding symplectic space $\W$. 
As an intermediate step, we introduce a notion of a Cartan subspace and a notion of an almost semisimple element in the symplectic space $\W$. We
prove that the almost semisimple elements are dense in $\W$. Finally, we provide estimates for the orbital integrals associated with the different
Cartan subspaces
in $\W$.  
\end{abstract}
\maketitle

\tableofcontents
\section{Introduction}\label{sec.1}
Let $\W$ be a symplectic real vector space, $\Sp$ the corresponding symplectic
group and $\wt\Sp$ the metaplectic group. Let $\G, \G'\subset  \Sp$ be a
real reductive dual pair, and let $\wt\G$ and $\wt\G'$ be the preimages of $\G$ and $\G'$ in $\wt\Sp$,
respectively. Furthermore, let $\Pi \otimes \Pi'$ be an irreducible admissible representation
of $\wt\G \times \wt\G'$ in Howe's correspondence. Each such representation
$\Pi \otimes \Pi'$ is attached to a tempered distribution $f=f_{\Pi\otimes\Pi'}$ on $\W$, called the intertwining distribution of $\Pi \otimes \Pi'$, which is uniquely determined up to a scalar multiple. 
One may therefore expect that the properties of the intertwining distribution contain useful information on the representation itself. 

Let $\g$ and $\g'$ respectively denote the Lie algebras of $\G$ and $\G'$. Let $\U(\g)$ denote the universal enveloping algebra of $\g$, and let $\U(\g)^\G$ be the subalgebra of the $\G$-invariants in $\U(\g)$. Then the intertwining distribution turns out to be $\wt\G \times \wt\G'$-invariant and 
an eigendistribution of $\U(\g)^\G$ and $\U(\g')^{\G'}$, with eigenvalue equal to the infinitesimal characters of $\Pi$ and $\Pi'$, respectively;  see \cite{PrzebindaUnitary}. In other words, $f_{\Pi\otimes\Pi'}$ is an invariant eigendistibution on the symplectic space $\W$.

Harish-Chandra's method of descent is one of the main tools for studying questions involving the structure of invariant eigendistributions on a real reductive Lie algebra $\g$. Ultimately, it transfers problems of invariant harmonic analysis from the Lie algebra $\g$ to the Cartan subalgebras of $\g$, see for example \cite{BouazizOrb}. Our goal is to develop a ``method of descent'' allowing us to study invariant
eigendistributions on symplectic spaces. We use the fact that the symplectic space $\W$ is the odd part $\ss1$ of a classical real Lie superalgebra $\s$, which is constructed from the dual pair $(\G,\G')$. The method of descent
will consist in transferring problems of invariant harmonic analysis from $\ss1$ to suitably defined Cartan subspaces in $\ss1$. (Notice that the 
name ``Cartan subspace'' is also used in other contexts, for example in the theory of symmetric spaces. Our usage of this term is however different.)  

To single out the candidates for the Cartan subspaces in $\ss1$, one needs a detailed knowledge of 
the geometry of the adjoint action of $\Sg=\G\times \G'$ on $\ss1$. This step has been accomplished in \cite{PrzebindaLocal} for all irreducible reductive dual pairs except for the two following ortho-symplectic cases:
\begin{eqnarray}
&&\text{$\Og_{p,q}\times \Sp_{2n}(\R)$ with $p+q$ odd and $p+q< 2n$} \label{basic assumption-one}\\
&&\text{$\Og_{p}(\C)\times \Sp_{2n}(\C)$ with $p$ odd and $p<2n$.} 
\label{basic assumption-two}
\end{eqnarray} 
The results of this paper also include these two cases. This leads us to a definition of Cartan subspace which is slightly different from the one used in \cite{PrzebindaLocal}, as well as to the notion of almost semisimple elements in $\ss1$. The definitions agree with those from \cite{PrzebindaLocal} when the dual pair is not isomorphic to (\ref{basic assumption-one}) or (\ref{basic assumption-two}). In fact, most of the results of \cite{PrzebindaLocal} carry over.

Our first result is an analog of the Weyl Harish-Chandra formula on the symplectic space $\W$ (Theorem 
\ref{Weyl integration}). It gives the integral on $\W=\ss1$ of a continuous compactly supported function 
in terms of semisimple orbital integrals parametrized by the mutually non-conjugate Cartan subspaces
of $\ss1$. As an important intermediate step, we prove that the set of the semisimple elements in $\W$ is dense, unless the dual pair is isomorphic to (\ref{basic assumption-one}) or (\ref{basic assumption-two}). This property was stated in \cite[Proposition 6.6]{PrzebindaLocal} without a proof. Here we show explicitly how any nilpotent element of $\W$ may be approximated by semisimple elements. For the remaining cases we 
introduce the notion of almost semisimple elements and show that they are dense (Theorem \ref{density}).

Weyl Harish-Chandra formula provides an orbital integral decomposition for the integral on $\W=\ss1$ of  continuous compactly supported functions. For the analysis of tempered eigendistributions, one needs to know that this formula extends to Schwartz functions on $\W$.   
The main result of this article, Theorem \ref{3.30}, shows that any such orbital integral of a rapidly decreasing function on $\W$ is rapidly decreasing at infinity on $\h_{\overline 1}$, though it might have logarithmic growth near walls defined by some real roots. 
This guarantees the extension of Weyl Harish-Chandra formula to Schwartz functions. 
However, unlike the Lie algebra case, it is not necessarily true that a semisimple orbital integral of a Schwartz function on the symplectic space is a Schwartz function on the corresponding Cartan subspace, though it is certainly differentiable on the regular set. 
Since the radial components of invariant differential operators in this setting do not behave as well as in the case of a Lie algebra, it doesn't seem to be easy to estimate these derivatives.  

Our analysis follows closely Harish-Chandra's work. Indeed, the even part $\so$ of the Lie 
superalgebra $\s$ is the reductive Lie algebra $\g\times \g'$. So the map of $\ss1$ to $\so$ given 
by $z \mapsto z^2$ relates geometrical and analytical objects on $\ss1$ to the corresponding objects 
on $\so=\g\times \g'$. In particular, if $\hs1 \subseteq \ss1$ is a Cartan subspace, then $\hs1^2$, that is the linear span of the all the anticommutants of the elements of $\hs1$, is contained in some Cartan subalgebra of $\so$. Nevertheless we shall need several ingredients not included in the work of Harish-Chandra. For instance, in the proof of Proposition \ref{2.2pro}, we need Rossmann's theorem expressing some nilpotent orbital integrals as limits of derivatives of the Harish-Chandra semisimple integrals, \cite{RossmannNilpotent}.

The regularity properties of an invariant eigendistribution $f$ on $\W=\ss1$ can be obtained from the analysis of the system of partial differential equations it satisfies. By determining the characteristic variety of this system, we prove in Corollary \ref{4.13} that the restriction of $f$ to the set $\reg{\W}\subseteq\W$ of the regular almost semisimple elements is a smooth function. 
However, most of the times, this distribution is not regular on the entire $\W$.
We will provide some examples in our future work. 
   
Our paper is organized as follows. In section \ref{section:intertwining} we recall the construction of the intertwining distribution attached to an irreducible admissible representation in Howe's correspondence. In section \ref{section:orbital semisimple} we collect some results on orbital integrals on semisimple Lie algebras which will be needed in the Lie superalgebra situation. 
The study of the adjoint action of a dual reductive pair $\G \times \G'$ on the symplectic space $\W$ is developed in section \ref{section:Weyl-integration-W}. We use the Lie superalgebra approach from \cite{PrzebindaLocal} by considering $\W$ as the odd part $\ss1$ of a classical real Lie superalgebra $\s$. We introduce a general notion of Cartan subspace of $\ss1$ which for complex dual pairs coincides with the notion of Cartan subspace in \cite{DadokKac}. Then we prove the density in $\ss1$ of the almost semisimple elements (Theorem \ref{density}). This yields the Weyl Harish-Chandra formula on $\ss1$ (Theorem \ref{Weyl integration}). In section \ref{section:orbital symplectic} we make a detailed study of the orbital integrals through each Cartan subspace in $\ss1$. Theorem \ref{3.30} provides sharp estimates for the orbital integral of rapidly decreasing functions on $\W$. The structure theory developed in section \ref{section:Weyl-integration-W} reduces the proof of this theorem to the analysis of three special cases: the case of elliptic Cartan subspaces; the case a dual pair of the form $(\GL_1(\DD),\GL_1(\DD))$, where $\DD=\R,\C$ or $\Ha$ (the quaternions); the case of $(\GL_2(\R),\GL_2(\R))$ with centralizer of $\hs1^2$ in $\so$ equal to a fundamental Cartan subalgebra.
Finally, in section \ref{section:diffeq} we include a brief study of the regularity properties of the invariant eigendistributions from the differential equation point of view.      

Some proofs that have not been included in the main text have been collected in the Appendices A,B, C and D. 

\medskip

 \textbf{Acknowledgement:\;} The authors are grateful to the referee for the careful reading of this manuscript and for helping them remove some unnecessary assumptions from Theorem \ref{2.2}.

\section{Intertwining distributions}
\label{section:intertwining}

Let $\W$ be a vector space of finite dimension $2n$ over $\Bbb R$ with a non-degenerate symplectic form $\langle\cdot ,\cdot \rangle$. Denote by $\Sp$ the corresponding symplectic group and by $\s\p$ its Lie algebra. Recall the Cayley transform $c(y)=(y+1)(y-1)^{-1}$, \cite{HoweOscill}, and let
\begin{equation}\label{1.1}
\wt{\Sp}^c=\{\t g=(g,\xi)\in \Sp\times \C,\ \det(g-1)\ne 0,\ \xi^2=\det(i(g-1))^{-1}\}.
\end{equation}
For each $x\in\s\p$, $\langle x\cdot ,\cdot \rangle$ is a symmetric bilinear form on $\W$ with the signature $\sgn\langle x\cdot ,\cdot \rangle$ equal to the maximal dimension of a subspace where this form is positive definite minus the maximal dimension of a subspace where this form is negative definite. Set
\begin{equation}\label{1.2}
\chc(x)=2^n|\det(x)|^{-\frac{1}{2}}\,\exp(\frac{\pi}{4}i\,\sgn\langle x\cdot ,\cdot \rangle)\qquad (x\in\s\p,\ \det(x)\ne 0).
\end{equation}
(This is a Fourier transform of one of the two minimal non-zero nilpotent co-adjoint orbits in $\s\p^*$, \cite[Proposition 9.3]{PrzebindaCauchy} and \cite[Theorem 7.6.1]{Hormander}.) For two elements $(g_1,\xi_1), (g_2,\xi_2)\in \wt{\Sp}^c$, with $c(g_1)+c(g_2)$ invertible, define the product
\begin{equation}\label{1.3}
(g_1,\xi_1) (g_2,\xi_2)=(g_1g_2, \xi_1\xi_2\chc(c(g_1)+c(g_2))).
\end{equation}
\begin{thm}\label{1.4}\cite[sec. 16]{HoweOscill}
Up to a group isomorphism there is a unique connected group $\wt\Sp$ containing $\wt\Sp^c$ with the multiplication given by (\ref{1.3}) on the indicated subset of $\wt\Sp^c\times \wt\Sp^c$. The map
$$
\wt\Sp^c\ni \t g\mapsto g\in \Sp
$$ 
extends to a double covering homomorphism
$$
\wt\Sp\ni \t g\mapsto g\in \Sp.
$$
\end{thm}
Fix the unitary character $\chi(r)=e^{2\pi i r}$, $r\in \Bbb R$, and a suitably normalized Lebesgue measure $dw$ on $\W$. For $\phi_1, \phi_2\in S(\W)$, the Schwartz space  of $\W$, we have the twisted convolution, \cite{HoweOscill},
\begin{equation}\label{1.5.0}
\phi_1\natural \phi_2(w')=\int_\W\phi_1(w)\phi_2(w'-w)\chi(\frac{1}{2}\langle w,w'\rangle)\,dw \qquad (w'\in\W).
\end{equation}
Multiplication by the Lebesgue measure gives an embedding of $S(\W)$ into 
$S^*(\W)$, the space of the temperate distributions on $\W$. The twisted convolution (\ref{1.5.0}) extends to some, but not all, temperate distributions. Also, we have a $*$-operation on $S(\W)$, 
\begin{equation}\label{1.6}
\phi^*(w)=\overline{\phi(-w)}\qquad (w\in\W, \phi\in S(\W)),
\end{equation}
which does extend to $S^*(\W)$ by
\begin{equation}\label{1.7}
f^*(\phi)=f(\phi^*)\qquad (f\in S^*(\W), \phi\in S(\W)).
\end{equation}
Define the following functions
\begin{eqnarray*}
&&\Theta:\wt\Sp^c\ni \t g=(g,\xi)\to\xi\in\Bbb C,\\
&&T:\wt\Sp^c\ni \t g\to \Theta(\t g)\chi_{c(g)}\in S^*(\W),
\end{eqnarray*}
where, for $x=c(g)\in\s\p$, we have set $\chi_x(w)=\chi(\frac{1}{4}\langle x(w),w\rangle)$. 
\begin{thm}\label{1.8}\cite{HoweOscill}
The map $T$ extends to an injective continuous map $T:\wt\Sp\to S^*(\W)$ and the following formulas hold
\begin{equation*}
\begin{array}{cccc}
T(\t g_1)\natural T(\t g_2)&=&T(\t g_1 \t g_2) \qquad &(\t g_1, \t g_2\in \wt\Sp^c,\ \det(c(g_1)+c(g_2))\ne 0),\\
T(g)^*&=&T(g^{-1}) \qquad &(g\in \wt\Sp),\\
T(1)&=&\delta\,, & 
\end{array}
\end{equation*}
where $\delta$ is the Dirac delta at the origin.
\end{thm}
The oscillator representation of $\wt\Sp$ may be realized as follows. Pick a complete polarization
\begin{equation}\label{1.8.5}
\W=\X\oplus \Y
\end{equation}
and recall the Weyl transform
\begin{eqnarray}\label{1.9}
&&\mathcal K:S^*(\W)\to S^*(\X\times \X),\\
&&\mathcal K(f)(x,x')=\int_\Y f(x-x'+y)\chi(\frac{1}{2}\langle y, x+x'\rangle)\,dy, \qquad f \in S(\W)\nn
\end{eqnarray}
with the inverse given by
\begin{equation}\label{1.10}
\mathcal K^{-1}(K)(x+y)=2^{-n}\int_\X K(\frac{x'+x}{2}, \frac{x'-x}{2})\chi(\frac{1}{2}\langle x',y\rangle)\,dx'\qquad (x,x'\in \X,\ y\in \Y).
\end{equation}
Each element $K\in S^*(\X\times \X)$ defines an operator
$OP(K)\in \Hom(S(\X),S^*(\X))$ by
\begin{equation}\label{1.11}
OP(K)v(x)=\int_\X K(x,x')v(x')\,dx'.
\end{equation}
The map $OP$ extends to an isomorphism of linear topological spaces, $S^*(\X\times \X)$ and $\Hom(S(\X),S^*(\X))$. This is known as the Schwartz Kernel Theorem, \cite[Theorem 5.2.1]{Hormander}.

The oscillator representation $\omega$ of $\wt\Sp$, corresponding to the character $\chi$, is defined by
\begin{equation}\label{1.12}
\omega(g)=OP(\mathcal K(T(g))) \qquad (g\in \wt\Sp).
\end{equation}
Each operator (\ref{1.12}) extends to a unitary operator on the Hilbert space $\mathcal H_\omega=L^2(\X)$ and to an automorphism of $S^*(\X)$. The space of the smooth vectors $\mathcal H_\omega^\infty=S(\X)$, and $\Theta$ coincides with the distribution character of $\omega$.

Let $\G, \G'\subseteq \Sp$ be a dual pair. Denote by $\wt\G, \wt\G'$  the preimages of $\G$ and $\G'$ in $\wt\Sp$. Let $\Pi\otimes\Pi'$ be an irreducible admissible representation of $\wt\G\times\wt\G'$ in Howe's correspondence. 
Then, up to infinitesimal equivalence, $\Pi\otimes\Pi'$ may be realized as a subspace of $\mathcal H_\omega^\infty{}^*=S^*(\X)$. Furthermore, Howe proved in \cite{HoweTrans} that
\begin{equation}\label{1.13}
\dim(\Hom_{\wt\G\wt\G'}(\mathcal H_\omega^\infty, \Pi\otimes \Pi'))=1.
\end{equation}
Thus for each representation $\Pi\otimes \Pi'$ there is a unique (up to a scalar multiple) distribution $f_{\Pi\otimes \Pi'}\in S^*(\W)$, called the intertwining distribution, such that
\begin{equation}\label{1.14}
\C\cdot OP(\mathcal K(f_{\Pi\otimes \Pi'}))=\Hom_{\wt\G\wt\G'}(\mathcal H_\omega^\infty, \Pi\otimes \Pi')).
\end{equation}
In this exposition we presented the Schr\"odinger model of the oscillator representation $\omega$ attached to the polarization (\ref{1.8.5}), but the distribution $f=f_{\Pi\otimes \Pi'}$ does not depend on this particular model. Furthermore, as shown in \cite{PrzebindaUnitary}, the intertwining distribution $f$ satisfies the following system of differential equations:
\begin{equation}\label{1.5}
T(z)\natural f=\gamma_\Pi(z)f \qquad (z\in \U(\g)^\G),
\end{equation}
where 
\begin{equation*}
T(z)\natural f=\frac{d}{dt}T(\exp(tz))\natural f|_{t=0} \qquad (z\in\g)
\end{equation*}
and $\gamma_\Pi:\U(\g)^\G\to \C$ is the infinitesimal character of $\Pi$.
(One obtains the same system of equations replacing $\U(\g)^\G$ by $\U(\g')^{\G'}$ in (\ref{1.5}).) Thus, $f$ is an invariant eigendistribution on the symplectic space $\W$. 


\section{Singular semisimple orbital integrals on a semisimple Lie algebra}
\label{section:orbital semisimple}

This section collects some properties which will be needed in section \ref{section:orbital symplectic} for studying the convergence of semisimple orbital integrals on a symplectic space. These properties adapt to our situation some well known results of Harish-Chandra; see \cite[sections 4 and 5]{HC-57Fourier} or \cite[Part 1, Ch. 3]{Varada}. 

As in \cite{HC-57Fourier}, $\g$ stands for a real semisimple Lie algebra and  
$\G$ is the (connected) adjoint group of $\g$. Fix a semisimple element $z\in\g$.
Let $\h\subseteq \g$ be a Cartan subalgebra containing $z$.  Let $\z=\g^z$ denote the centralizer of $z$ in $\g$ and let $\Zg\subseteq \G$ be the centralizer of $z$ in $\G$. Then $\Zg$ is a real reductive group with the Lie algebra $\z$. Let $\c$ denote the center of $\z$. Then $\c\subseteq \h\subseteq \z$ and $\c$ does not depend on $\h$. 
Recall that a root $\alpha$ for the pair $(\h_\C, \g_\C)$ is said to be real if $\alpha(\h)\subseteq \R$ and imaginary if $\alpha(\h) \subseteq i\R$;
if $\alpha$ is neither real nor imaginary, it is said to be complex. We fix a set of positive roots of $(\h_\C, \g_\C)$. For $\h \subseteq \q\subseteq \g$, let $\pi_{\g/\q}$ denote the product of all positive roots such that the corresponding root spaces do not occur in $\q_\C$.
Thus $\pi_{\g/\h}$ is the product of all positive roots and
\begin{equation}\label{2.1}
\pi_{\g/\z}=\prod_{\alpha>0,\,\alpha(z)\ne 0}\alpha\,.
\end{equation}
Let $\reg\c\subseteq \c$ be the subset where $\pi_{\g/\z}$ does not vanish. Moreover, denote by $d(g\Zg)$ a positive $\G$-invariant measure on the quotient space $\G/\Zg$. This measure does exist because both $\G$ and $\Zg$ are unimodular. 

In the proofs of this section we shall need some notation concerning constant coefficients differential operators on $\g$. 
For $x\in\g$ let $\partial(x)$ denote the corresponding directional derivative
\begin{equation}\label{2.3}
\partial(x)\psi(y)=\frac{d}{dt}\psi(y+tx)|_{t=0} \qquad (x,y\in\g,\ \psi\in C^\infty(\g))\,.
\end{equation}
Then $\partial$ extends to an isomorphism from the symmetric algebra $\mathcal S(\g)$ of $\g_\C$ onto the algebra of the constant coefficient differential operators on $\g$. 
If $u\in \mathcal S(\g)$ is homogeneous of degree $d$, then the differential operator $\partial(u)$ is said to be of degree $d$.
Fix a $\g$-invariant non-degenerate symmetric bilinear form $B(\cdot ,\cdot )$ on the real 
vector space $\g$. (So $B$ is a nonzero constant multiple of the Killing form.) This form extends to the complexification $\g_\C$ and provides an isomorphism of $\g_\C$ with the dual $\g_\C^*$, and hence a $\G_\C$-equivariant identification of $\mathcal S(\g)$ with $\P(\g)$, the algebra of the complex valued polynomial functions on $\g$. Here $\G_\C$ denotes the adjoint group of $\g_\C$.
Furthermore, let $D$ be a differential operator on $\g$ and let $x \in \g$. Then there exists a unique element $p \in \mathcal S(\g)$ such that 
$$D\psi(x)=\big(\partial(p)\psi\big)(x) \qquad (\psi\in C^\infty(\g))\,.$$
The constant coefficient differential operator $\partial(p)$ is called the local expression of $D$ at $x$, and is denoted by $D|_x$.

In the following, $|\cdot |$ is a fixed norm on the real vector space $\g$ and $g.x$ denotes the adjoint action of $g\in\G$ on $x\in\g$.  
\begin{thm}\label{2.2}
For every $M\geq 0$ there is $N\geq 0$ such that 
\begin{equation*}
\sup_{x\in \reg\c} (1+|x|)^M |\pi_{\g/\z}(x)| \int_{\G/\Zg}(1+|g.x|)^{-N}\,d(g\Zg)<\infty.
\end{equation*}
\end{thm}
\begin{prf}
Let $\z''=[\z,\z]$ be the derived algebra of $\z$. Then $\z=\c\oplus \z''$ and $z\in\c$. Let $\b\subseteq \z''$ be a fundamental Cartan subalgebra.
Then $\h=\c\oplus \b$ is a Cartan subalgebra of $\g$, see \cite[Part 1, Ch 1, Lemma 1]{Varada}.

Let $\pi_{\z/\h}$ denote the product of all the positive roots for $(\h_\C,\g_\C)$ such that the corresponding root spaces occur in $\z_\C$. 
According to Harish-Chandra, \cite[Theorem 3, page 568]{HC-64c}, there is a non-zero constant $C\in\R$ such that for all $f\in C_c^\infty(\z)$ and all $x\in\c$,
\begin{equation}\label{2.4}
\partial(\pi_{\z/\h})\left.\left(\pi_{\z/\h}(x+y)\int_{\Zg/\H}f(x+g''.y)\,d(g''\H))\right)\right|_{y=0}=C f(x).
\end{equation}
Here $y\in \b$, $\pi_{\z/\h}(x+y)=\pi_{\z/\h}(y)$ and the differential operator $\partial(\pi_{\z/\h})$ acts with respect to the variable $y$. 
Let $\psi\in C_c^\infty(\g)$. Suppose $x+y\in \reg\h$. Then
\begin{eqnarray}\label{2.50}
&&\pi_{\z/\h}(x+y)\int_{\G/\H}\psi(g.(x+y))\,d(g\H)\\
&=&\int_{\G/\Zg}\left(\pi_{\z/\h}(x+y)\int_{\Zg/\H}\psi(g.(x+g''.y))\,d(g''\H))\right)\,d(g\Zg).\nn
\end{eqnarray}
Since $\h$ is fundamental in $\z$, there are no real roots, so the expression in the parenthesis in \eqref{2.50} is the Harish-Chandra orbital integral, which depends on the parameter $g\Zg$. 
By a theorem of Harish-Chandra, \cite[Theorem 3, page 225]{HC-57Fourier}, this  is a smooth compactly supported function of the variable $g\Zg$, while $x+y$ is allowed to vary over a compact subset of $\reg\h$. 
Hence, Leibnitz Rule, \cite[8.11.2]{DieudonneElements}, shows that
\begin{eqnarray}\label{2.51}
&&\partial(\pi_{\z/\h})\left(\pi_{\z/\h}(x+y)\int_{\G/\H}\psi(g.(x+y))\,d(g\H)\right)\\
&=&\int_{\G/\Zg}\partial(\pi_{\z/\h})\left(\pi_{\z/\h}(x+y)\int_{\Zg/\H}\psi(g.(x+g''.y))\,d(g''\H))\right)\,d(g\Zg).\nn
\end{eqnarray}
By taking the limit if $y\to 0$ in \eqref{2.51} and applying (\ref{2.4}) we see that for $x\in\reg\c$,
\begin{eqnarray}\label{2.52}
&&\partial(\pi_{\z/\h})\left.\left(\pi_{\z/\h}(x+y)\int_{\G/\H}\psi(g.(x+y))\,d(g\H)\right)\right|_{y=0}\\
&=&\int_{\G/\Zg}\left.\partial(\pi_{\z/\h})\left(\pi_{\z/\h}(x+y)\int_{\Zg/\H}\psi(g.(x+g''.y))\,d(g''\H))\right)\right|_{y=0}\,d(g\Zg)\nn\\
&=&C\int_{\G/\Zg} \psi(g.x)\,d(g\Zg).\nn
\end{eqnarray}
Let $W(\g,\h)$ denote the Weyl group of $(\h_\C,\g_\C)$ and let $\P(\h)^{W(\g,\h)}$ be the space of the $W(\g,\h)$-invariants in $\P(\h)$. Let $\eta\in \P(\h)^{W(\g,\h)}$ be viewed as a differential operator of degree zero, i.e. multiplication by $\eta$. Then
\begin{equation}\label{2.6}
\left.\left(\partial(\pi_{\g/\h})\eta\right)\right|_{0}= \eta(0)\partial(\pi_{\g/\h}),
\end{equation}
where the left hand side  is the local expression at zero of the composition of the two differential operators and $\eta(0)\in\C$ is the value of $\eta$ at $0\in\h$. Indeed, there is a unique $\zeta\in \P(\h)$ such that the left hand side of (\ref{2.6}) is equal to $\partial(\zeta)$. Since $\pi_{\g/\h}$ is $W(\g,\h)$-skew invariant, so is $\zeta$. Lemma 10 in \cite[page 100]{HC-57DifferentialOperators} implies that there is $\xi\in \P(\h)^{W(\g,\h)}$ such that $\zeta=\xi\pi_{\g/\h}$. But, by the definition of $\zeta$, the degree of $\zeta$ is less or equal to the degree of $\pi_{\g/\h}$. Hence $\xi$ is a constant. 

Let $\eta^+$ be the sum of the homogeneous components  of positive degrees of $\eta$, so that $\eta=\eta(0)+\eta^+$. Then
$\partial(\pi_{\g/\h})(\eta^+\pi_{\g/\h})(0)=0$. Hence, $\partial(\pi_{\g/\h})(\eta\pi_{\g/\h})(0)=\eta(0)\big(\partial(\pi_{\g/\h})\pi_{\g/\h}\big)(0)$.
Therefore, 
$$
\xi\big(\partial(\pi_{\g/\h})\pi_{\g/\h}\big)(0)=\partial(\zeta)\pi_{\g/\h}(0)=\partial(\pi_{\g/\h})(\eta\pi_{\g/\h})(0)=\eta(0)\big(\partial(\pi_{\g/\h})\pi_{\g/\h}\big)(0),
$$
which shows that $\xi=\eta(0)$, and (\ref{2.6}) follows.

We apply (\ref{2.6}) with $\g$ replaced by $\z$ and $\eta$ replaced by $\pi_{\g/\z}$, which is $W(\z,\h)$-invariant. Then for $x\in\reg\c$ and $y\in \h''$,
\begin{equation}\label{2.7}
\left.\left(\partial(\pi_{\z/\h})\pi_{\g/\z}(x+y)\right)\right|_{y=0}=
\left.\left(\pi_{\g/\z}(x)\partial(\pi_{\z/\h})\right)\right|_{y=0}.
\end{equation}
Since $\pi_{\g/\h}=\pi_{\g/\z}\pi_{\z/\h}$, equations (\ref{2.52}) and (\ref{2.7})  show that for $\psi\in C_c^\infty(\g)$ and $x\in\reg \c$,
\begin{eqnarray}\label{2.8}
\partial(\pi_{\z/\h})\left.\left(\pi_{\g/\h}(x+y)\int_{\G/\H}\psi(g.(x+y))\,d(g\H)\right)\right|_{y=0}=C\pi_{\g/\z}(x)\int_{\G/\Zg} \psi(g.x)\,d(g\Zg).
\end{eqnarray}
Since $C\ne 0$, (\ref{2.8}) together with a theorem of Harish-Chandra, \cite[Theorem 3, page 225 and Lemma 25, page 232]{HC-57Fourier} imply that for any $M\geq 0$ there is a continuous seminorm $\nu$ on $S(\g)$ such that for all $\psi\in C_c^\infty(\g)$ and all $x\in \reg\c$
\begin{eqnarray}\label{2.9}
(1+|x|)^M\left|\pi_{\g/\z}(x)\int_{\G/\Zg} \psi(g.x)\,d(g\Zg)\right|
\leq \nu(\psi).
\end{eqnarray}
Notice that for each $x\in \reg\c$ the formula
\begin{eqnarray*}
(1+|x|)^M|\pi_{\g/\z}(x)|\int_{\G/\Zg} \psi(g.x)\,d(g\Zg)
\end{eqnarray*}
defines a positive Borel measure on $\g$. Hence, \cite[Lemma 8, page 37]{Varada} shows that there are finite non-negative constants $C$, $N$ such that
\begin{eqnarray}\label{estimate:elliptic}
(1+|x|)^M|\pi_{\g/\z}(x)|\int_{\G/\Zg} (1+|g.x|)^{-N}\,d(g\Zg)
\leq C \qquad (x\in\reg\c).
\end{eqnarray}
This completes the proof of the theorem.
\end{prf}
\begin{cor}
Theorem \ref{2.2} extends to the case where $G$ is a real reductive Lie group with finitely many connected components. 
\end{cor}
\begin{prf}
Suppose first that $\G$ is connected. The Lie algebra of $\G$ decomposes as $\g=\z(\g)\oplus \g''$ where $\z(\g)$ is the center of $\g$ and
$\g''=[\g,\g]$ is semisimple. Write $x=x_{\z(\g)}+x''$ for the corresponding decomposition of $x \in \g$. 
The center $\Zg(\G)$ of $\G$ has Lie algebra $\z(\g)$, and $\G/\Zg(\G) \cong \Ad(\G)$ is a connected semisimple Lie group with Lie 
algebra $\g''$. A Cartan subalgebra in $\g$ is of the form $\h=\z(\g)\oplus \h''$ where $\h''$ is a Cartan subalgebra of $\g''$. The roots of $(\h_\C,\g_\C)$ coincide with those of $(\h''_\C,\g''_\C)$ and vanish on $\z(\g)$.
Let $z=z_{\z(\g)}+z'' \in \h$. If $\Zg$ is the centralizer of $z$ in $\G$, then $\Zg/\Zg(\G)$ is the centralizer of $z''$ in $\G/\Zg(\G)$. The Lie algebra of $\Zg/\Zg(\G)$ is $\z'=(\g'')^{z''}$. Write $\z'=\c'\oplus \z''$, where $\c'$ and $\z''$ are respectively the center and the derived Lie algebra of $\z'$. Then the Lie algebra $\z=\g^z$ of $\Zg$ decomposes as $\z=\z(\g)\oplus \z'=\c\oplus \z''$ where $\c=\z(\g)\oplus \c'$ is the center of $\z$. For $x=x_{\z(\g)}+x'' \in \c$ we have 
$\pi_{\g/\z}(x)=\pi_{\g''/\z'}(x'')$. In particular, the regularity of $x$ depends only on the regularity of its semisimple part $x''$. Furthermore, $g.x=x_{\z(\g)}+ g.x''$. Thus the inequality stated in Theorem \ref{2.2} for $\G$ follows from the corresponding inequality for $\G/\Zg(\G)$.

Suppose now that $\G$ has finitely many connected components. Let $\G_0$ be the connected component of the identity. Then $\G_0$ is an open normal subgroup of $\G$. Let $\Zg$ and $\Zg^0$ be respectively the centralizers of $z$ in $\G$ and in $\G_0$. Then $\Zg^0=\Zg \cap \G_0$ is normal in $\Zg$ and $\G_0\Zg$ is normal in $\G$. 
Let $\pi:\G \to \G/\Zg$ be the canonical projection. Then $\iota:\G_0/\Zg^0 \to \G/\Zg$ given by $\iota(g\Zg^0)=g\Zg$ is an open embedding so that
$\pi^{-1} \iota(\G_0/\Zg^0)=\G_0\Zg$. In particular, we can normalize the invariant measure on $\G_0/\Zg^0 \cong \iota(\G_0/\Zg^0)$ so that it agrees with the one induced by $\G/\Zg$. Moreover, if $x_1 \G_0 \Zg, \dots, x_n\G_0 \Zg$, with $x_j \in \G$, are distinct elements of $\G/\G_0 \Zg$, 
then $\G/\Zg=\cup_{j=1}^n s_j\iota(\G_0/\Zg^0)$ (disjoint union). Hence, if $f$ is a sufficiently regular measurable function on $\g$, then 
$$
\int_{\G/\Zg} f(g.x)\, d(d\Zg)=\sum_{j=1}^n \int_{\G_0/\Zg^0} f(s_jg.x)\, d(g\Zg^0)\,.
$$
If $\|\cdot\|$ denotes the operator norm in some finite dimensional representation of $\g$ on a Hilbert space, then for all $j$ there is a constant $C_j>0$ so that for all $g \in \G$ we have
$\|g.x\|=\|s_j^{-1}s_j g.x\| \leq C_j \|s_jg.x\|$\,. Since all norms on $\g$ are equivalent, we obtain for a suitable constant $C>0$,
\begin{eqnarray*}
\pi_{\g/\z}(x) \int_{\G/\Zg} (1+|g.x|)^{-N} \, d(g\Zg) &=&\pi_{\g/\z}(x) \sum_{j=1}^n \int_{\G_0/\Zg^0} (1+|s_jg.x|)^{-N} \, d(g\Zg)\\
 &\leq&  C \pi_{\g/\z}(x) \int_{\G_0/\Zg^0} (1+|g.x|)^{-N} \, d(g\Zg^0)\,.  
\end{eqnarray*}
Hence the inequality of Theorem \ref{2.2} for $\G$ follows from the corresponding inequality for $\G_0$. 
\end{prf}

We keep the notation introduced at the beginning of this section. Recall the derived Lie algebra $\z''=[\z,\z]$ of $\z$. Assume $\z''$ is isomorphic to a symplectic Lie algebra $\s\p_{2n}(\R)$ or $\s\p_{2n}(\C)$. Let $\Zg''$ be the commutator subgroup of $\Zg$. Then $\Zg''$ is a Lie subgroup of $\Zg$ with the Lie algebra $\z''$. Let $\Og''\subseteq \z''$ be a non-zero minimal nilpotent orbit. Fix an element $n\in\Og''$ and let $\Zg^n\subseteq \Zg$ be the centralizer of $n$. Then
$$
\Zg/\Zg^n=\Zg''/\Zg''{}^n=\Og''.
$$
Since there is an invariant measure on $\Og''$ we see that the group $\Zg^n$ is unimodular. Let $d(g\Zg^n)$ denote an invariant measure on the quotient $\G/\Zg^n$. 
\begin{pro}\label{2.2pro}
Suppose the Cartan subalgebra $\h\subseteq \g$ is fundamental. Then, under the above assumptions, for every $M\geq 0$ there is $N\geq 0$ such that 
\begin{equation*}
\sup_{x\in \reg\c} (1+|x|)^M |\pi_{\g/\z}(x)| \int_{\G/\Zg^n}(1+|g.(x+n)|)^{-N}\,d(g\Zg^n)<\infty.
\end{equation*}
\end{pro}
\begin{prf}
The proof is entirely analogous to the proof of Theorem \ref{2.2}. For reader's convenience, we repeat it with necessary modifications.

We consider first the case when $\z''$ is isomorphic to $\s\p_{2n}(\R)$. The Cartan subalgebra $\h\cap \z''$ is elliptic.
Let $\pi_{\z/\h}^{short}$ be the product of the positive short roots of $\h$ in $\z_\C$.
A result of Rossmann, \cite[Corollary 5.4, page 283]{RossmannNilpotent}, shows that there is a non-zero constant $C\in\C$ such that for all $f\in C_c^\infty(\z)$ and all $x\in\c$,
\begin{eqnarray}\label{2.4pro}
&&\partial(\pi_{\z/\h}^{short})\left.\left(\pi_{\z/\h}(x+y)\int_{\Zg/\H}f(x+g''.y)\,d(g''\H))\right)\right|_{y=0}\\
&=&C \int_{\Zg/\Zg^n}f(x+g''.n)\,d(g''\Zg^n)\,.\nn
\end{eqnarray}
Here $y \in \h\cap \z''$ approaches zero from a fixed connected component of the set of the regular elements of $\h\cap \z''$, and the constant $C$, which is not relevant for us, depends on that component. Moreover, the differential operator $\partial(\pi_{\z/\h}^{short})$ acts with respect to the variable $y\in \h \cap \z''$. 
(We shall verify the assumptions of  \cite[Corollary 5.4, page 283]{RossmannNilpotent} in Appendix A.)

Let $\psi\in C_c^\infty(\g)$. Then (\ref{2.4pro}) shows that for $x\in\reg\c$, 
\begin{eqnarray}\label{2.5pro}
&&\partial(\pi_{\z/\h}^{short})\left.\left(\pi_{\z/\h}(x+y)\int_{\G/\H}\psi(g.(x+y))\,d(g\H)\right)\right|_{y=0}\\
&=&\int_{\G/\Zg}\left.\partial(\pi_{\z/\h}^{short})\left(\pi_{\z/\h}(x+y)\int_{\Zg/\H}\psi(g.(x+g''.y))\,d(g''\H)\right)\right|_{y=0}\,d(g\Zg)\nn\\
&=&C\int_{\G/\Zg^n} \psi(g.(x+n))\,d(g\Zg^n).\nn
\end{eqnarray}
Let $\eta\in \P(\h)^{W(\z,\h)}$ be viewed as a differential operator of degree zero, i.e. multiplication by $\eta$. Then
\begin{equation}\label{2.6pro}
\left.\left(\partial(\pi_{\z/\h}^{short})\eta\right)\right|_{0}= \eta(0)\partial(\pi_{\z/\h}^{short}),
\end{equation}
where the left hand side  is the local expression at zero of the composition of the two differential operators and $\eta(0)\in\C$ is the value of $\eta$ at $0\in\h$. Indeed, there is a unique $\zeta\in \P(\h)$ such that the left hand side of (\ref{2.6pro}) is equal to $\partial(\zeta)$. The polynomial $\zeta$ transforms under the Weyl group $W(\z,\h)$ the same way $\pi_{\z/\h}^{short}$ does. In particular $\zeta$ times the product of the long roots is $W(\z,\h)$-skew invariant.  Lemma 10 in \cite[page 100]{HC-57DifferentialOperators} implies that there is $\xi\in \P(\h)^{W(\z,\h)}$ such that $\zeta$ times the product of the long roots is equal to $\xi\pi_{\z/\h}$. Thus $\zeta=\xi\pi_{\z/\h}^{short}$.
But, by the definition of $\zeta$, the degree of $\zeta$ is less or equal to the degree of $\pi_{\z/\h}^{short}$. Hence $\xi$ is a constant. 

Let $\eta^+$ be the sum of the homogeneous components  of positive degrees of $\eta$, so that $\eta=\eta(0)+\eta^+$. Then
$\partial(\pi_{\z/\h}^{short})(\eta^+\pi_{\z/\h}^{short})(0)=0$. Hence, 
$$
\partial(\pi_{\z/\h}^{short})(\eta\pi_{\z/\h}^{short})(0)=\eta(0)\partial(\pi_{\z/\h}^{short})\pi_{\z/\h}^{short}(0).
$$
Therefore, 
$$
\xi\partial(\pi_{\z/\h}^{short})\pi_{\z/\h}^{short}(0)=\partial(\zeta)\pi_{\z/\h}^{short}(0)=\partial(\pi_{\z/\h}^{short})(\eta\pi_{\z/\h}^{short})(0)=\eta(0)\partial(\pi_{\z/\h}^{short})\pi_{\z/\h}^{short}(0),
$$
which shows that $\xi=\eta(0)$, and (\ref{2.6pro}) follows.

We apply (\ref{2.6pro}) with $\eta=\pi_{\g/\z}$, which is $W(\z,\h)$-invariant. Then for $x\in\reg\c$ and $y\in \z''\cap \h$,
\begin{equation}\label{2.7pro}
\left.\left(\partial(\pi_{\z/\h}^{short})\pi_{\g/\z}(x+y)\right)\right|_{y=0}=
\left.\left(\pi_{\g/\z}(x+y)\partial(\pi_{\z/\h}^{short})\right)\right|_{y=0}.
\end{equation}
Since $\pi_{\g/\h}=\pi_{\g/\z}\pi_{\z/\h}$, (\ref{2.5pro}) and (\ref{2.7pro})  show that for $\psi\in\C_c^\infty(\g)$ and $x\in\reg \c$,
\begin{eqnarray}\label{2.8pro}
&&\partial(\pi_{\z/\h}^{short})\left.\left(\pi_{\g/\h}(x+y)\int_{\G/\H}\psi(g.(x+y))\,d(g\H)\right)\right|_{y=0}\\
&=&C\pi_{\g/\z}(x)\int_{\G/\Zg^n} \psi(g.(x+n))\,d(g\Zg^n).\nn
\end{eqnarray}
Since $C\ne 0$, (\ref{2.8pro}) together with a theorem of Harish-Chandra, \cite[Theorem 2, page 207 and Lemma 25, page 232]{HC-57Fourier} imply that for any $M\geq 0$ there is a seminorm $\nu$ on $S(\g)$ such that for all $\psi\in C_c^\infty(\g)$ and all $x\in \reg\c$
\begin{eqnarray}\label{2.9pro}
(1+|x|)^M\left|\pi_{\g/\z}(x)\int_{\G/\Zg^n} \psi(g.(x+n))\,d(g\Zg^n)\right|
\leq \nu(\psi).
\end{eqnarray}
Notice that for each $x\in \reg\c$ the formula
\begin{eqnarray*}
(1+|x|)^M|\pi_{\g/\z}(x)|\int_{\G/\Zg^n} \psi(g.(x+n))\,d(g\Zg^n)
\end{eqnarray*}
defines a positive Borel measure on $\g$. Hence, \cite[Lemma 8, page 37]{Varada} shows that there are finite non-negative constants $C$, $N$ such that
\begin{eqnarray}\label{2.10}
(1+|x|)^M|\pi_{\g/\z}(x)|\int_{\G/\Zg^n} (1+|g.(x+n)|)^{-N}\,d(g\Zg^n)
\leq C \qquad (x\in\reg\c).
\end{eqnarray}
This ends the proof for the case when $\z''$ is isomorphic to $\s\p_{2n}(\R)$.

Suppose now that $\z''$ is isomorphic to $\s\p_{2n}(\C)$. Let $\Og\subseteq \z''$ be the unique non-zero minimal nilpotent orbit and let $\mu_\Og$ be a (non-zero) invariant  measure on $\Og$. Let $f\to \hat f$ denote the Fourier transform on $\z''$ defined with respect to the Killing form. Then 
$\hat\mu_\Og$ may be computed as in the real case; see for instance \cite[proof of Proposition 9.3]{PrzebindaCauchy}. In particular the restriction of $\hat\mu_\Og$ to $\h''=\h\cap \z''$ is the reciprocal of the product of the long roots of $\h''$ in $\z''_\C$. Hence, $\pi_{\z''/\h''}\hat\mu_\Og$ is the product of the short roots. Therefore, with an appropriate normalization of all the measures involved, for a test function $\psi\in S(\z'')$ 
\begin{eqnarray}\label{complex sp}
\int_{\z''}\psi\,d\mu_{\Og}&=&\int_{\z''}\hat\mu_{\Og}(x)\hat\psi(x)\,dx\\
&=&\frac{1}{|W(\Zg'',\H'')|}\int_{\h''}(\hat\mu_{\Og}(x)\pi_{\z''/\h''}(x))\left(\overline\pi_{\z''/\h''}(x)\int_{\Zg''/\H''}\hat\psi(g.x)\,d(g\H'')\right)\,dx\nn\\
&=& C\;\partial(\hat\mu_{\Og}\pi_{\z''/\h''})\left(\overline\pi_{\z''/\h''}(y)\int_{\Zg''/\H''}\psi(g.y)\,d(g\H'')\right)_{y=0},\nn
\end{eqnarray}
where $C$ is a constant. The last equality follows from the fact that for complex semisimple groups the Harish-Chandra orbital integral commutes with Fourier transform, see \cite[ Lemma 35.1, page 198]{HC-75}. The proof for the real symplectic case carries over, with (\ref{2.4pro}) replaced by (\ref{complex sp}).
\end{prf}
\section{A Weyl Harish-Chandra formula on the odd part of an ordinary classical Lie superalgebra}
\label{section:Weyl-integration-W}

In this section we restrict ourselves to real reductive dual pairs $(\G,\G')$ which are irreducible, that is no nontrivial direct sum decomposition 
of the symplectic space $\W$ is simultaneously preserved by $\G$ and $\G'$. Irreducible reductive dual pairs have been classified by Howe 
\cite{howetheta} and are of two types:
\smallskip

Type I: \quad 
$(\Og_{p,q}, \Sp_{2n}(\R))$\,, \quad  $(\Og_p(\C), \Sp_{2n}(\C))$\,, \quad  $(\Ug_{p,q}, \Ug_{r,s})$\,, \quad        
$(\Og^*_{2n}, \Sp_{p,q})$

Type II: \quad $(\GL_n(\Dc), \GL_m(\Dc))\,  \quad (\Dc\in \{\R,\C,\Ha\})$  

\smallskip 
\noindent (In the type I case, one should add to each $(\G,\G')$ the pair $(\G',\G)$, if not already listed.)  

According to \cite[section 2]{PrzebindaLocal}, we may view an irreducible reductive dual pair $(\G,\G')$ acting on the symplectic space $\W$ as a supergroup $(\Sg,\s)$. Here $\Sg$ is a Lie group isomorphic to the direct product $\G\times \G'$, and $\s=\so\oplus \ss1$ is a Lie superalgebra with even part $\so$ equal to the Lie algebra of $\Sg$ and odd part $\ss1$ equal to $\W$. The Lie superalgebra $\s$ can be realized as a subalgebra of the Lie superalgebra $\End(\V)$ of the endomorphisms of a finite dimensional 
$(\Zb/2\Zb)$-graded vector space $\V=\V_{\overline 0}\oplus \V_{\overline 1}$ over $\Dc$. The vector space $\V=\V_{\overline 0}\oplus \V_{\overline 1}$ is called the defining module for $(\Sg,\s)$.

We denote by $[\cdot,\cdot]$ the superalgebra Lie bracket on $\s$. Its restriction to $x,y \in \ss1$ coincides with the anticommutator
$$
\{x,y\}=xy+yx \in \so
$$
of the endomorphisms $x$ and $y$ of $\V$. If $x \in \so$ and $y \in \s_\alpha$ ($\alpha \in \Zb/2\Zb$) is homogeneous, then $[x,y]=xy-yx \in \s_\alpha$ is the usual commutator in $\End(\V)$. The adjoint action of $\Sg$ on $\so$, $\ss1$ and $\s$ is given by conjugation by elements of $\Sg$. In the following, $\Sg$-orbits always mean adjoint orbits. The $\Sg$-orbit of $x$ is denoted by $\Sg.x$.

Recall that there is a non-degenerate, $\Sg$-invariant, bilinear form $\langle\cdot ,\cdot \rangle$ on $\s$ such that 
$\langle\cdot ,\cdot \rangle|_{\so\times\so}$ is symmetric, 
$\langle\cdot ,\cdot \rangle|_{\ss1\times\ss1}$ is skew-symmetric, and $\so$ is orthogonal to $\ss1$ with respect to $\langle\cdot ,\cdot \rangle$. Furthermore,  there is an automorphism $\theta$ of $\s$, unique up to conjugation by $\Sg$, such that the restriction of $\theta$ to $\so$ is a Cartan involution assiciated with $\langle\cdot ,\cdot \rangle|_{\so\times\so}$ and the restriction of $\theta$ to $\ss1$ is a negative compatible complex structure, \cite[Theorem 2.19]{PrzebindaLocal}. In particular
\begin{equation}\label{3.17}
-\langle\theta\cdot ,\cdot \rangle
\end{equation} 
is a non-degenerate positive definite symmetric bilinear form on $\s$. Moreover $\theta$ preserves each member of the dual pair.

An element $x \in \s$ is called \textit{semisimple} (resp., \textit{nilpotent}) if $x$ is semisimple (resp., nilpotent) as 
an endomorphism of $\V$. Thus $x$ is nilpotent if there is $k \in \Nb$ so that $x^k=0$. If $\DD=\R$ or $\C$, then 
$x$ is semisimple provided it is diagonalizable over $\C$;  if $\DD=\Ha$, we consider $\V$ as a vector space over $\C$ by identifying $\Ha$ as a $2$-dimensional right vector space over $\C$ (see e.g. \cite[p. 61]{knappLie}).

Every element $x \in \End(\V)$ admits a Jordan decomposition $x=x_s+x_n$ with $x_s$ semisimple, $x_n$ nilpotent and so that $x_sx_n=x_nx_s$;
see e.g. \cite[Ch. VII, \S 5, no. 9]{BourbakiAlgebraII}. If $x,y \in \ss1$, then $x_s \in \ss1$ and $x_n \in \ss1$. Moreover, an element $y\in \ss1$ anti-commutes with $x$ if and only if it anti-commutes with $x_s$ and $x_n$. See \cite[Theorem 4.1]{PrzebindaLocal}.

Suppose $\g$ is a semisimple Lie algebra consisting of endomorphisms of a vector space $\M$. It is a classical property (see e.g. \cite[\S 6, n$\null^0 3$]{BourbakiGALCh1}) that $x \in \g$ is nilpotent (or semisimple) if and only if $\ad x$ is nilpotent (or semisimple) in $\End(\g)$.
The elements of $\ss1$ have a similar property.  

\begin{lem} \label{lemma:ad-nilpotent}
Let $x \in \ss1$. Then $x$ is nilpotent (resp. semisimple) if and only if $\ad x$ is a nilpotent (resp. semisimple) endomorphism of $\s$. 
\end{lem}
\begin{prf}
If $x$ is nilpotent, there is $k \in \Nb$ so that $x^k=0$. We prove that $(\ad x)^{2k-1}=0$. Indeed, $[x,y]=xy\pm yx$ for $y\in \s_\beta$ 
($\beta\in \Zb/2\Zb$). Hence, for $m \in \Nb$,
$$(\ad x)^{m}(y)=\sum_{\stackrel{a,b>0}{a+b=m}} C_{a,b} x^a y x^b\,,$$ 
where $C_{a,b} \in \Bbb Z$.
If $m=2k-1$ and $b<k$, then $a=2k-1-b>k-1$ and $x^a=0$. So $(\ad x)^{2k-1}(y)=0$.
Thus $\ad x$ is nilpotent in $\End(\s)$.

If $x$ is semisimple, there is a basis $\mathcal B=\{e_1,\dots,e_n\}$ of $\V_\C$ and elements $\lambda_1,\dots, \lambda_n\in\C$ so that $x(e_j)=\lambda_j e_j$ for all $j=1,\dots,n$. Let $E_{i,j}$ denote the $n\times n$-matrix with all entries $0$ but the $ij$-th which is equal to $1$. Then $x$ admits matrix representation $x=\sum_{j=s}^n \lambda_s E_{s,s}$ with respect to $\mathcal B$. Hence
 $$\ad x(E_{i,j})=\sum_{j=s}^n \lambda_s [E_{s,s},E_{i,j}]=(\lambda_i - \lambda_j) E_{i,j}\,.$$
Indeed $E_{s,s} \in \End(\V)_{\overline 0}$, so $[E_{s,s},E_{i,j}]$ is the usual commutator of matrices. 
It follows that $\ad x$ is diagonalizable in $\End(\V_\C)$, hence semisimple in $\End(\V)$.

Suppose now that $x \in \ss1$ is semisimple. Let $x=x_s+x_n$ be its Jordan decomposition with $x_s\in \ss1$ semisimple and $x_n \in \ss1$ nilpotent. 
Then $\ad x=\ad x_s+\ad x_n$. By the above, $\ad x_s$ is semisimple and $\ad x_n$ is nilpotent. 
Then $\ad x_n=\ad x-\ad x_s$ is at the same time semisimple and nilpotent, so $\ad x_n=0$.  
This means that $x_n$ anticommutes with all elements of $\ss1$, i.e. $x_n \in \anticomm{\ss1}{\ss1}$. By \cite[Lemma 13.0']{PrzebindaLocal}, 
$\anticomm{\ss1}{\ss1}=0$. Thus $x_n=0$, and $x=x_s$ is semisimple. 

The case where $\ad x$ is nilpotent can be proven in a similar way. In fact, one does not need to invoke \cite[Lemma 13.0']{PrzebindaLocal}.
In fact, one now gets $\ad x_s=0$, i.e. $x_s$ anticommutes with all elements of $\ss1$. In particular, it anticommutes with itself.
So $x_s^2=1/2\{x_s,x_s\}=0$. Thus $x_s$ is at the same time nilpotent and semisimple, i.e. $x_s=0$ and $x=x_n$ is nilpotent.
\end{prf}

As in the classical Lie algebra case, the property of being semisimple or nilpotent can be stated in terms of adjoint orbits. 

\begin{lem} \label{lemma:ss-n-orbits}
Let $x \in \ss1$. Then $x$ is semisimple if and only if the orbit $S.x$ is closed.
It is nilpotent if and only if $0 \in \overline{S.x}$.
\end{lem} 
\begin{prf}
The first part is \cite[Theorem 4.3]{PrzebindaLocal}. To prove the second, suppose first that $x\in \ss1$ is nilpotent. Then $0=x_s \in \overline{S.x}$ by \cite[Theorem 4.2]{PrzebindaLocal}. Conversely, suppose $0 \in \overline{S.x}$. Let $g_n\in S$ so that $0=\lim_{n\to\infty} g_n.x$. Then $0=\lim_{n\to \infty} (g_n.x)^2=\lim_{n\to \infty} g_n.x^2$. Hence $0 \in \overline{S.x^2}$. Since $x^2=1/2\{x,x\} \in \so$, the classical property of nilpotent adjoint orbits in reductive Lie algebras, e.g. \cite[Lemma 4]{Varada}, proves that $x^2$ is nilpotent. Thus $x$ is nilpotent.  
\end{prf}

Lemma \ref{lemma:ss-n-orbits} shows that our definitions of semisimple and nilpotent elements agree with those in \cite{DadokKac}. 
Following \cite{DadokKac}, we say that a semisimple element $x \in \ss1$ is \textit{regular} if nonzero and $\dim(S.x) \geq \dim(S.y)$ for all semisimple $y \in \ss1$. 

Let $x \in \ss1$ be fixed. The \textit{anticommutant} of $x$ in $\ss1$ is 
$$
\anticomm{x}{\ss1}=\{y \in \ss1:\{x,y\}=0\}\,.
$$
The \textit{double anticommutant} of $x$ in $\ss1$ is 
\begin{equation} \label{eq:doubleanticomm}
\danticomm{x}{\ss1}= \bigcap_{y \in \anticomm{x}{\ss1}} \anticomm{y}{\ss1}\,.
\end{equation}
Notice that $x \in \danticomm{x}{\ss1}$ and $\anticomm{x}{\ss1}=\anticomm{(\danticomm{x}{\ss1})}{\ss1}$.
We define a \textit{Cartan subspace} $\hs1$ of $\ss1$ as the double anticommutant of a regular semisimple element $x \in \ss1$.
We denote by $\reg{\hs1}$ the set of regular elements in $\hs1$.
The \textit{Weyl group} $\W(\Sg,\hs1)$ is the quotient of the stabilizer of $\hs1$ in $\Sg$ by the subgroup which acts trivially on 
$\hs1$. 

The following lemma relates our definition of Cartan subspace in $\ss1$ with the one given in \cite{DadokKac}.
\begin{lem} \label{lemma:CartanDadokKac}
Let $x \in \ss1$ and let $c_x=\{y \in \ss1:[\so, y] \subseteq [\so,x]\}$\,. Then $c_x=\danticomm{x}{\ss1}$.
\end{lem}
\begin{prf}
By definition, $y \in \danticomm{x}{\ss1}$ if and only if $\{z,y\}=0$ for all $z \in \anticomm{x}{\ss1}$, i.e. if and only if 
$\anticomm{x}{\ss1} \subseteq \anticomm{y}{\ss1}$. By \cite[Lemma 3.5]{PrzebindaLocal}, for any $y \in \ss1$, we have 
$[\so,y]=(\anticomm{y}{\ss1})^\perp$, where $\perp$ denotes the orthogonal with respect to $\langle\cdot,\cdot\rangle|_{\ss1\times\ss1}$.
Hence $\anticomm{x}{\ss1} \subseteq \anticomm{y}{\ss1}$ if and only if $[\so,y] \subseteq [\so,x]$.
\end{prf}

Notice that the above definition of regular semisimple elements in $\ss1$ is slightly different from the one given in \cite{PrzebindaLocal}. Correspondingly, our notion of Cartan subspaces is slightly less restrictive, and allows us to treat also the cases (\ref{basic assumption-one}) and 
(\ref{basic assumption-two}). The study of the density properties of semisimple elements given in Proposition 6.6 of that paper
needs to be modified. It is carried over in this section and the main density result is Theorem \ref{density} below. Observe however that outside the two cases (\ref{basic assumption-one}) and (\ref{basic assumption-two}), all our definitions agree with those in \cite{PrzebindaLocal}.
Most of the results in \cite[section 6]{PrzebindaLocal} carry over for all dual pairs. 
Exceptions are the second equality in Lemma 6.9(b), and Proposition 6.10(b) and (c).
Also, in the proof of \cite[Corollary 6.21]{PrzebindaLocal} the reference to (4.4) and (6.10.b) has to be replaced by (6.11.c). 

The next proposition collects the results from \cite{PrzebindaLocal} we shall need. 
As in \cite{PrzebindaLocal}, we exclude from our considerations the dual pair $(\Og_1,\Sp_{2n})$ over $\R$ or $\C$, as $\ss1$ does not contain any nonzero semisimple element in this case. (Indeed, for every $0\neq x \in \ss1$, we have $x^2 \in \so=0$, so $x$ is nilpotent.)

\begin{pro} \label{pro:summaryLocal}
Suppose $(\Sg,\s)$ corresponds to a dual pair different from $(\Og_1,\Sp_{2n})$ over $\R$ or $\C$. Then:
\begin{enumerate}
\thmlist
\item There are finitely many $\Sg$-conjugacy classes of Cartan subspaces in $\ss1$. 
\item Any two elements of a Cartan subspace $\hs1 \subseteq \ss1$ commute as endomorphisms of $\V$. 
\item Every element of $\hs1$ is semisimple and every semisimple element of $\ss1$ is contained in a Cartan subspace. 
\item The Weyl group $\W(\Sg,\hs1)$ is a finite group of linear automorphisms of $\hs1$.
\end{enumerate}
\end{pro}
\begin{prf} (Sketch)\; Let $x \in \ss1$ and let $\W \subseteq \V$ be a $(\Zb/2\Zb)$-graded subspace preserved by $x$ (i.e. $x(\W)\subseteq \W$). The element $x$, or the pair $(x,\W)$, is said to be decomposable if there are two proper non-zero $(\Zb/2\Zb)$-graded subspaces $\W_1, \W_2$ of $\W$ which are preserved by $x$ and so that $\W=\W_1\oplus\W_2$ (orthogonal direct sum if $(\Sg,\s)$ is of type I). Otherwise, $(x,\W)$ is said to be indecomposable. For every $x \in \ss1$, the pair $(x,\V)$ is the direct sum of indecomposable elements, see \cite[Theorem 5.1]{PrzebindaLocal}. The complete list (up to conjugation by an element of $\Sg$) of Cartan subspaces $\hs1 \subseteq\ss1$ and Weyl groups $\W(\Sg,\hs1)$ such that $\hs1$ contains a regular semisimple indecomposable element is given in \cite[Proposition 6.2]{PrzebindaLocal}. In the general case, a Cartan subspace $\hs1\subseteq \ss1$ induces a direct sum decomposition of $\V$ into indecomposable pairs $(x_j, \V^j)$ with $x_j \in \hs1$. Then $\hs1$ will be direct sum of Cartan subspaces in each of these indecomposable pieces. They can be read off from the above lists. Likewise, the Weyl group $\W(\Sg,\hs1)$ can be reconstructed out of those of the indecomposable pieces; see \cite[page following the proof of Proposition 6.2]{PrzebindaLocal}. The description 
of each Cartan subspace as double anticommutant of a regular semisimple element is given in \cite[pp. 487--498]{PrzebindaLocal}, based on the classification of \cite[Proposition 6.2]{PrzebindaLocal}. From the explicit expressions obtained from this case-by-case analysis, one also deduces the remaining statements.
\end{prf}

The special nature of the dual pairs (\ref{basic assumption-one}) and (\ref{basic assumption-two}) is explained in the next proposition.

Let $\V=\V_{\overline{0}} \oplus \V_{\overline{1}}$ be the defining module of $(\Sg,\s)$. For a fixed semisimple $x \in \ss1$ we can decompose 
$\V=\V^0 \oplus \V^+$, where $\V^0=\ker(x)$ and $\V^+=x\V$ are $(\Zb/2\Zb)$-graded subspaces. We shall denote by $(\Sg(\V^+),\s(\V^+))$ and $(\Sg(\V^0),\s(\V^0))$ be the supergroups obtained by restricting $(\Sg,\s)$ to $\V^+$ and $\V^0$.
\begin{pro} \label{prop:s1V0}
Let $\V=\V^0 \oplus \V^+$ and $x\in \ss1$ regular semisimple be as above. Then
$\ss1(\V^0)\neq \{0\}$ if and only if $(\Sg,\s)$ corresponds to one of the dual pairs in (\ref{basic assumption-one}) or 
(\ref{basic assumption-two}). In cases (\ref{basic assumption-one}) and (\ref{basic assumption-two}), we have $\dim_\DD\ss1(\V^0)=
2n-(p+q-1)$.
\end{pro}
\begin{prf}
Suppose first that $(\G,\G')$ is not $(\Og_{p,q},\Sp_{2n}(\R))$ or $(\Og_{p}(\C),\Sp_{2n}(\C))$. We shall prove that $\ss1(\V^0)=0$.

We begin by proving that $\ss1(\V^0)$ does not contain any nonzero semisimple element. For this, it suffices to prove that $\ss1(\V^0)$ does not contain any nonzero semisimple $y$ so that $(y,\V^0)$ is indecomposable. (This notion has been introduced in the proof of the Proposition \ref{pro:summaryLocal} above.) Indeed, by \cite[Theorem 5.1]{PrzebindaLocal}, each pair $(y, \V^0)$ with $0\neq y$ semisimple is a direct sum of nonzero indecomposable elements, which are themselves nonzero semisimple elements in $\ss1(\V^0)$.
 
We argue by contradiction. Suppose there is $0\neq y\in \ss1(\V^0)$ semisimple so that $(y,\V^0)$ is indecomposable. Then, up to similarity,
$y$ is a non-zero element of the Cartan subspaces listed in \cite[Proposition 6.2]{PrzebindaLocal}. A case-by-case inspection shows that 
there exists an element in $\so(\V^0)$ not commuting with $y$. Hence 
\begin{equation}\label{eq:dim-less}
\dim \so(\V^0)^y < \dim \so(\V^0)\,.
\end{equation}
Notice that $\so^x \subseteq \so^{x^2}=\so(\V^0)\oplus \so(\V^+)^{x^2}$ and $\so(\V^0) \subseteq \so^x$.
Hence 
\begin{equation}\label{eq:dimsox}
\so^x =\so(\V^0)\oplus (\so(\V^+)^{x^2})^x\,.
\end{equation} 
As $y|_{\V^+}=0$, we have $\so(\V^+) \subseteq \so^y$.  
Therefore 
\begin{equation}\label{eq:dimsoxy}
\so^{x+y}=\so^x \cap \so^y=\so(\V^0)^y\oplus (\so(\V^+)^{x^2})^x\,.
\end{equation}
As $x|_{\V^0}=0$, the elements $x$ and $y$ commute. So $x+y$ is semisimple. From (\ref{eq:dim-less}), (\ref{eq:dimsox}) and (\ref{eq:dimsoxy}) , 
we get $\dim \so^{x+y} < \dim \so^x$. 
Using the relation $\dim \so^z=\dim \so - \dim [\so,z]$, we therefore obtain that $\dim [x,\so] <\dim [x+y,\so]$, contradicting that $x$ is regular.
Thus $\ss1(\V^0)$ cannot contain any nonzero semisimple element.

Now $(\Sg(\V^0),\s(\V^0))$ is a supergroup coming from a dual pair. If $\ss1(\V^0)\neq 0$, the only case where $\ss1(\V^0)$ contains no 
nonzero semisimple elements is when $\Sg(\V^0)=\Og_1 \times \Sp_{2m}$ over $\R$ or $\C$ for some $m \in \Bbb N$. This means that $(\Sg,\s)$ corresponds to $(\Og_{p,q},\Sp_{2n}(\R))$ or $(\Og_{p}(\C),\Sp_{2n}(\C))$. This case was excluded. Hence $\ss1(\V^0)=0$.

We now turn to the case where $(\G,\G')=(\Og_{p,q},\Sp_{2n}(\R))$ or $(\Og_{p}(\C),\Sp_{2n}(\C))$. 

Notice that the proof above also shows that if $\ss1(\V^0)\neq 0$ then $\Sg(\V^0)=\Og_1 \times \Sp_{2m}$ over $\DD=\R$ or $\C$, for some $1\leq m \leq n$. Thus $\ss1(\V^0)=\DD^{2m}$. (The precise value of $m$ as a function of the given $p$, $q$ and $n$ will be computed at the end of the proof.) Moreover $\dim_\DD (\V^0 \cap \V_{\overline 0})=1$ and $\dim_\DD (\V^0 \cap \V_{\overline 1})=2m$.
Since $x|_{\V^+}$  is injective, we have $x(\V_{\overline 0} \cap \V^+) \subseteq \V_{\overline 1} \cap \V^+$ and 
$x(\V_{\overline 1} \cap \V^+) \subseteq \V_{\overline 0} \cap \V^+$. Hence $\dim_\DD (\V_{\overline 0} \cap \V^+)=\dim_\DD (\V_{\overline 1} \cap \V^+)$. It follows that
\begin{equation}\label{eq:dimV+}
\dim_\DD(\V_{\overline 1} \cap \V^+)=\dim_\DD (\V_{\overline 0} \cap \V^+)=\dim_\DD \V_{\overline 0}- \dim_\DD (\V_{\overline 0} \cap \V^0)=p+q-1\,.
\end{equation}
(Here and in the following, we set $q=0$ when $\DD=\C$.)
Therefore
$$2n=\dim_\DD \V_{\overline 1}\geq \dim_\DD (\V_{\overline 1} \cap \V^+) =p+q-1\,.$$
This proves that $\ss1(\V^0)= 0$ if $p+q>2n$.

Suppose now $p+q <2n$. Notice that $\ker(x^2)=\ker(x)$ as $x$ is semisimple. Hence $\ker(x^2|_{\V_{\overline 0}})=\ker(x|_{\V_{\overline 0}})=
\V_{\overline 0} \cap \V^0$. The element $x$ is regular and $x^2|_{\V_{\overline 0}}$ belongs to a Cartan subalgebra 
of $\so(\V_{\overline 0})$, the Lie algebra
of $\G=\Og_{p,q}$ or $\Og_{p}(\C)$. So $\dim\big(\ker(x^2|_{\V_{\overline 0}})\big)=\dim(\V_{\overline 0} \cap \V^0)=1$ if and only if $p+q$ is odd.

Finally, if $p+q <2n$ and $p+q$ is odd, then 
\begin{equation} \label{eq:dimV01}
\dim_\DD (\V^0 \cap \V_{\overline 1})= \dim_\DD \V_{\overline 1} - \dim_\DD (\V^+\cap \V_{\overline 1})=2n-(p+q-1)>0\,.
\end{equation}
Thus $\V^0 \cap \V_{\overline 1}\neq 0$ and $\V^0 \cap \V_{\overline 1}\neq 0$, which yields $\ss1(\V^0)\neq 0$.

Formula (\ref{eq:dimV01}) also implies that the constant $2m$ obtained before as $\dim_\DD (\V^0 \cap \V_{\overline 1})$ is equal to 
$2n-(p+q-1)$. Thus $\dim_\DD \ss1(\V^0)=2n-(p+q-1)$.
\end{prf}

Recall that the rank of a reductive Lie group $\G$, denoted by ${\rm rank\,} \G$, is the complex dimension of a Cartan subalgebra of $\g_\C$, the complexification of the Lie algebra of $\G$.  

\begin{cor} \label{cor:dim+basic-assumptions}
If $x \in \hs1$ is regular, then $\hs1=\danticomm{x}{\ss1}$. Moreover,
\begin{equation}\label{eq:relation-dim}
\dim [\so,x]+\dim \hs1+\dim \ss1(\V^0)=\dim \ss1\,.
\end{equation}
All Cartan subspaces of $\ss1$ have the same dimension. It is equal to $\min\{{\rm rank\,} \G, {\rm rank\,} \G'\}$ if $(\G,\G')$ is the dual pair associated with $(\Sg,\s)$. 
\end{cor} 
\begin{prf}
The first part is an immediate consequence of Lemma \ref{lemma:CartanDadokKac}. Indeed $\dim S.x=\dim[\so,x]$. Hence, if $\hs1=\anticomm{y}{\ss1}$ and $x \in \hs1$ is regular, then $[\so,x]=[\so,y]$ and $c_x=c_y$. For the dimension formula (\ref{eq:relation-dim}), we have 
$\hs1=\ss1(\V^+)^x$ by \cite[Proposition 6.10(a)]{PrzebindaLocal}. Since $\ker(x|_{\V^+})=0$, \cite[Theorem 4.4(a)]{PrzebindaLocal} gives $\dim \ss1(\V^+)^x=\dim \anticomm{x}{(\ss1(\V^+))}$. Moreover, $\anticomm{x}{\ss1}=\ss1(\V^0)\oplus \anticomm{x}{(\ss1(\V^+))}$. By \cite[Theorem 3.5]{PrzebindaLocal}, $\dim[\so,x]=\dim \ss1 - \dim \anticomm{x}{\ss1}$. This proves (\ref{eq:relation-dim}). 
The formula for the dimension of the Cartan subspaces is from \cite[Theorem 4.4(c)]{PrzebindaLocal}.
\end{prf}

\begin{rem}
Because of Lemmas \ref{lemma:ss-n-orbits} and \ref{lemma:CartanDadokKac}, the fact that the elements of a Cartan subspace are semisimple can also be deduced from \cite[Lemma 2.1]{DadokKac}. Since \cite{DadokKac} considers actions on complex vector spaces, this lemma can be directly applied when $\ss1$ admits a complex structure commuting with the group action of $\Sg$. This is equivalent to the fact that $\V$ is a vector space over $\C$. (Recall that we are considering $\V$ as a vector space over $\C$ if $\DD=\Ha$.) 

When $\V$ is a vector space over $\R$, it is necessary to consider the complexification $(\Sg_\C,\s_\C)$ of the supergroup $(\Sg,\s)$. 
 The list of irreducible dual pairs gives the following two possibilities:
$$
\begin{array}{ll}
\Sg=\GL_n(\R)\times \GL_m(\R) \qquad & \Sg_\C=\GL_n(\C)\times \GL_m(\C) \\
\Sg=\Og_{p,q}\times \Sp_{2n}(\R) & \Sg_\C=\Og_{p+q}(\C)\times \Sp_{2n}(\C)
\end{array}
$$
Notice that if $\V$ is the defining module of $(\Sg,\s)$, then the defining module of $(\Sg_\C,\s_\C)$ is $\V_\C$. Hence, for an element 
$x\in \ss1$, to be semisimple means that $x$ is semisimple as an element of $(\ss1)_\C=(\s_\C)_{\overline 1} \subseteq \End(\V_\C)$.

Let $\hs1=\danticomm{x}{\ss1}$ be a Cartan subspace in $\ss1$. Then $\hs1 \subseteq \danticomm{x}{(\s_\C)_{\overline{1}}}$. We can deduce from \cite{DadokKac} that $\hs1$ consists of semisimple elements if we prove that $\danticomm{x}{(\s_\C)_{\overline{1}}}$ is a Cartan subspace in $(\s_\C)_{\overline{1}}$, i.e. that $x$ is regular in $(\s_\C)_{\overline{1}}$. Since $\Sg_\C$ is a complexification of the Lie group $\Sg$, the last part of Corollary \ref{cor:dim+basic-assumptions} shows that the (real) dimension of a Cartan subspace in $(\s_\C)_{\overline 1}$ is twice the dimension of a Cartan subspace of $\ss1$. Apply 
formula (\ref{eq:relation-dim}) to both $(\Sg,\s)$ and $(\Sg_\C,\s_\C)$. Corollary \ref{cor:dim+basic-assumptions} also shows that 
$\dim \ss1(\V^0)$ is the same for all spaces $\V^0=\ker(x)$ with $x$ regular semisimple in $\ss1$. This constant doubles in the complexification. It follows that a semisimple element $y \in (\s_\C)_1$ is regular if and only if $\dim [(\s_\C)_{\overline 0},y]=2\dim [\so,x]$, where $x$ is the fixed regular semisimple element of $\ss1$. Since $(\s_\C)_{\overline 0}=(\so)_\C$, it follows that $x$ is regular in $(\s_\C)_{\overline 1}$.

Even if this argument allows us to deduce the semisimplicity of elements in the Cartan subspaces from the general results of \cite{DadokKac}, it uses many properties \cite{PrzebindaLocal}. So it still (strongly) depends of the case-by-case analysis made in \cite{PrzebindaLocal}. 
\end{rem}

Recall the defining module $\V=\V_{\overline 0}\oplus \V_{\overline 1}$ for $(\Sg,\s)$. 
Let $\hs1\subseteq \ss1$ be a Cartan subspace, and let $\V^0\subseteq \V$ be the intersection of the 
kernels of all the elements on $\h_{\overline 1}$. Equivalently, $\V^0=\ker(x)$ if $\hs1=\danticomm{x}{\ss1}$. 
Define
\begin{equation} \label{eq:tildeh1}
\t\h_{\overline 1}=\ss1(\V^0)\oplus \h_{\overline 1}\,.
\end{equation}
By Corollary \ref{cor:dim+basic-assumptions}, $\t\h_{\overline 1} \neq \h_{\overline 1}$ if and only if $(\Sg,\s)$ corresponds to a dual pair in (\ref{basic assumption-one}) or (\ref{basic assumption-two}).

We say that an element $x\in\ss1$ is \textit{almost semisimple}  if it is conjugate to an element of $\t\h_{\overline 1}$ for a Cartan subspace $\h_{\overline 1}\subseteq \ss1$.

Our next main result (Theorem \ref{density}) shows that the almost semisimple elements are dense in $\ss1$. Notice that semisimple elements are almost semisimple by Proposition \ref{pro:summaryLocal}. We now study how nilpotent elements can be approximated by almost semisimple elements. 
For this, we first consider some suitable root spaces decompositions in $\s$.
  
Let $\hs1\subseteq \ss1$ be a Cartan subspace and let $\hs1^2\subseteq \so$ be the linear span of all the anticommutants of the elements of $\hs1$. Then $\hs1^2$ is generated by the pairwise commuting semisimple elements $x^2=1/2\{x,x\}$ with $x \in \hs1$. Moreover, $\hs1^2$ is contained in some Cartan subalgebra of $\so$. 
\begin{lem}\label{theta stable csa}
There is an element $s\in \Sg$ such that $\theta(s.\hs1^2)=s.\hs1^2$.
\end{lem}
\begin{prf}
It follows from \cite[Proposition 6.10(f)]{PrzebindaLocal} that we may choose $\overline j, \overline k\in\{\overline 0, \overline 1\}$ so that  $\hs1^2|_{\V_{\overline j}}$ is a Cartan subalgebra of $\so|_{\V_{\overline j}}$ and $\hs1^2|_{\V_{\overline k}}$ is a Cartan subalgebra of $\so|_{\V_{\overline k}\cap\V^+}$. (Here $\so|_{\V_{\overline j}}$ is the Lie algebra of the isometries of $\V_{\overline j}$ in the type I case and all endomorphisms of $\V_{\overline j}$ in the type II case, and similarly for $\so|_{\V_{\overline k}\cap\V^+}$.) Furtheremore, by \cite[Theorem 4.4(b)]{PrzebindaLocal}, there is a direct sum decomposition
\[
\V_{\overline k}=(\V_{\overline k}\cap \V^0)\oplus (\V_{\overline k}\cap \V^+),
\]
which is orthogonal in the type II case. By a theorem of Chevalley, \cite[Corollary, page 100]{HC-56}, there is an element $s_{\overline k}\in \Ad(\Sg|_{\V_{\overline k}\cap\V^+})$ such that 
\[
\theta|_{V_{\overline k}}(s_{\overline k}\hs1^2|_{\V_{\overline k}})=s_{\overline k}\hs1^2|_{\V_{\overline k}}\,.
\]
Let $\epsilon$ be the linear isomorphism of $\V_{\overline k}$ which is multiplication by $1$ on $\V_{\overline k}\cap \V^+$ and 
by $-1$ on $\V_{\overline k}\cap \V^0$. Then conjugation by $\epsilon$ is an involutive automorphism $\sigma$ of $\so|_{\V_{\overline k}}$ having set of the $\sigma$-fixed points
equal to
\[
(\so|_{\V_{\overline k}})^\sigma=\so|_{\V_{\overline k}\cap\V^0}\oplus \so|_{\V_{\overline k}\cap\V^+}.
\]
Since $s_{\overline k}\hs1^2|_{\V_{\overline k}}$ is a $\sigma\theta|_{V_{\overline k}}$-stable
subspace of $ \so|_{\V_{\overline k}\cap\V^+}$, we know from the theory of the symmetric Lie algebras, \cite[Lemma 3, p. 337]{Matsuki79}, that there is an element $t_{\overline k}\in \Ad(\Sg|_{\V_{\overline k}})$ such that the Cartan involution $\theta_{\overline k}=t_{\overline k}^{-1}\theta|_{\V_{\overline k}}t_{\overline k}$ commutes with $\sigma$ and 
$t_{\overline k}(s_{\overline k}\hs1^2|_{\V_{\overline k}})=
s_{\overline k}\hs1^2|_{\V_{\overline k}}$\,. Hence 
$\theta_{\overline k}(s_{\overline k}\hs1^2|_{\V_{\overline k}})=s_{\overline k}\hs1^2|_{\V_{\overline k}}$\,. Thus,
\[
\theta(t_{\overline k}s_{\overline k}\hs1^2|_{\V_{\overline k}})=t_{\overline k}\theta_{\overline{k}}(s_{\overline k}\hs1^2|_{\V_{\overline k}})=t_{\overline k}s_{\overline k}\hs1^2|_{\V_{\overline k}}\,.
\]

Again, by the theorem of Chevalley, \cite[Corollary, page 100]{HC-56}, there is an element $s_{\overline j}\in \Ad(\Sg|_{\V_{\overline j}})$ such that 
\[
\theta(s_{\overline j}\hs1^2|_{\V_{\overline j}})=s_{\overline j}\hs1^2|_{\V_{\overline j}}\,.
\]
Let $s\in\Sg$ be such that the adjoint action of $s$ on $\so|_{\V_{\overline j}}$ coincides with $s_{\overline j}$ and  the adjoint action of $s$ on $\so|_{\V_{\overline k}}$ coincides with $t_{\overline k}s_{\overline k}$. Then
\[
\theta(s.\hs1^2)=s.\hs1^2.
\]
\end{prf}
(By looking at the classification of the Cartan subspaces in \cite{PrzebindaLocal} one can see that in fact the element $s\in\Sg$ in the statement of Theorem \ref{2.2} can be chosen so that $s.\hs1=\hs1$.)
We  assume from now on that $\hs1^2$ is $\theta$-stable. Define
\begin{equation}\label{definition of a}
\a=\{x\in\hs1^2: \theta x=-x\}.
\end{equation}
For each $x\in\a$ the operator $\ad\,x\in \End(\s)$ is symmetric with respect to the form (\ref{3.17}).  
Hence, the Spectral Theorem implies that $\s$ decomposes into an orthogonal direct sum of the eigenspaces for the action of $\a$:
\begin{equation}\label{s decomposition under a}
\s=\sum_\lambda \s_\lambda,\ \text{where}\ \s_\lambda=\{y\in\s:\ [x,y]=\lambda(x)y\ \text{for all}\ x\in\a\}.
\end{equation}
Each  $\s_\lambda$ is $(\Zb/2\Zb)$-graded.
Let $\l=\s_0$ denote the zero eigenspace, i.e. the centralizer of $\a$ in $\s$. 
The non-zero eigenfunctionals  $\lambda\in\a^*$ are called roots (of $\a$ in $\s$). 
Thus for a root $\lambda$ we have a non-zero $y\in\s$ such that 
\begin{equation}\label{3.17'}
[x,y]=\lambda(x)y \qquad (x\in\a).
\end{equation}
Notice also that if $y\in\s$ is a root vector corresponding to $\lambda\in\a^*$, then $\theta y$ is a root vector corresponding to $-\lambda$. Furthermore, $[\s_\lambda,\s_\mu]\subseteq \s_{\lambda+\mu}$.

Let $\reg\a\subseteq\a$ be the subset where no roots vanish.
Fix an element $x_0\in\reg\a$. We say that a root $\lambda$ is positive, written $\lambda>0$, if $\lambda(x_0)>0$. Set
$$
\n=\sum_{\lambda>0} \s_\lambda.
$$ 
Then $\n\subseteq\s$ is a Lie sub-superalgebra and 
$$
\theta\n=\sum_{\lambda>0} \s_{-\lambda}.
$$
Hence (\ref{s decomposition under a}) may be rewritten as
\begin{equation}\label{3.18}
\s=\l\oplus \n\oplus\theta\n,
\end{equation}
Notice that $\s_\lambda\perp\s_\mu$ if $\lambda>0$ and $\mu>0$. Indeed, if $y\in\s_\lambda$ and $z\in\s_\mu$ then
$$
\lambda(x_0)\langle y,z\rangle = \langle [x_0,y],z\rangle = -\langle y, [x_0,z],\rangle= -\mu(x_0) \langle y, z\rangle,
$$
which is possible if and only if $\langle y, z\rangle=0$. 
Thus the restriction of the form $\langle\ \cdot,\cdot \rangle$ to $\n$ is zero. Similarly, we check that $\l\perp\n$ and $\l\perp\theta\n$. 

On the other hand, the Spectral Theorem implies that the restriction of the form (\ref{3.17}) to $\n$ is non-degenerate. 
Hence  the form $\langle\cdot ,\cdot \rangle$  provides a non-degenerate pairing between $\n $ and $\theta\n$. Thus the restriction of the form $\langle\cdot ,\cdot \rangle$ to $\theta\n\oplus \n$ is non-degenerate. Similarly, the restriction of the form $\langle\cdot,\cdot \rangle$ to $\l$ is non-degenerate. Finally we notice that $\n$ is the radical of the form $\langle\cdot ,\cdot \rangle$ restricted to $\l+\n$.

From the structure of the Cartan subspaces  \cite[Theorems 6.2, 5.1 and 4.4]{PrzebindaLocal},  
we see that the following lemma holds.
\begin{lem}\label{3.19}
The Lie superalgebra $\l$ decomposes as
\begin{equation*}
\l=\l_0\oplus\sum_{i=1}^n\l_i\oplus\sum_{i=n+1}^{n+m}\l_i,
\end{equation*}
where $(\l_0)_{\overline 1}=0$ in the type II case and $\l_{0}$ is of the same type as $\s$ in the type I case. For $1\leq i\leq n$, $\l_i$ is isomorphic to $(\g\l_1(\Bbb D), \g\l_1(\Bbb D))$ as a dual pair, where $\Bbb D=\R$, $\C$ or $\Bbb H$. 
If $\Bbb D\ne \R$ then $m=0$. If $\Bbb D=\R$ then, for $n+1\leq i\leq n+m$, $\l_i$ is isomorphic to $(\g\l_2(\Bbb R), \g\l_2(\Bbb R))$ as a dual pair. In all cases $\hs1^2 \cap (\l_0)_{\overline 0}$ (which might be zero) is elliptic and $\hs1^2 \cap \big(\sum_{i=1}^{n+m}(\l_i)_{\overline 0}\big)$ is fundamental.
\end{lem}
\begin{cor}\label{3.20}
With the notation of (\ref{3.18}) we have
$
\dim\,\n_{\overline 0}=\dim\,\n_{\overline 1}.
$
\end{cor}
\begin{prf}
We see from (\ref{3.18}) that it suffices to check that
\begin{equation}\label{3.21}
\dim\,\so - \dim\,\l_{\overline 0} =\dim\,\ss1 - \dim\,\l_{\overline 1}.
\end{equation}
For this we use the classification of the irreducible real reductive dual pairs, see e.g. \cite[p. 548]{HoweTrans}.
Let $\V=\V_{\overline 0}\oplus \V_{\overline 1}$ be the defining module of the pair, and set $d=\dim_\DD \V_{\overline 0}$ and $d'=\dim_\DD \V_{\overline 1}$.
  
Suppose $\s$ is of type II. Then 
\begin{eqnarray*}
\dim\,\so&=&\dim_\R\,\Bbb D \cdot (d^2+d'{}^2),\\
\dim\,\ss1&=&\dim_\R\,\Bbb D \cdot 2dd',\\
\dim\,\l_{\overline 0}&=&\dim_\R\,\Bbb D \cdot (2n+8m+(d-d')^2),\\
\dim\,\l_{\overline 1}&=&\dim_\R\,\Bbb D \cdot (2n+8m),
\end{eqnarray*}
and (\ref{3.21}) follows.

Suppose $\s$ is of type I. Let $\iota=1$ if $\Dc\ne \Bbb H$ and $\iota=\frac{1}{2}$ if $\Dc=\Bbb H$.
We may choose $\epsilon, \epsilon'\in\{0,\frac{\iota}{2}, \iota\}$ such that $\epsilon+\epsilon'=\iota$ and
\begin{eqnarray*}
\dim\,\so&=&\dim_\R\,\Bbb D \cdot \left(\frac{d'(d'-\iota)}{2}+d'\epsilon'
+\frac{d(d-\iota)}{2}+d\epsilon\right),\\
\dim\,\ss1&=&\dim_\R\,\Bbb D \cdot dd',\\
\dim\,\l_{\overline 0}&=&\dim_\R\,\Bbb D \cdot 
\left(\frac{(d'-2n-4m)(d'-2n-4m-1)}{2}+(d'-2n-4m)\epsilon'\right.\\
&+&\left.\frac{(d-2n-4m)(d-2n-4m-1)}{2}+(d-2n-4m)\epsilon
+2n+8m\right),\\
\dim\,\l_{\overline 1}&=&\dim_\R\,\Bbb D \cdot ((d-2n-4m)(d'-2n-4m)+
2n+8m).
\end{eqnarray*}
Hence, (\ref{3.21}) follows.
\end{prf}
By definition, $\l$ and $\n$ are Lie superalgebras and $\l$ normalizes $\n$. Thus,
\begin{equation}\label{3.22}
[\n_{\overline 0}, \l_{\overline 1}]\subseteq \n_{\overline 1},\  
\{\n_{\overline 1}, \l_{\overline 1}\}\subseteq \n_{\overline 0}\ \text{and}\ 
[\n_{\overline 0}, \n_{\overline 1}]\subseteq \n_{\overline 1},.
\end{equation}

Because of Corollary \ref{3.20}, the vector spaces $\n_{\overline 0}$ and $\n_{\overline 1}$ are isomorphic. Let $T:\n_{\overline 1}\to \n_{\overline 0}$ be any fixed isomorphism. Let $A_x:\n_{\overline 0}\to \n_{\overline 1}$ and $B_x:\n_{\overline 1}\to \n_{\overline 0}$
be linear maps depending on a variable $x$. Then $T \circ A_x$ and $B_x \circ T$ are endomorphisms of $\n_{\overline 0}$ and $\n_{\overline 1}$,
respectively. Define
\begin{eqnarray}
\det(A_x)_{\n_{\overline 0}\to \n_{\overline 1}}&=&\det(T \circ A_x) \label{eq:determinant}\\
\det(B_x)_{\n_{\overline 1}\to \n_{\overline 0}}&=&\det(B_x \circ T^{-1})\,. \notag
\end{eqnarray}
Then $\det(B_x \circ A_x)=\det(B_x)_{\n_{\overline 1}\to \n_{\overline 0}} \det(A_x)_{\n_{\overline 0}\to \n_{\overline 1}}$.
Moreover, $B_x \circ A_x$ is an automorphism of $\n_{\overline 0}$ if and only if $A_x$ and $B_x$ are isomorphisms. This is in turn equivalent to the condition that $\det(A_x)_{\n_{\overline 0}\to \n_{\overline 1}}\neq 0$ and
 $\det(B_x)_{\n_{\overline 1}\to \n_{\overline 0}}\neq 0$.
\begin{lem}\label{3.23}
Suppose $x\in \l_{\overline 1}$ is such that $\det(\ad\,x^2)_{\n_{\overline 0}}\ne 0$. Then the map
\begin{eqnarray}\label{3.23.2}
\n_{\overline 0}\ni z\to \exp(z).x-x\in \n_{\overline 1}
\end{eqnarray}
is a well defined analytic diffeomorphism of $\n_{\overline 0}$ onto $\n_{\overline 1}$. Define the maps
\begin{eqnarray}
&&\n_{\overline 0}\ni z\to [z,x]\in \n_{\overline 1}\label{3.23.1}\\
&&\n_{\overline 1}\ni y\to \{y,x\}\in \n_{\overline 0}.\label{3.23.3}
\end{eqnarray}
The composition (\ref{3.23.3})$\circ$(\ref{3.23.1}) coincides with 
\begin{eqnarray}\label{3.23.4}
\n_{\overline 0}\ni z\to [z,x^2]\in \n_{\overline 0}.
\end{eqnarray}
If we identify the vector spaces $\n_{\overline 0}$ and $\n_{\overline 1}$ by fixing a linear isomorphism between them, as above (\ref{eq:determinant}), then the Jacobian of the map (\ref{3.23.2}) at $z$ is equal to $|\det(\ad\,x)_{\n_{\overline 0}\to\n_{\overline 1}}|$, the absolute value of the determinant of the map (\ref{3.23.1}). Also, 
\begin{equation}\label{3.23.5}
|\det(\ad\,x^2)_{\n_{\overline 0}}|=|\det(\ad\,x)_{\n_{\overline 1}\to\n_{\overline 0}}\,\det(\ad\,x)_{\n_{\overline 0}\to\n_{\overline 1}}|.
\end{equation}
\end{lem}
\begin{prf}
We see from (\ref{3.22}) that the maps (\ref{3.23.1}) and (\ref{3.23.3}) are well defined. Since the adjoint action is a derivation and since $x^2=\frac{1}{2}\{x,x\}$, the composition (\ref{3.23.3})$\circ$(\ref{3.23.1}) coincides with (\ref{3.23.4}). Hence, (\ref{3.23.5}) follows.

Also, we see from (\ref{3.22}) that for $z \in \n_{\overline 0}$ the element
\begin{eqnarray*}
\exp(z).x-x=[z,x]+\frac{1}{2}[z,[z,x]]+\frac{1}{3!}[z,[z,[z,x]]]+ \dots
\end{eqnarray*}
belongs to $\n_{\overline 1}$.  Thus the map (\ref{3.23.2}) is well defined. 

Let $z,y\in \n_{\overline 0}$. As is well known, see e.g. \cite[\S1.2, p. 15]{Rossmannbook}, 
\begin{eqnarray*}
\exp(-z)\exp(z+y)=I+\frac{I-\exp(-\ad\,z)}{\ad\,z}y+\text{higher order terms in}\ y.
\end{eqnarray*}
Hence, the derivative of the map (\ref{3.23.2}) at $z$ coincides with the following map
\begin{eqnarray}\label{3.24}
\n_{\overline 0}\ni y\to \exp(z).\left[\left(\frac{I-\exp(-\ad\,z)}{\ad\,z}y\right),x\right]\in \n_{\overline 1}. 
\end{eqnarray}
Since the map (\ref{3.23.1}) is injective and since $z$ is nilpotent, it is easy to see that the map (\ref{3.24}) is injective and in fact has the determinant equal to $\pm \det(\ad\,x)_{\n_{\overline 0}\to\n_{\overline 1}}$. 
Hence, the map (\ref{3.23.2}) is everywhere locally bijective. 

Suppose $\exp(z).x=\exp(z').x$. Then $\exp(z).x^2=(\exp(z).x)^2=(\exp(z').x)^2
=\exp(z').x^2$. Hence, by \cite[page 218]{HC-57Fourier}, $z=z'$. Therefore the map (\ref{3.23.2}) is injective. 

We already know that the image of the map (\ref{3.23.2}) is open. If we show that this image is also closed then the 
surjectivity will follow. As in \cite[page 218]{HC-57Fourier}, it suffices to show that if $\exp(z).x-x$ is bounded then $z$ is bounded. But if $\exp(z).x$ is bounded, then $\exp(z).x^2=(\exp(z).x)^2$ is bounded. In turn, this implies that $z$ is bounded, see \cite[page 218]{HC-57Fourier}.
\end{prf}

\begin{lem}\label{nilpotent elements of l+n}
Any element of $\n_{\overline 1}$ is nilpotent.
Any nilpotent element $z\in\l_{\overline 1}+\n_{\overline 1}$ is of the form $z=z_\l+z_\n$, where both $z_\l\in\l_{\overline 1}$ and $z_\n\in \n_{\overline 1}$ are nilpotent. 
\end{lem}
\begin{prf}

Since $\n_{\overline 1}\subseteq \sum_{\lambda>0} \s_\lambda$ and $[\s_\lambda,\s_\mu]\subseteq \s_{\lambda+\mu}$, the endomorphism $\ad x$ is nilpotent for every $x \in \n_{\overline 1}$. 
By Lemma \ref{lemma:ad-nilpotent}, $x$ is nilpotent. 

Suppose now that $z=z_\l+z_\n$ is nilpotent, and $z_\l\in\l_{\overline 1}$ and $z_\n\in \n_{\overline 1}$. 
By Lemma \ref{lemma:ad-nilpotent}, $\ad z$ is a nilpotent endomorphism of $\End(\s)$. Let $k \in \Bbb N$ so that $(\ad z)^k=0$.
Let $\Delta$ denote the set of positive roots. Then $z \in \l+\n= \sum_{\lambda \in \Delta \cup \{0\}} \s_\lambda$ and 
$z_\n\in \n=\sum_{\lambda \in \Delta} \s_\lambda$. Let $\mu$ be either a root or $\mu=0$. For $a, b \in {\Bbb N}$ we have:
\begin{eqnarray*}
(\ad z_\n)^a(\s_\mu)&\subseteq & \sum_{\lambda \in \mu +a \Delta} \s_\lambda\,,\\
(\ad z)^b(\s_\mu)&\subseteq & \sum_{\lambda \in \mu + \Delta} \s_\lambda\,,\\
(\ad z)^b(\ad z_\n)^a(\s_\mu)&\subseteq & \sum_{\lambda \in \mu + \Delta+ a \Delta} \s_\lambda\,,\\
\end{eqnarray*}
where we have set $a\Delta=\{ \lambda_1+\dots+\lambda_a: \lambda_j \in \Delta\}$ and $\s_\lambda=0$ if $\lambda\neq 0$ is not a root.   
Since the set of roots is finite, there is $h \in \Bbb N$  such that for every root $\mu $ and all integers $a \geq h$ the set
 $\mu + a \Delta$ does not contain any root. 
For $p \in \Bbb N$ we have
\begin{equation*}
(\ad z_\l)^p=(\ad z-\ad z_\n)^p=\sum \pm (\ad z)^{b_1}(\ad z_\n)^{a_1}\dots (\ad z)^{b_s}(\ad z_\n)^{a_s},
\end{equation*}
where the summation is over positive integers $a_j$, $b_j$, $1\leq i\leq s$, with $\sum_{j=1}^s(a_j+b_j)=p$.
Fix $\mu$ and $x\in \s_\mu$, and look at 
$$z=(\ad z)^{b_1}(\ad z_\n)^{a_1}\dots (\ad z)^{b_s}(\ad z_\n)^{a_s}(x)\,.$$
Then $z \in \sum_{\lambda \in \mu + \Delta+ a\Delta} \s_\mu$ where $a=a_1+\dots+a_s$.
If $s \geq h$ then $a \geq h$. So $\s_\lambda =0$ for all $\lambda\in \mu + \Delta+ a\Delta$, and thus $z=0$. 
Suppose then $s<h$. If $a\geq h$, then, as above, $z=0$. We can therefore assume $a<h$. Set $b=b_1+\dots+b_s$.
Then $b=p-a>p-h$. If we have chosen $p$ sufficiently big (specifically, $\geq h(k+1)$), then there must be an index $j$ so that $b_j\geq k$. Otherwise $b< sk$, which would give $hk \geq sk >p-h$. But if $b_j\geq k$ then $(\ad z)^{b_j}=0$. We conclude that $z=0$ also in this case.

This proves that $\ad z_\l$ is nilpotent. Because of Lemma \ref{lemma:ad-nilpotent}, $z_\l$ must be nilpotent too.
\end{prf}

\begin{cor}\label{density for gl1 gl2}
Consider a nilpotent element $z\in\ss1$ with the decomposition $z=z_\l+z_\n$ as in Lemma \ref{nilpotent elements of l+n}. Suppose $\l'_{\overline 1}\subseteq\l_{\overline 1}$ is a subset containing $z_\l$ in its closure and such that  $\det(\ad\,x^2)_{\n_{\overline 0}}\ne 0$ for all $x\in \l'_{\overline 1}$.
Then $z$ may be approximated by elements of the $\Sg$-orbits passing through $\l'_{\overline 1}$.
\end{cor}
\begin{prf} 
Choose a small element $y\in \l_{\overline 1}$ such that $z_\l-y\in\l'_{\overline 1}$. Let $x=z_\l-y$. 
We see from Lemma \ref{3.23} that there is an element $u\in \n_{\overline 0}$ such that
$$
z_\n=\exp(u).x-x.
$$
Hence,
$$
(z_\n+z_\l)-\exp(u).x=y
$$
is small.
\end{prf}
Since our proof of the following proposition requires a significant amount of additional notation, we shall present it in Appendix B. (In fact what we are missing here is an analog of the Jacobson-Morozov Theorem.)
\begin{pro}\label{good nilpotent elements}
For any nilpotent element  $z\in \ss1$ there is a Cartan subspace $\hs1\subseteq\ss1$ such that, in terms of (\ref{3.18}), $z$ belongs to an $\Sg$-orbit passing through $\l_{\overline 1}+\n_{\overline 1}$. Moreover, $\l_{\overline 1}\supseteq \t\h_{\overline 1}$.

Suppose $z\in\l_{\overline 1}$ is nilpotent. Let $\L\subseteq \Sg$ be the subgroup with the Lie algebra $\l_{\overline 0}$, so that $(\L,\l)$ is the direct product of classical Lie supergroups, as in \cite[sec.2]{PrzebindaLocal}.
Then $z_\l$ is the limit of elements in the $\L$-orbits passing through $\t\h_{\overline 1}$.

\end{pro}

\begin{thm} \label{thm:nilpotent-almostSS}
Every nilpotent element $z \in \ss1$ is a limit of almost semisimple elements. It is a limit of semisimple elements if and only if $(\Sg,\s)$ is not associated with (\ref{basic assumption-one}) or (\ref{basic assumption-two}).
\end{thm}
\begin{prf}
By Proposition \ref{good nilpotent elements}, $z$ singles out a Cartan subspace $\hs1\subseteq \ss1$. Moreover, if $\s=\l\oplus \n\oplus \theta\n$ is the corresponding decomposition (\ref{3.18}), then $z\in\Sg. \wt{z}$ with $\wt{z}=\wt{z}_\l+\wt{z}_\n \in \l_{\overline 1}\oplus \n_{\overline 1}$ necessarily nilpotent. Lemma  \ref{nilpotent elements of l+n} shows that $\wt{z}_\l \in \l_{\overline 1}$ is nilpotent. So, by Proposition 
\ref{good nilpotent elements}, there are elements $l_n \in L$ and $x_n\in \t\h_{\overline 1}$ so that $\wt{z}_\l=\lim_{n\to \infty} l_n.x_n$. 
By Proposition \ref{good nilpotent elements}, $\t\h_{\overline 1} \subset \l_{\overline 1}$.
The map $x \to \det(\ad x^2)_{\n_{\overline 0}}$ is polynomial on $\l_1$. Since it is not identically zero on $\hs1$, the set 
$\t\h_{\overline 1}'=\{x \in \t\h_{\overline 1}:\det(\ad x^2)_{\n_{\overline 0}}\neq 0\}$ is dense in $\t\h_{\overline 1}$. So we can suppose that $x_n \in \t\h_{\overline 1}'$. Corollary \ref{density for gl1 gl2} with $\l_{\overline 1}'=\t\h_{\overline 1}'$ therefore assures that $\wt z$ is a limit of elements of the $S$-orbits passing through $\t\h_{\overline 1}'$. The same property holds thus for $z$.  
\end{prf}

Combining Proposition \ref{pro:summaryLocal}, Theorem \ref{thm:nilpotent-almostSS} with Jordan decomposition of an arbitrary element 
of $\ss1$ into the sum of its semisimple and nilpotent components, we obtain the following theorem.

\begin{thm}\label{density}
The set of the almost semisimple elements of $\ss1$ is dense in $\ss1$. 
The set of the semisimple elements in $\ss1$ is dense in $\ss1$ if and only if $(\Sg,\s)$ is not associated with (\ref{basic assumption-one}) or (\ref{basic assumption-two}).
\end{thm}
\begin{prf}
Let $x=x_s+x_n$ be the Jordan decomposition of an element in $x\in \ss1$. 
Because of Proposition \ref{pro:summaryLocal}(c) and Theorem \ref{thm:nilpotent-almostSS}, we can assume that $x$ is neither semisimple nor nilpotent. In particular $x_s\neq 0$. Since $x$, $x_s$ and $x_n$ commute in $\End(\V)$, they all belong to $\ss1^{x_s^2}$. As in the proof of Theorem 4.4(b) in \cite[section 13]{PrzebindaLocal}, let 
$\V=\V^0\oplus \V^1 \oplus \cdots \oplus \V^k$ be the isotypic decomposition of $\V$ with respect to $x_s$. Hence $\V^0=\ker(x_s)$ and 
$\V^1 \oplus \cdots \oplus \V^k=\V^+=x_s\V$. Isotypic means that each $(x_s|_{\V^j},\V^j)$ decomposes into mutually similar indecomposable pieces.
Let $(\Sg(\V^j),\s(\V^j))$ denote the restriction of $(\Sg,\s)$ to $\V^j$. Then 
\begin{equation}\label{eq:isotypic}
\ss1^{x_s^2}=\ss1(\V^0) \oplus \ss1(\V^1)^{x_s^2} \oplus \cdots \oplus \ss1(\V^k)^{x_s^2}\,.
\end{equation}
Decompose $x=x_0+x_1+\dots+x_k$ according to (\ref{eq:isotypic}). If we prove that for $j=1,\dots,k$ there is a sequence $y_{m,j}$ of semisimple elements so that $x_j=\lim_{m\to \infty} y_{m,j}$, then $x=\lim_{m\to \infty} x_0+y_{m,1}+\dots+y_{m,k}$ will be a limit of almost semisimple
elements. Notice that $\ker(x_s|_{\V^j})=0$. Hence, by \cite[Theorem 4.4(a)]{PrzebindaLocal}, each $\ss1(\V^1)^{x_s^2}$ is a supergroup with corresponding dual pair isomorphic either to $(\Ug_n,\Ug_n)$ or to $(\GL_n(\DD),\GL_n(\DD))$. Here the division algebra $\DD$ may be different from the one over which $\V$ is defined. This reduces the proof of our Theorem to the case where $x=x_s+x_n$ belongs to one of these isotypic components. A further reduction allows us to select $x_s$ inside non-conjugate Cartan subalgebras $\hs1\subseteq \ss1$. Indeed suppose $x_s=g.x_s'$ with $g \in \Sg$ and $x_s' \in \hs1$. Write $x_n=g.x_n'$. Then $x_n'$ is nilpotent and $x_n'x_s'=x_s'x_n'$. If $x_s'+x_n'=\lim y_m$ with $y_m$ semisimple, then $x=\lim g.y_m$ with $g.y_m$ semisimple. This case-by-case analysis will be carried out in Appendix C.
\end{prf}

Recall that the almost semisimple elements in $\ss1$ are the elements of the $\Sg$-orbits of the sets 
$$\t{\h}_{\overline 1}=\hs1\oplus \ss1(\V^0)$$
when $\hs1$ ranges among the (conjugacy classes of) Cartan subspaces in $\hs1$. Here $\V^0$ denotes the intersection of the 
kernel of all elements in $\hs1$. Recall also that if $x\in \hs1$ is regular, then $\V^0=\ker(x)$ and we have the decomposition
$\V=\V^0 \oplus \V^+$ with $\V^+=x\V$. 
Let  $\Sg^{\hs1}$ denote the centralizer of $\hs1$ in $\Sg$. Then 
\begin{equation}\label{eq:sghs}
\Sg^{\hs1}=\Sg(\V^0)\times \Sg(\V^+)^{\hs1}\,.
\end{equation} 
The group $\Sg^{\hs1}$ is reductive and hence unimodular. Therefore there is an invariant measure $d(g\Sg^{\hs1})$ on the quotient
space $\Sg/\Sg^{\hs1}$.

The title of this section refers to the following theorem.
\begin{thm}\label{Weyl integration}
For a continuous compactly supported function $\phi:\ss1\to\C$ the following integration formula holds,
\begin{eqnarray*}\label{Weyl int}
&&\int_{\ss1}\phi(x)\,dx=\sum_{\hs1}\frac{1}{|W(\Sg,\h_{\overline 1})|}\int_{\hs1^2}
|\pi_{\so/\hs1^2}(x^2)|\int_{\Sg/\Sg^{\hs1}}\int_{\ss1(\V^0)} \phi(g.(y+x))\,dy\,d(g\Sg^{\hs1}) 
\,d x^2,
\end{eqnarray*}
where the summation is over a maximal family of mutually non-conjugate Cartan subspaces $\hs1\subseteq \ss1$.
\end{thm}
For a fixed $x\in\h_{\overline 1}$ the integral $\int_{\ss1(\V^0)} \phi(g.(y+x))\,dy$ is a function of $g\in\Sg$, invariant under the right translations by elements of $\Sg^{\hs1}$.
Also, the integral over $\Sg/\Sg^{\hs1}$ is constant on the fibers of the square map $x\mapsto x^2$, because the group $W(\Sg,\hs1)$ acts transitively on these fibers. (See \cite[Corollary 6.25]{PrzebindaLocal}.)
\begin{prf}
Notice that for any $y_0\in \ss1(\V^0)\setminus\{0\}$ the $\Sg(\V^0)$-orbit through $y_0$ is equal to $\ss1(\V^0)\setminus\{0\}$. 
Fix any $y_0\in \ss1(\V^0)\setminus\{0\}$, then
$\Sg^{\hs1+y_0}$, the centralizer of $y_0+\hs1$ in $\Sg$, is a unimodular group. 
Furthermore, by (\ref{eq:sghs}), $\Sg^{y_0+\hs1}=\Sg(\V^0)^{y_0}\times \Sg(\V^+)^{\hs1}$.
Hence
\begin{equation} \label{orbital integral}
\int_{\Sg/\Sg^{\hs1}}\int_{\ss1(\V^0)} \phi(g.(y+x))\,dy\,d(g\Sg^{\hs1}) 
=\int_{\Sg/\Sg^{y_0+\hs1}}\phi(g.(y_0+x))\,d(g\Sg^{y_0+\hs1}), 
\end{equation}
Also, if the function $\phi$ is supported in the union of the orbits passing through $\t\h_{\overline 1}=\ss1(\V^0)\times \h_{\overline 1}$, then
\begin{equation*}
\int_{\ss1}\phi(x)\,dx=\frac{1}{|W(\Sg,\h_{\overline 1})|}\int_{\t\h_{\overline 1}}\int_{\Sg/\Sg^{y_0+\hs1}}\phi(g.(y_0+x))\,d(g\Sg^{y_0+\hs1})\,\J(y_0+x)\,dy_0\,dx, 
\end{equation*}
where $\J$ is the Jacobian, which may be computed as in the proof of \cite[Corollary 6.21]{PrzebindaLocal}. In particular, $\J(y_0+x)$ is equal to $|\pi_{\so/\hs1^2}(x^2)|$ times the Jacobian of the square map $\h_{\overline 1}\ni x\mapsto x^2\in \h_{\overline 1}^2$. Since, as we mentioned previously, the integral is constant on the fibers of this map, the formula follows.
\end{prf}

\section{A semisimple orbital integral on the odd part of an ordinary classical Lie superalgebra}
\label{section:orbital symplectic}

In this section we prove that the integration formula of Theorem \ref{Weyl integration} extends to Schwartz functions on the symplectic space
$\W=\ss1$. For this, we study whether the orbital integral of Theorem \ref{Weyl integration} is integrable on $\hs1^2$ when 
$\phi(x)=(1+|x|)^{-N}$ and $N>0$ is sufficiently large. Here $|\cdot |$ denotes a fixed norm on the real vector space $\s$.  Recall that
the Cartan subspace $\hs1\subseteq\ss1$ is called elliptic if and only if all the eigenvalues of the $\hs1^2$ acting on $\so$ are imaginary.
\begin{pro}\label{3.3}
Suppose the Cartan subspace $\hs1\subseteq\ss1$ is elliptic. Then for any constant $M\geq 0$ there are positive constants $C, N$ such that
\begin{equation*}
|\pi_{\so/\hs1^2}(x^2)|\int_{\Sg/\Sg^{\hs1}}\int_{\ss1(\V^0)} (1+|g.(x+y)|)^{-N}\,dy\,d(g\Sg^{\hs1})  \leq C(1+|x|)^{-M} \qquad (x\in\reg{\hs1}).
\end{equation*}
Here the integral over $\ss1(\V^0)$ is omitted if the dual reductive pair is not isomorphic to (\ref{basic assumption-one}) or (\ref{basic assumption-two}).
\end{pro}
\begin{prf}
We may assume that the norm $|\cdot |$ is chosen so that the maximum of $|x^2|$, when $x\in\ss1$ varies through the unit sphere, is $1$. Then
\begin{equation}\label{3.3.1}
|x^2|\leq |x|^2 \qquad (x\in\ss1). 
\end{equation}
Also, it is clear from the description of the Cartan subspaces in \cite{PrzebindaLocal} that there is a constant $1\leq C<\infty$ such that
\begin{equation}\label{3.3.2}
|x|^2\leq C |x^2| \qquad (x\in\hs1). 
\end{equation}
Thus for any non-negative constants $m$, $N$ and $g\in\Sg$,
\begin{equation*}
(1+|x|^2)^M(1+|g.(x+y)|^2)^{-N}\leq C^M(1+|x^2|)^M (1+|g.(x^2+y^2)|)^{-N} \quad (x\in\hs1,\ y\in\ss1(\V^0)). 
\end{equation*}
Notice that since the Cartan subspace $\hs1\subseteq\ss1$ is elliptic, the quotient $(\Sg/\Sg^{\hs1})/(\Sg/\Sg^{\hs1^2})$ is compact and therefore has finite volume. Hence our integral is dominated by the identical integral with the $x$ replaced by $x^2$ and the $\Sg/\Sg^{\hs1}$ replaced by $\Sg/\Sg^{\hs1^2}$.

Now we apply Theorem \ref{2.2} or Proposition \ref{2.2pro} and (\ref{2.4pro}) with $\G=\Sg$, $\g=\so$, $z=x^2$ for any $x\in\reg{\hs1}$ so that $\Zg=\Sg^{\hs1^2}$, $\z=\so^{\hs1^2}$, $\hs1^2\subseteq \c$ and $\pi_{\g/\z}|_{\hs1^2}=\pi_{\so/\hs1^2}$. 
\end{prf}

\begin{rem}
The inequalities (\ref{3.3.1}) and (\ref{3.3.2}) show that the estimate of Proposition \ref{3.3} is as sharp as the estimate of Theorem \ref{2.2}. In case when $\Zg=\H\subseteq \G$ is a Cartan subgroup, there is a characterization of the regular semisimple orbital integrals on a semisimple Lie algebra due to Bouaziz, \cite{BouazizOrb}, from which one can see that the estimate of Theorem \ref{2.2} is sharp.
\end{rem}

For $N\geq 0$ define
\begin{equation}\label{3.4}
I_N(t)=
\begin{cases}
&t^{-N}\ \ \text{if}\ \ t\geq 1,\\
&1-N\,\ln(t)\ \ \text{if}\ \ 0< t\leq 1
\end{cases}
\end{equation}
Then
\begin{equation}\label{A.2}
(I_N(t))^{1/2}\leq I_{N/2}(t)\qquad (t>0)
\end{equation}
and
\begin{equation}\label{A.3}
I_N(s)\geq I_{N}(t)\qquad (0<s<t).
\end{equation}
\begin{lem}\label{3.5}
Suppose that the supergroup $(\Sg,\s)$ is such that the corresponding dual pair is isomorphic to $(\GL_1(\Bbb D), \GL_1(\Bbb D))$, $\Bbb D=\R, \C$ or $\Bbb H$. Then for each $N\geq 1$ there is a non-negative, finite constant $C$ such that
\begin{equation*}
\int_{\Sg/\Sg^{\hs1}} (1+|s.x|^2)^{-N}\,d(s\Sg^{\hs1})  \leq CI_N(|x|) \qquad (x\in\reg{\hs1}).
\end{equation*}
\end{lem}
\begin{prf}
We may realize our supergroup in terms of matrices as follows,
\begin{equation*}
\Sg=\left\{
s=\left(
\begin{array}{cc}
g & 0\\
0 & g'
\end{array}
\right);\ g,g'\in \Bbb D^\times
\right\}
,\ \ 
\ss1=\left\{
\left(
\begin{array}{cc}
0 & z_1\\
z_2 & 0
\end{array}
\right);\ z_1, z_2\in \Bbb D
\right\}.
\end{equation*}
Then, by the classification of Cartan subspaces, \cite{PrzebindaLocal}, we may assume that
\begin{equation*}
\hs1=\left\{
x=\left(
\begin{array}{cc}
0 & z\\
\epsilon z & 0
\end{array}
\right);\ z\in \wt{\Bbb D}
\right\},
\end{equation*}
where $\wt{\Bbb R}=\R$,  $\wt{\Bbb C}=\C$ and $\wt{\Bbb H}$ is a copy of $\C$ in $\Bbb H$. Moreover, $\epsilon=\pm 1$ if $\Bbb D=\R$ and $\epsilon=1$ in $\Bbb D\ne \R$. Since,
\begin{equation*}
\left(
\begin{array}{cc}
g & 0\\
0 & g'
\end{array}
\right)
\left(
\begin{array}{cc}
0 & z\\
\epsilon z & 0
\end{array}
\right)
\left(
\begin{array}{cc}
g & 0\\
0 & g'
\end{array}
\right)^{-1}=
\left(
\begin{array}{cc}
0 & gzg'{}^{-1}\\
g'\epsilon zg^{-1} & 0
\end{array}
\right)
\end{equation*}
we see that
\begin{equation*}
\Sg^{\hs1}
=\left\{
\left(
\begin{array}{cc}
g & 0\\
0 & g
\end{array}
\right);\ 
g\in\wt{\Bbb D}
\right\}.
\end{equation*}
The following formula defines a norm $|\cdot |$ on $\ss1$:
\begin{equation*}
\left|
\left(
\begin{array}{cc}
0 & z_1\\
z_2 & 0
\end{array}
\right)
\right|^2= |z_1|^2+|z_2|^2,
\end{equation*}
where $|z|^2=z\overline z$, $z\in\wt{\Bbb D}$. Let
\begin{equation*}
U(\Bbb D)=\{g\in\Bbb D^\times;\ |g|=1\},\ U(\wt{\Bbb D})=\wt{\Bbb D}\cap U(\Bbb D).
\end{equation*}
Then, as a homogeneous space,
\begin{equation*}
\Sg/\Sg^{\hs1}=\Bbb D^\times\times(\Bbb D^\times/\wt{\Bbb D}^\times)=\Bbb D^\times\times (U(\Bbb D)/U(\wt{\Bbb D})).
\end{equation*}
Hence, we may normalize all the measures involved so that
\begin{eqnarray*}
&&\int_{\Sg/\Sg^{\hs1}} (1+|s.x|^2)^{-N}\,d(s\Sg^{\hs1})\\
&=&\int_{\Bbb D^\times}\int_{U(\Bbb D)}(1+|gzg'{}^{-1}|^2+|g'zg^{-1}|^2)^{-N}\,dg'\,dg\\
&=&\int_0^\infty(1+(a^2+a^{-2})|z|^2)^{-N}\,da/a\leq \frac{2}{N}I_N(|z|), 
\end{eqnarray*}
where the last inequality follows from (\ref{A.4.a}).
\end{prf}
\begin{lem}\label{3.7}
Suppose that the supergroup $(\Sg,\s)$ is such that the corresponding dual pair is isomorphic to $(\GL_2(\R), \GL_2(\R))$ and that the centralizer of $\hs1^2$ in $\so$ is a fundamental Cartan subalgebra. Then for each $N>8$ there is a finite constant $C$ such that
\begin{equation*}
|\pi_{\so/\hs1^2}(x^2)|\int_{\Sg/\Sg^{\hs1}} (1+|s.x|^2)^{-N}\,d(s\Sg^{\hs1})  \leq C(1+|x|)^{-\frac{N}{2}} \qquad (x\in\reg{\hs1}).
\end{equation*}
\end{lem}
\begin{prf}
We realize our supergroup in terms of matrices as follows,
\begin{equation*}
\Sg=\left\{
s=\left(
\begin{array}{cc}
g & 0\\
0 & g'
\end{array}
\right);\ g,g'\in \GL_2(\R)
\right\}
,\ \ 
\ss1=\left\{
\left(
\begin{array}{cc}
0 & z_1\\
z_2 & 0
\end{array}
\right);\ z_1, z_2\in \g\l_2(\R)
\right\}.
\end{equation*}
Then, by the classification of Cartan subspaces, \cite{PrzebindaLocal}, we may assume that
\begin{equation*}
\hs1=\left\{
x=\left(
\begin{array}{cc}
0 & z\\
 z & 0
\end{array}
\right);\ 
z=\left(
\begin{array}{cc}
\xi & \eta\\
-\eta & \xi
\end{array}
\right),\ 
\xi, \eta\in \R
\right\}.
\end{equation*}
For $x\in\hs1$ as above, we shall write $\xi=\Re(z)$ and $\eta=\Im(z)$. As in the proof of Lemma \ref{3.5} we check that
\begin{equation*}
\Sg^{\hs1}
=\left\{
\left(
\begin{array}{cc}
g & 0\\
0 & g
\end{array}
\right);\ 
g=
\left(
\begin{array}{cc}
\xi & \eta\\
-\eta & \xi
\end{array}
\right)\in \GL_2(\R)
\right\},
\end{equation*}
which is isomorphic to $\C^\times$ embedded into $\GL_2(\R)$, as a Lie group.
Hence, as a homogeneous space,
\begin{equation*}
\Sg/\Sg^{\hs1}=\GL_2(\R)\times(\GL_2(\R)/\C^\times)=\GL_2(\R)\times(\SL_2^\pm(\R)/\SO_2).
\end{equation*}
Define a norm $|\cdot |$ on $\ss1$ by
\begin{equation*}
\left|
\left(
\begin{array}{cc}
0 & z_1\\
z_2 & 0
\end{array}
\right)
\right|^2= tr(z_1^tz_1)+tr(z_2^tz_2).
\end{equation*}
we may normalize all the measures involved so that
\begin{eqnarray}\label{3.10}
&&\int_{\Sg/\Sg^{\hs1}} (1+|s.x|^2)^{-N}\,d(s\Sg^{\hs1})\\
&=&\int_{\GL_2(\R)}\int_{\SL_2^\pm(\R)}(1+|gzg'{}^{-1}|^2+|g'zg^{-1}|^2)^{-N}\,dg'\,dg\nn\\
&=&\int_{\SL_2(\R)}\int_{\SL_2(\R)}\int_0^\infty(1+a^2|gzg'{}^{-1}|^2
+a^{-2}|g'zg^{-1}|^2)^{-N}\,da/a\,dg'\,dg\nn\\ 
&\leq&\int_{\SL_2(\R)}\int_{\SL_2(\R)}\frac{2}{N} I_N((|gzg{}^{-1}|
|g'zg^{-1}|)^{1/2})\,dg'\,dg,\nn
\end{eqnarray}
where the inequality follows from (\ref{A.4.b}). Notice that
\begin{equation*}
|gzg'{}^{-1}||g'zg^{-1}|\geq |gzg'{}^{-1}g'zg^{-1}|=|gz^2g^{-1}|.
\end{equation*}
Hence, by (\ref{A.3}) and (\ref{A.2}),
\begin{eqnarray*}
&&I_N((|gzg'{}^{-1}||g'zg^{-1}|)^{1/2})\leq (I_N((|gz^2g{}^{-1}|^{1/2})
I_N((|g'z^2g'{}^{-1}|^{1/2}))^{1/2}\\
&\leq&I_{N/2}((|gz^2g{}^{-1}|^{1/2})I_{N/2}((|g'z^2g'{}^{-1}|^{1/2}).
\end{eqnarray*}
Therefore (\ref{3.10}) is less or equal to
\begin{equation}\label{3.11}
\frac{2}{N}\left(
\int_{\SL_2(\R)}I_{N/2}((|gz^2g^{-1}|^{1/2})\,dg\right)^2.
\end{equation}
We perform the integration (\ref{3.11}) in terms of Cartan decomposition ($\G=\K\A\K$):
\begin{eqnarray}\label{3.12}
&&\int_{\SL_2(\R)}I_{N/2}(|gz^2g^{-1}|^{1/2})\,dg\\
&=&\int_1^\infty I_{N/2}((2( \Re (z^2))^2+(a^4+a^{-4})(\Im(z^2))^2)^{1/4})\frac{1}{2}(a^2-a^{-2})\,da/a\nn\\
&\leq&\int_1^\infty I_{N/2}(a|\Im(z^2)|^{1/2})a^2\,da/a
\leq \int_0^\infty I_{N/2}(a|\Im(z^2)|^{1/2})a^2\,da/a\nn\\
&=&|\Im(z^2)|^{-1}\int_0^\infty I_{N/2}(a)a\,da,\nn
\end{eqnarray}
where that last integral is finite because $\frac{N}{2}>2$. Therefore (\ref{3.10}) is less or equal to 
\begin{equation}\label{3.13}
\frac{2}{N}\left(\int_0^\infty I_{N/2}(a)a\,da\right)^2|\Im(z^2)|^{-2}.
\end{equation}
Notice that, by (\ref{A.2}), in terms of (\ref{3.12}) we have
\begin{eqnarray*}
&& I_{N/2}((2(\Re(z^2))^2+(a^4+a^{-4})(\Im(z^2))^2)^{1/4})\\
&\leq&\left(I_{N/4}((2(\Re(z^2))^2+(a^4+a^{-4})(\Im(z^2))^2)^{1/4})\right)^2\\
&\leq&I_{N/4}((|z^2|^{1/2})I_{N/4}((2(\Re(z^2))^2+(a^4+a^{-4})(\Im(z^2))^2)^{1/4})\\
&\leq&I_{N/4}((|z^2|^{1/2})I_{N/4}(a|\Im(z^2)|^{1/2}).
\end{eqnarray*}
Hence, (\ref{3.12}) is less or equal to 
\begin{equation}\label{3.14}
|\Im(z^2)|^{-1}I_{N/4}(|z^2|^{1/2})\int_0^\infty I_{N/4}(a)a\,da,
\end{equation}
where the integral is finite because $N/4>2$. Therefore (\ref{3.10}) is less or equal to
\begin{equation}\label{3.15}
\frac{2}{N}\left(I_{N/4}(|z^2|^{1/2})\right)^2\left(\int_0^\infty I_{N/4}(a)a\,da\right)^2|\Im(z^2)|^{-2}.
\end{equation}
Clearly (\ref{3.10}) is less or equal to the minimum of (\ref{3.13}) and 
(\ref{3.15}). Thus there is a finite constant $C$ such that (\ref{3.10}) is less or equal to
\begin{equation}\label{3.16}
C(1+|z|)^{-N/2}|\Im(z^2)|^{-2}.
\end{equation}
Since $|\pi_{\so/\hs1^2}(x^2)|$ is proportional to $|\Im(z^2)|^{2}$, we see that (\ref{3.16}) completes the proof.
\end{prf}
Let $\L,\N\subseteq \Sg$ be the Lie subgroups with the Lie algebras $\l_{\overline 0}$, $\n_{\overline 0}$ respectively and let $\K\subseteq \Sg$ be the maximal compact subgroup corresponding to the Cartan involution 
$\theta$. 
\begin{cor}\label{3.25}
For any non-negative measurable function $\phi:\ss1\to\R$ and $x\in \hs1^{reg}$ 
\begin{eqnarray*}
&&\left|\pi_{\so/\hs1^2}(x^2)\right|
\int_{\Sg/\Sg^{\hs1}} \phi(s.x)\,d(s\Sg^{\hs1})\\
&=&\left|\pi_{\l_{\overline 0}/\hs1^2}(x^2)\right|\left|\det(\ad\,x)_{\n_{\overline 1}
\to\n_{\overline 0}}\right|\int_{\L/\L^{\hs1}} \int_{\n_{\overline 1}}\int_\K 
\phi(k.(l.x+y)))\,dk\,dy\,d(l\L^{\hs1})\\
&=&\left|\pi_{\l_{\overline 0}/\hs1^2}(x^2)
\right|
\left|\det(\ad\,x^2)_{\n_{\overline 0}}\right|^{\frac{1}{2}}
\int_{\L/\L^{\hs1}} \int_{\n_{\overline 1}}\int_\K 
\phi(k.(l.x+y)))\,dk\,dy\,d(l\L^{\hs1}).
\end{eqnarray*}
\end{cor}
\begin{prf}
Notice that $\Sg=\K\N\L$ and that the Haar measure on $\Sg$ may be written as $ds=dk\,dn\,dl$. Furthermore, Lemma \ref{3.23} implies that 
\begin{equation*}
\int_\N\phi(n.x)\,dn=
\frac{1}{|\det(\ad\,x)_{\n_{\overline 0}\to\n_{\overline 1}}|}
\int_{\n_{\overline 1}}\phi(x+y)\,dy.
\end{equation*}
We see from (\ref{3.23.5}) that
\begin{eqnarray*}
&&\frac{|\pi_{\so/\hs1^2}(x^2)|}{|\det(\ad\,x)_{\n_{\overline 0}\to\n_{\overline 1}}|}
=\frac{|\pi_{\l_{\overline 0}/\hs1^2}(x^2)||\det(\ad\,x^2)_{\n_{\overline 0}}|}{|\det(\ad\,x)_{\n_{\overline 0}\to\n_{\overline 1}}|}
=|\pi_{\l_{\overline 0}/\hs1^2}(x^2)\,\det(\ad\,x)_{\n_{\overline 1}\to\n_{\overline 0}}|.
\end{eqnarray*}
This implies the first equality.
Since the maps (\ref{3.23.1}), (\ref{3.23.3}) are adjoint to each other (see \cite[Proposition 2.18]{PrzebindaLocal})
$$
|\det(\ad\,x)_{\n_{\overline 1}\to\n_{\overline 0}}|
=|\det(\ad\,x)_{\n_{\overline 0}\to\n_{\overline 1}}|.
$$
Hence, by (\ref{3.23.5}),
$$
|\det(\ad\,x)_{\n_{\overline 1}\to\n_{\overline 0}}|
=\left|\det(\ad\,x^2)_{\n_{\overline 0}}\right|^{\frac{1}{2}}
$$
and the last equality follows.
\end{prf}
Fix a $\K$-invariant norm $|\cdot |$ on the real vector space $\s$. (For example $|x|^2=-\langle \theta x,x\rangle$, $x\in \Sg$.)
The following theorem guarantees that the Weyl Harish-Chandra integration formula of Theorem \ref{Weyl integration} extends to Schwartz functions on $\W=\ss1$.
\begin{thm}\label{3.30}
For any constant $0\leq M<\infty$ there are constants $0\leq N, C<\infty$ such that for all $x\in\reg{\hs1}$,
\begin{eqnarray*}
&&|\pi_{\so/\hs1^2}(x^2)|\int_{\Sg/\Sg^{\hs1}} \int_{\ss1(\V^0)} (1+|g.(x+y)|)^{-N}\,dy\,d(g\Sg^{\hs1})\\
\leq &&C\cdot \left|\det(\ad\,x^2)_{\n_{\overline 0}}\right|^{\frac{1}{2}}\cdot(1+\left|x|_{\l_{0\overline 1}}\right|)^{-M}
\cdot\prod_{i=1}^nI_M(\left|x|_{\l_{i\overline 1}}\right|)\cdot 
\prod_{i=n+1}^{n+m}(1+\left|x|_{\l_{i\overline 1}}\right|)^{-M}.
\end{eqnarray*}
\end{thm}
\begin{prf}
This is immediate from Corollary \ref{3.25}, 
Lemmas \ref{3.5}, \ref{3.7} and Proposition \ref{3.3}.
\end{prf}
\section{Some properties of the invariant eigendistributions on the symplectic space.}
\label{section:diffeq}

The group $\Sg$ acts on $\ss1$ via the adjoint action. Hence, we have the permutation representation on the Schwartz space $S(\ss1)$:
\begin{equation}\label{4.1}
\pi(s)\phi(x)=\phi(s^{-1}.x) \qquad (s\in \Sg, x\in \ss1, \phi\in S(\ss1)).
\end{equation}
with the derivative
\begin{equation}\label{4.2}
\pi(z)\phi(x)=\frac{d}{dt}\phi(\exp(-tz).x)|_{t=0}
=\partial(-[z,x])\phi(x) \qquad (z\in\so, x\in \ss1, \phi\in S(\ss1)).
\end{equation}
As mentioned in the Introduction, a distribution $f\in S^*(\ss1)$ shall be called an invariant eigendistribution if 
\begin{eqnarray}\label{4.3}
\pi(s)f=f\qquad (s\in \Sg)
\end{eqnarray}
and
\begin{eqnarray}\label{4.4}
T(z)\natural f=\gamma(z) f\qquad (z\in \U(\so)^{\Sg}),
\end{eqnarray}
where $\gamma: \U(\so)^{\Sg}\to\Bbb C$ is an algebra homomorphism. Then (\ref{4.4}) is the system of equations explained in the introduction. Notice that if $\Sg$ is connected then (\ref{4.3}) is equivalent to
\begin{eqnarray}\label{4.5}
\partial([z,x])f(x)=0 \qquad (z\in\so, x\in\ss1).
\end{eqnarray}
Thus in this case $f$ is an invariant eigendistribution if it satisfies the two systems of differential equations (\ref{4.4}) and (\ref{4.5}).

Let us identify $\ss1$ with the dual $\ss1^*$ by
\begin{equation*}
x'(x)=\langle x,x'\rangle \qquad (x,x'\in\ss1).
\end{equation*}
This leads to an identification of the cotangent bundle to $\ss1$ with $\ss1\times\ss1$:
\begin{equation*}
T^*\ss1=\ss1\times\ss1^*=\ss1\times\ss1.
\end{equation*}
\begin{lem}\label{4.11}
The characteristic variety of the system of equations (\ref{4.4}) and (\ref{4.5}) is equal to the set of all the pairs $(x,y)\in \ss1\times\ss1$ such that the anticommutant $\{x,y\}=0$ and $y$ is nilpotent.
\end{lem}
\begin{prf}
The principal symbol of the differential operator (\ref{4.5}) may be computed as follows
\begin{equation*}
\underset{t\to\infty}{\lim}\,t^{-1}e^{-it\langle x,y\rangle}\partial([z,x])
e^{it\langle x,y\rangle}=i\langle [z,x],y\rangle.
\end{equation*}
Hence, the characteristic variety of the joined system (\ref{4.5}) consists of pairs $(x,y)$ such that
\begin{equation*}
\langle [z,x],y\rangle=0 \qquad (z\in \so).
\end{equation*}
In other words, $y$ is orthogonal to the space $[\so,x]$. By \cite[Lemma 3.5]{PrzebindaLocal} this last condition is equivalent to the vanishing of the anticommutant, $\{x,y\}=0$. But formula (7.22) in \cite{PrzebindaUnipotent} implies that the characteristic variety of the system (\ref{4.4}) consists of the pairs $(x,y)$ such that $y$ is annihilated by all the non-constant, complex valued, $\Sg$-invariant polynomials on $\ss1$, i.e. $y$ is nilpotent.
\end{prf}

An element $x\in\ss1$ is called \textit{regular} if the semisimple part of its Jordan decomposition is regular.    We shall denote by $\reg{\ss1}\subseteq \ss1$ the subset of all the regular elements.  

\begin{lem}\label{4.12}
The characteristic variety of the system of equations (\ref{4.4}) and (\ref{4.5}) over the set of the regular elements $\reg{\ss1}\subseteq \ss1$ is equal to $\reg{\ss1}\times\{0\}$.
\end{lem}
\begin{prf}
We may assume that $(\Sg,\s)$ is a complex supergroup. 
Let $x\in \reg{\ss1}$ and let $x=x_s+x_n$ be Jordan decomposition of $x$. Then, by \cite[Theorems 6.10(b) and 4.4(a)]{PrzebindaLocal}, the supergroup $(\Sg^{x_s^2}, \ss1)$ is isomorphic to the direct product of several copies of the supergroup whose underlying dual pair coincides with $(\GL_1(\C), \GL_1(\C))$, and one copy of $(\Og_1(\C), \Sp_{2n}(\C))$, if $x_n\ne 0$.

Let $y\in \ss1$ be nilpotent and anticommute with $x$. We need to show that $y=0$. It will suffice to show that $y$ is semisimple. Notice that $y$ commutes with $x^2=x_s^2+x_sx_n+x_nx_s+x_n^2=x_s^2$. 
Hence, we may assume that the underlying dual pair is either $(\GL_1(\C), \GL_1(\C))$ or $(\Og_1(\C), \Sp_{2n}(\C))$. In the first case the defining module $\V=\V_{\overline 0}\oplus \V_{\overline 1}$, with $dim\,\V_{\overline 0}=dim\,\V_{\overline 1}=1$. Let $v_0\in \V_{\overline 0}$ and $v_1\in \V_{\overline 1}$. We may assume that 
\begin{equation*}
x(v_0)=\xi v_1\ \text{and}\ x(v_1)=\pm \xi v_0,
\end{equation*}
where $\xi\in\C^\times$. Then the anticommutant of $x$ consists of elements $x'$ such that
\begin{equation*}
x'(v_0)=\xi' v_1\ \text{and}\ x'(v_1)=\mp \xi' v_0,
\end{equation*}
where $\xi'\in\C$, see \cite[(13.47)]{PrzebindaLocal}. In particular $x'$, and hence $y$, is semisimple.

In the second case $x$ is any non-zero element of $\ss1$. Recall the following formula
\begin{equation*}
\langle z,\{x,y\}\rangle=\langle y, [z,x]\rangle \qquad (z\in \so).
\end{equation*}
Since the set of all the commutators $[z,x]$ coincides with $\ss1$ (because the set of the non-zero elements of $\ss1$ is a single $\Sp_{2n}(\C)$-orbit) and $\{x,y\}=0$, we see that $y$ is orthogonal to $\ss1$ and hence equal to zero.
\end{prf}
By combining Lemmas \ref{4.11} and \ref{4.12} with \cite{Hormander} we obtain the following corollary.
\begin{cor}\label{4.13}
The restriction of an invariant eigendistribution on $\ss1$ to $\reg{\ss1}$ is a smooth function.
\end{cor}

\section*{Appendix A}
\label{appendix:Rossmann}
\setcounter{thh}{0}
\renewcommand{\thethh}{A.\fontindex{thh}}
\setcounter{equation}{0}
\renewcommand{\theequation}{A.\fontindex{equation}}
\renewcommand{\theequation}{A.\fontindex{equation}}
In this appendix we verify the assumptions of \cite[Corollary 5.4, page 283]{RossmannNilpotent} for a minimal 
non-zero nilpotent orbit in a real symplectic Lie algebra 
\begin{equation}\label{sp}
\s\p_{2n}(\R)=\left\{X=
\left(
\begin{array}{cc}
A & B\\
C & -A^t
\end{array}
\right);\ A, B=B^t, C=C^t\in \g\l_n(\R)
\right\}
\end{equation}
Let $\h$ be the elliptic Cartan subalgebra consisting of the matrices $X=X(x_1, x_2,\dots, x_n)$, as in (\ref{sp}), with $A=0$, $B$ diagonal with diagonal entries $x_1$, $x_2$,\dots, $x_n$, and $C=-B$. Also, let $\k\subseteq \s\p_{2n}(\R)$ be the subset of the skew-symmetric matrices.

The set of regular elements $\reg\h$ in $\h$ consists of the $X$ such that the $x_j$ are all distinct and non-zero. Let $\h^+\subseteq \reg\h$ be the Weyl chamber defined by the condition $x_1>x_2>\dots >x_n>0$.  For $t>0$ let
\begin{equation*}
a(t)=
\left(
\begin{array}{cccc|cccc}
t & 0 &\dots & 0&0&0& \dots &0\\
0 & 1 &\dots & 0&0&0& \dots &0\\
\vdots & \vdots &\ddots & \vdots & \vdots & \vdots &\ddots & \vdots\\
0 & 0 &\dots & 1&0&0& \dots&0\\ \hline
0 & 0 &\dots & 0&t^{-1}&0& \dots&0\\
0 & 0 &\dots & 0&0&1& \dots&0\\
\vdots & \vdots &\ddots & \vdots & \vdots & \vdots & \ddots & \vdots\\
0 & 0 &\dots & 0&0&0& \dots&1\\
\end{array}
\right).
\end{equation*}
Then $a(t)\in\Sp_{2n}(\R)$ and
\begin{equation*}
a(t)Xa(t)^{-1}=
\left(
\begin{array}{cccc|cccc}
0 & 0 &\dots & 0&t^2x_1&0& \ddots &0\\
0 & 0 &\dots & 0&0&x_2& \ddots &0\\
\vdots & \vdots &\ddots & \vdots &\vdots & \vdots & \ddots &\vdots \\
0 & 0 &\dots & 0&0&0& \dots &x_n\\ \hline
-t^{-2}x_1 & 0 &\dots & 0&0&0& \dots &0\\
0 & -x_2 &\dots & 0&0&0& \dots &0\\
\vdots & \vdots &\ddots & \vdots &\vdots & \vdots & \ddots &\vdots \\
0 & 0 &\dots & -x_n&0&0& \dots &0\\
\end{array}
\right).
\end{equation*}
Hence,
\begin{equation}\label{nilpotent element}
\lim\limits_{t \to \infty} t^{-2}a(t)Xa(t)^{-1}=
\left(
\begin{array}{cccc|cccc}
0 & 0 & \dots & 0 &x_1&0 &\dots &0\\
0 & 0 & \dots & 0 &0 & 0 &\dots &0\\
\vdots & \vdots &\ddots & \vdots &\vdots & \vdots & \ddots &\vdots \\
0 & 0 & \dots & 0 &0 & 0 &\dots &0\\ \hline
0 & 0 & \dots & 0 &0 & 0 &\dots &0\\ 
0 & 0 & \dots & 0 &0 & 0 &\dots &0\\
\vdots & \vdots &\ddots & \vdots &\vdots & \vdots & \ddots &\vdots \\
0 & 0 & \dots & 0 &0 & 0 &\dots &0\\
\end{array}
\right).
\end{equation}
If $x_1>0$ then (\ref{nilpotent element}) belongs to one minimal nilpotent orbit $\Og$, and if $x_1<0$ it belongs to the other nilpotent orbit equal to $-\Og$. Let $\Og_\C$ be the complexification of $\Og$. This is the unique minimal nilpotent orbit for the complex symplectic group $\Sp_{2n}(\C)$. We see that
\begin{equation*}\label{O}
\Og\subseteq\Og_\C\cap \big(\text{closure of}\ \Sp_{2n}(\R).\R^+X\big) \qquad (X\in\h^+).
\end{equation*}
Let
\[
J=\left(
\begin{array}{rrr}
0 & -I\\
I & 0
\end{array}
\right)
\]
The map $Y\to JY$ transforms elements of the symplectic Lie algebra into the elements the space of the symmetric matrices and intertwines the conjugation action of $\Sp_{2n}(\R)$ on the Lie algebra with the natural action on the symmetric matrices. In particular the image of an $\Sp_{2n}(\R)$-orbit of an element $X\in \h^+$ consists of positive definite matrices. Hence the image of $\text{closure of}\ \Sp_{2n}(\R).\R^+X$ consists of positive semidefinite matrices. Therefore
\begin{equation*}\label{O}
-\Og\nsubseteq\text{closure of}\ \Sp_{2n}(\R).\R^+X \qquad (X\in\h^+).
\end{equation*}
Since $\Og_\C\cap\s\p_{2n}(\R)=\Og\cup(-\Og)$, we see that
\begin{equation}\label{O}
\Og=\Og_\C\cap \big(\text{closure of}\ \Sp_{2n}(\R).\R^+X\big) \qquad (X\in\h^+).
\end{equation}
This verifies \cite[Hypothesis (O), page 280]{RossmannNilpotent}.

Let us choose a positive root system for the roots of $\h$ in $\g_\C$ so that
$\pi_{\s\p_{2n}(\R)/\h}^{short}(X)=\prod_{1\leq j<k\leq n}(ix_k-ix_j)(ix_k+ix_j)$. The Weyl group $W(\s\p_{2n}(\R),\h)$ consists of all the permutations and the sign changes of the $x_j$. This group acts on $\C[\h]$,  the space of the complex valued polynomials on $\h$, in the obvious way. The representation $\rho_{\Og_\C}$ of $W(\s\p_{2n}(\R),\h)$ attached to the orbit $\Og_\C$ via Springer correspondence may be realized on the subspace of the polynomials generated by $\pi_{\s\p_{2n}(\R)/\h}^{short}$. (One way to see it is to use Wallach's theorem \cite{WallachSpringer} saying that $\rho_{\Og_\C}$ is generated by the product of $\hat\mu_{\Og}$, the Fourier transform of the invariant measure on $\Og$ restricted to $\h^+$, times $\pi_{\s\p_{2n}(\R)/\h}$ the product of all the roots of $\h$ in $\s\p_{2n}(\R)$. The Fourier transform is well known and is equal to a constant multiple of the function $\chc$ given in (\ref{1.2}), see for example \cite[proof of Proposition 9.3]{PrzebindaCauchy}. In particular the restriction to $\h^+$ coincides with the reciprocal of the product of the long roots. Therefore $\hat\mu_{\Og}|_{\h^+}\pi_{\s\p_{2n}(\R)/\h}|_{\h^+}$ is a non-zero constant multiple of $\pi_{\s\p_{2n}(\R)/\h}^{short}$.) The representation $\rho_{\Og_\C}$ is one dimensional and its restriction to $W(\k,\h)$, the subgroup of all the permutations of the $x_j$, coincides with
the sign representation. Therefore the latter is contained in $\rho_{\Og_\C}$ with multiplicity one. This verifies the assumptions of \cite[Corollary 5.4]{RossmannNilpotent}.
\section*{Appendix B}
\label{appendix:good nilpotent elements}
\setcounter{thh}{0}
\renewcommand{\thethh}{B.\fontindex{thh}}
\setcounter{equation}{0}
\renewcommand{\theequation}{B.\fontindex{equation}}
\renewcommand{\theequation}{B.\fontindex{equation}}
Here we verify Proposition \ref{good nilpotent elements}.  
Recall that our dual pair acts on a defining module $\V=\V_{\overline 0}\oplus \V_{\overline 1}$, which is  a finite dimensional $\Zb/2\Zb$-graded vector space over $\Bbb D=\R$, $\C$ or $\Bbb H$.  For $z\in \ss1$, the pair $(z,\V)$ is called decomposable if there are two $\Zb/2\Zb$-graded non-zero subspaces $\V', \V''\subseteq \V$ preserved by $z$ (and orthogonal if the corresponding dual pair is of type I) such that $\V=\V'\oplus \V''$. Otherwise  $(z,\V)$ is called indecomposable. 

We shall verify Proposition \ref{good nilpotent elements} under the assumption that the nilpotent element $z\in\ss1$ is such that $(z,\V)$ is indecomposable. (We shall see at the end of this appendix that this implies the result for general nilpotent elements.) 

For dual pairs of type I such elements are classified in \cite[Theorem 5.2]{DaszKrasPrzebindaK-S2}. There are seven cases to consider (Case I.a - Case I.g below). For pairs of type II there are two cases, see \cite[Theorem 5.5]{PrzebindaLocal} (Case II.a  and  Case II.b below). (Some misprints which occur in \cite[Theorem 5.2]{DaszKrasPrzebindaK-S2} have been corrected in \cite[Theorem 5.4]{PrzebindaLocal}. There is one remaining misprint in \cite[Theorem 5.4]{PrzebindaLocal} which we correct in what follows.) In each case we describe the graded vector space $\V$, the automorphism $\theta$ (\ref{3.17}), the Cartan subspace $\h_{\overline 1}\subseteq \ss1$, the abelian Lie algebra $\a\subseteq \h_{\overline 1}^2$ and the decomposition of $z$ into root spaces for the action of $\a$. 

\medskip

We begin with the dual pairs of type I. In this case there is a possibly trivial involution $\Bbb D\ni a\to \overline a\in\Bbb D$ on the division algebra $\Bbb D$.  We shall assume that the space $\V_{\overline 0}$ is equipped with a non-degenerate  hermitian (or symmetric) form $(\ ,\ )_{\overline 0}$ and that $\V_{\overline 1}$ with a non-degenerate  skew-hermitian (or skew-symmetric) form $(\ ,\ )_{\overline 1}$. Let $(\ ,\ )=(\ ,\ )_{\overline 0}\oplus (\ ,\ )_{\overline 1}$.
Then, in particular, $x\in\End(\V)$ belongs to $\ss1$ if and only if $(xu,v)=(u,sxv)$ for all $u,v\in\V$, where $s(v_{even}+v_{odd})=v_{even}-v_{odd}$.
Recall that the signature of a hermitian (or symmetric) form is the difference of the dimension of a maximal subspace where the form is positive definite and the dimension of a maximal subspace where the form is negative definite. We shall write $\sgn((\cdot ,\cdot )_{\overline 0})=1$ if the form $(\cdot ,\cdot )_{\overline 0}$ has positive signature, $\sgn((\cdot ,\cdot )_{\overline 0})=-1$ if it has negative signature and $\sgn((\cdot ,\cdot )_{\overline 0})=0$ if the signature is zero. In each case we'll describe an element $T\in\Sg$ such that for a specific basis of $\V$ consisting of vectors $v_k$,  $(Tv_k,v_k)>0$. Consequently the conjugation by $T$ coincides with the automorphism $\theta$ (\ref{3.17}), unless the involution $\Bbb D\ni a\to \overline a\in\Bbb D$ is trivial.
In that case  $\theta$ is the conjugation by $T$ composed with the conjugation on $\s$ induced by the complex conjugation on $\V$ which leaves the vectors $v_k$ fixed.
Set 
\begin{equation}\label{delta(k)}
\delta(k)=(-1)^{k(k-1)/2}=\left\{
\aligned
&1\ \text{if}\ k\in 4\Bbb Z\ \text{or}\ 4\Bbb Z + 1,\\
-&1\ \text{if}\ k\in 4\Bbb Z+2\ \text{or}\ 4\Bbb Z + 3.
\endaligned
\right.
\end{equation}
Then, in particular, $\delta(k+1)=(-1)^k\delta(k)$, if $m$ is divisible by $4$ then $\delta(m+2-k)=-(-1)^k\delta(k)$,  if $m$ is even but is not divisible by $4$ then $\delta(m-k)=(-1)^{k+1}\delta(k)$ .

\newpage
\enlargethispage*{3truecm}

\underline{Case I.a: $\Sg=\Og_{p+1,p}\times \Sp_{2p}(\R)$, $\Og_{2p+1}(\C)\times \Sp_{2p}(\C)$, $\Ug_{p+1,p}\times \Ug_{p,p}$ or $\Sp_{p+1,p}\times \Og^*_{4p}$ }

\begin{eqnarray*}
&&m\in 4\Bbb Z,\ m>0;\\
&&\V=\sum_{k=0}^m\Bbb D v_k,\ v_{even}\in \V_{\overline 0},\ v_{odd}\in \V_{\overline 1};\\
&&v_k=z^kv_0\ne 0, 0\leq k\leq m, zv_m=0;\\ 
&&(v_k,v_l)=0\ \text{if}\ l\ne m-k,
( v_k, v_{m-k}) =\delta(k)\delta(\frac{m}{2}) \sgn((\cdot ,\cdot )_{\overline 0}),\\
&&\text{where, by definition,}\ \sgn((\cdot ,\cdot )_{\overline 0})=1\\ 
&&\text{if}\ \Bbb D=\Bbb C\ \text{and the involution}\ \Bbb D\ni a\to \overline a\in\Bbb D\ \text{is trivial}; \hfill\null
\end{eqnarray*}

Here,
\begin{eqnarray*}
&&Tv_k=t_kv_{m-k},\ 
t_k=(-1)^k\delta(k)\delta(\frac{m}{2}) \sgn((\cdot ,\cdot )_{\overline 0}),\ 0\leq k\leq m.
\end{eqnarray*}
The Cartan subspace $\h_{\overline 1}$ 
consists of the linear maps $x$ defined by
\begin{eqnarray*}
&&xv_{2j}=x_{j}v_{2j+1},\ xv_{2j+1}=x_{j}v_{2j},\ xv_{\frac{m}{2}}=0,\\ 
&&xv_{m-2j}=-x_{j}v_{m-2j-1},\ xv_{m-2j-1}=x_{j}v_{m-2j},\\
&&\text{if}\ \Bbb D\ne \Bbb H,\ \text{then}\ \ x_j\in\Bbb D,\\
&&\text{if}\ \Bbb D=\Bbb H,\ \text{then}\ x_j\in \C = \text{the centralizer of  $i$ in}\ \Bbb H,\  0\leq 2j<\frac{m}{2}.
\end{eqnarray*}
This is the direct sum of the indecomposable Cartan subspaces which occur in \cite[Theorem 6.2(a)]{PrzebindaLocal}.

The Lie algebra $\a$ consists of the linear maps $a$ defined by
\begin{eqnarray*}
&&av_k=a_kv_{k},\ a_k\in\R,\ 0\leq k\leq m,\\
&&a_k=-a_{m-k},\ 0\leq k< \frac{m}{2},\ a_{\frac{m}{2}}=0,\\
&&a_{2j}=a_{2j+1}\ \text{for}\ 0\leq 2j< \frac{m}{2},
\end{eqnarray*}

For $0\leq k<\frac{m}{2}$ define a linear map $z(k)$ by
\begin{eqnarray*}
&&z(k)v_k=v_{k+1},\ z(k)v_{m-k-1}=v_{m-k},\\
&&z(k)v_j=0,\ j\notin\{k, m-k-1\}.
\end{eqnarray*}
Then $z=\sum_{k=0}^{\frac{m}{2}-1}z(k)$, $z(k)\in\ss1$ and for $a\in \a$, $[a,z(k)]=(a_{k+1}-a_k)z(k)$. In particular $z(even)\in \sum_{i=1}^m\l_{i\overline 1}$ and $z(odd)\in\n_1$.

For a fixed $j$ with $0\leq 2j<\frac{m}{2}$ and for $t>0$ define $b=b(t)\in\Sg$ by
\begin{eqnarray*}
&&bv_{2j}=tv_{2j},\ bv_{2j+1}=t^{-1}v_{2j+1},\\
&&bv_{m-2j}=t^{-1}v_{m-2j},\ bv_{m-2j-1}=tv_{m-2j-1},\\
&&bv_k=v_k\ \text{for}\ k\notin\{2j, 2j+1, m-2j, m-2j-1\}.
\end{eqnarray*}
Also, let $x=x(t)\in\h_{\overline 1}$ be such that $x_j=t^2$ and $x_i=0$ for $i\ne j$. Then $\lim_{t\to 0} bxb^{-1}=z(2j)$.
Hence, $z_\l=\sum_j z(2j)$ is the limit of elements of the $\L$-orbits passing through $\h_{\overline 1}$.

\newpage
\enlargethispage*{3truecm}
\underline{Case I.b: $\Sg=\Ug_{p,p}\times \Ug_{p+1,p}$ or $\Sp_{p,p}\times \Og^*_{2(2p+1)}$}

\begin{eqnarray*}
&&m\in 4\Bbb Z,\ m>0,\ \Bbb D\ne \R,\ \text{the involution}\ \Bbb D\ni a\to\overline a\in\Bbb D\ \text{is not trivial};\\
&&\V=\sum_{k=1}^{m+1}\Bbb D v_k,\ v_{even}\in \V_{\overline 0},\ v_{odd}\in \V_{\overline 1};\\
&&v_{k+1}=z^kv_1\ne 0, 0\leq k\leq m, zv_{m+1}=0;\\ 
&&(v_k,v_l)=0\ \text{if}\ l\ne m+2-k,
( v_k, v_{m+2-k}) =\delta(k)( v_1, v_{m+1}),\\
&&( v_1, v_{m+1})=i\,\sgn(-i(\cdot ,\cdot )_{\overline 1})\delta(1+\frac{m}{2}),\ \text{if}\ \Bbb D=\Bbb C,\\ 
&&( v_1, v_{m+1})=i,\ \text{if}\ \Bbb D=\Bbb H;
\end{eqnarray*}
Here,
\begin{eqnarray*}
&&Tv_k=t_kv_{m+2-k},\ 
t_k=(-1)^k\delta(k)\overline{( v_1, v_{m+1})},\ 1\leq k\leq m+1.
\end{eqnarray*}
The Cartan subspace $\h_{\overline 1}$ 
consists of the linear maps $x$ defined by
\begin{eqnarray*}
&&xv_{2j+1}=x_{j}v_{2j+2},\ xv_{2j+2}=x_{j}v_{2j+1},\ xv_{\frac{m+2}{2}}=0,\\ 
&&xv_{m+1-2j}=-\overline{x_{j}}v_{m-2j},\ xv_{m-2j}=\overline{x_{j}}v_{m+1-2j},\\
&&x_j\in\C=\ \text{centralizer of}\ i\ \text{in}\ \Bbb D,\  1\leq 2j+1<\frac{m}{2}.
\end{eqnarray*}
This is the direct sum of the indecomposable Cartan subspaces which are isomorphic to those which occur in \cite[Theorem 6.2(a)]{PrzebindaLocal}. (One has to adjust the sesquilinear forms $\tau_0$ and $\tau_1$ listed in \cite[Theorem 6.2(a)]{PrzebindaLocal} in order to get the forms $(\ ,\ )_{\overline 0}$ and $(\ ,\ )_{\overline 1}$ we are working with here.)

The Lie algebra $\a$ consists of the linear maps $a$ defined by
\begin{eqnarray*}
&&av_k=a_kv_{k},\ a_k\in\R,\ 1\leq k\leq m+1,\\
&&a_k=-a_{m+2-k},\ 1\leq k< \frac{m+2}{2},\ a_{\frac{m+2}{2}}=0,\\
&&a_{2j+1}=a_{2j+2}\ \text{for}\ 1\leq 2j+1< \frac{m}{2}.
\end{eqnarray*}

For $1\leq k\leq \frac{m}{2}$ define a linear map $z(k)$ by
\begin{eqnarray*}
&&z(k)v_k=v_{k+1},\ z(k)v_{m+1-k}=v_{m+2-k},\\
&&z(k)v_j=0,\ j\notin\{k, m+1-k\}.
\end{eqnarray*}
Then $z=\sum_{k=1}^{\frac{m}{2}}z(k)$, $z(k)\in\ss1$ and for $a\in \a$, $[a,z(k)]=(a_{k+1}-a_k)z(k)$. In particular $z(odd)\in \sum_{i=1}^m\l_{i\overline 1}$ and $z(even)\in\n_1$.

For a fixed $j$ with $1\leq 2j+1<\frac{m}{2}$ and for $t>0$ define $b=b(t)\in\Sg$ by
\begin{eqnarray*}
&&bv_{2j+1}=tv_{2j+1},\ bv_{2j+2}=t^{-1}v_{2j+2},\\
&&bv_{m+1-2j}=t^{-1}v_{m+1-2j},\ bv_{m-2j}=tv_{m-2j},\\
&&bv_k=v_k\ \text{for}\ k\notin\{2j+1, 2j+2, m+1-2j, m-2j\}.
\end{eqnarray*}
Also, let $x=x(t)\in\h_{\overline 1}$ be such that $x_j=t^2$ and $x_i=0$ for $i\ne j$. Then $\lim_{\to 0} bxb^{-1}=z(2j+1)$.
Hence, $z_\l=\sum_j z(2j+1)$ is the limit of elements of the $\L$-orbits passing through $\h_{\overline 1}$.
\newpage

\enlargethispage*{3truecm}
\underline{Case I.c: 
$\Sg=\Og_{p,p}\times \Sp_{2(p+1)}(\R)$, $\Og_{2p}(\C)\times \Sp_{2(p+1)}(\C)$, $\Ug_{p+1,p}\times \Ug_{p,p}$ or $\Sp_{p+1,p}\times \Og^*_{2p}$}

\begin{eqnarray*}
&&m\in 4\Bbb Z,\ m>0,\ \Bbb D\ne \Bbb H,\ \text{the involution}\ \Bbb D\ni a\to\overline a\in\Bbb D\ \text{is trivial};\\
&&\V=\sum_{k=1}^{m+1}(\Bbb D v_k+\Bbb D v'_k);\ v_{even},\ v'_{even}\in \V_{\overline 0};\ v_{odd},\ v'_{odd}\in \V_{\overline 1};\\
&&v_{j+1}=z^jv_1\ne 0,\ v'_{j+1}=z^jv'_1\ne 0,\ 0\leq j\leq m,\ zv_{m+1}=zv'_{m+1}=0;\\ 
&&( v_k, v'_{m+2-k}) =\delta(k),\ ( v'_k, v_{m+2-k})=-\delta(k),\ 1\leq k\leq m+1\\ 
&&\text{and all other pairings are zero};
\end{eqnarray*}
Here,
\begin{eqnarray*}
&&Tv_k=t_kv'_{m+2-k},\ Tv'_k=t'_kv_{m+2-k},\\ 
&&t_k=(-1)^k\delta(k),\ t'_k=-(-1)^k\delta(k),\ 1\leq k\leq m+1.
\end{eqnarray*}
The Cartan subspace $\h_{\overline 1}$ 
consists of the linear maps $x$ defined by
\begin{eqnarray*}
&&xv_{2j+1}=x_{j}v_{2j+2},\ xv_{2j+2}=x_{j}v_{2j+1},\ xv_{m+1}=0,\\ 
&&xv'_{m+1-2j}=-x_{j}v'_{m-2j},\ xv'_{m-2j}=x_{j}v'_{m+1-2j},\ xv'_1=0,\\
&&x_j\in\Bbb D,\  1\leq 2j+1<{m+1}{}.
\end{eqnarray*}
This is the direct sum of the indecomposable Cartan subspaces which are isomorphic to those which occur in \cite[Theorem 6.2(a)]{PrzebindaLocal}. (One has to adjust the sesquilinear forms $\tau_0$ and $\tau_1$ listed in \cite[Theorem 6.2(a)]{PrzebindaLocal} in order to get the forms $(\ ,\ )_{\overline 0}$ and $(\ ,\ )_{\overline 1}$ we are working with here.)

The Lie algebra $\a$ consists of the linear maps $a$ defined by
\begin{eqnarray*}
&&av_k=a_kv_{k},\ av'_{m+2-k}=-a_{k}v'_{m+2-k},\ a_k\in\R,\ 1\leq k\leq m,\\
&&av_{m+1}=0,av'_{1}=0,\\
&&a_{2j+1}=a_{2j+2}\ \text{for}\ 1\leq 2j+1< m+1.
\end{eqnarray*}

For $1\leq k\leq m$ define a linear map $z(k)$ by
\begin{eqnarray*}
&&z(k)v_k=v_{k+1},\ z(k)v'_{m+1-k}=v'_{m+2-k},\\
&&z(k)v_j=0,\ z'(k)v'_{m+1-j}=0,\ j\ne k.
\end{eqnarray*}
Then $z=\sum_{k=1}^{m}z(k)$, $z(k)\in\ss1$ and for $a\in \a$, $[a,z(k)]=(a_{k+1}-a_k)z(k)$. In particular $z(odd)\in \sum_{i=1}^m\l_{i\overline 1}$ and $z(even)\in\n_1$.

For a fixed $j$ with $1\leq 2j+1<m+1$ and for $t>0$ define $b=b(t)\in\Sg$ by
\begin{eqnarray*}
&&bv_{2j+1}=tv_{2j+1},\ bv_{2j+2}=t^{-1}v_{2j+2},\\
&&bv_{m+1-2j}=t^{-1}v_{m+12j},\ bv_{m-2j}=tv_{m-2j},\\
&&bv_k=v_k\ \text{for}\ k\notin\{2j+1, 2j+2, m+1-2j, m-2j\}.
\end{eqnarray*}
Also, let $x=x(t)\in\h_{\overline 1}$ be such that $x_j=t^2$ and $x_i=0$ for $i\ne j$. Then $\lim_{t\to 0}bxb^{-1}=z(2j+1)$.
Hence, $z_\l=\sum_j z(2j+1)$ is the limit of elements of the $\L$-orbits passing through $\h_{\overline 1}$.
\newpage

\underline{Case I.d: $\Sg=\Ug_{p,p}\times \Ug_{p,p-1}$ or $\Sp_{p,p}\times \Og^*_{2(2p-1)}$}

\begin{eqnarray*}
&&m\in 2\Bbb Z\setminus 4\Bbb Z,\ m>0,\ \Bbb D\ne \R,\ \text{the involution}\ \Bbb D\ni a\to \overline a\in\Bbb D\ \text{is not trivial};\\
&&\V=\sum_{k=0}^m\Bbb D v_k,\ v_{even}\in \V_{\overline 0},\ v_{odd}\in \V_{\overline 1};\\
&&v_k=z^kv_0\ne 0, 0\leq k\leq m, zv_m=0;\\ 
&&(v_k,v_l)=0\ \text{if}\ l\ne m-k,
( v_k, v_{m-k}) =\delta(k)ic,\ c=\pm 1.
\end{eqnarray*}
Here,
\begin{eqnarray*}
&&Tv_k=t_kv_{m-k},\ 
t_k=(-1)^k\delta(k)ic,\ 0\leq k\leq m.
\end{eqnarray*}
The Cartan subspace $\h_{\overline 1}$ 
consists of the linear maps $x$ defined by
\begin{eqnarray*}
&&xv_{2j}=x_{j}v_{2j+1},\ xv_{2j+1}=x_{j}v_{2j},\ xv_{\frac{m}{2}}=0,\\ 
&&xv_{m-2j}=-\overline{x_{j}}v_{m-2j-1},\ xv_{m-2j-1}=\overline{x_{j}}v_{m-2j},\\
&&x_j\in \C = \text{the centralizer of  $i$ in}\ \Bbb D,\  0\leq 2j<\frac{m-2}{2},\\
&&xv_{\frac{m}{2}}=x_{\frac{m}{2}}\frac{1}{\sqrt{2}}(v_{\frac{m}{2}-1}+iv_{\frac{m}{2}+1}),\ 
x\frac{1}{\sqrt{2}}(v_{\frac{m}{2}-1}+iv_{\frac{m}{2}+1})=x_{\frac{m}{2}}v_{\frac{m}{2}},\\
&&x_{\frac{m}{2}}\in \frac{1}{\sqrt{2}}(i-1)\R,\\
&&x\frac{1}{\sqrt{2}}(v_{\frac{m}{2}-1}-iv_{\frac{m}{2}+1})=0.
\end{eqnarray*}
This is the direct sum of the indecomposable Cartan subspaces which occur in \cite[Theorem 6.2(a) and (b)]{PrzebindaLocal}.

The Lie algebra $\a$ consists of the linear maps $a$ defined by
\begin{eqnarray*}
&&av_k=a_kv_{k},\ a_k\in\R,\ 0\leq k\leq m,\\
&&a_k=-a_{m-k},\ 0\leq k< \frac{m}{2},\ a_{\frac{m}{2}}=0,\\
&&a_{2j}=a_{2j+1}\ \text{for}\ 0\leq 2j< \frac{m}{2},
\end{eqnarray*}

For $0\leq k<\frac{m-2}{2}$ define a linear map $z(k)$ by
\begin{eqnarray*}
&&z(k)v_k=v_{k+1},\ z(k)v_{m-k-1}=v_{m-k},\\
&&z(k)v_j=0,\ j\notin\{k, m-k-1\}.
\end{eqnarray*}
Then $z=\sum_{k=0}^{\frac{m-2}{2}}z(k)$, $z(k)\in\ss1$ and for $a\in \a$, $[a,z(k)]=(a_{k+1}-a_k)z(k)$. In particular $z(even)\in \sum_{i=1}^m\l_{i\overline 1}$ and $z(odd)\in\n_1$.

Set $p=\frac{m}{2}$, $u=\frac{-1}{\sqrt{2}}(1+i)$ and for $t\in\R^\times$ define
\begin{eqnarray*}
&&bv_{p-1}=tv_{p-1},\ bv_p=uv_p,\ bv_{p+1}=t^{-1}v_{p+1},\\
&&bv_k=v_k\ \text{if}\ k\notin\{p-1,p,p+1\}.
\end{eqnarray*}
Then $b\in \Sg$ and
\begin{eqnarray*}
&&bxb^{-1}v_{p-1}=t^{-1}u\frac{x_p}{\sqrt{2}}v_p,\ bxb^{-1}v_{p}=u^{-1}\frac{x_p}{\sqrt{2}}(tv_{p-1}+t^{-1}iv_{p+1}),\\
&&bxb^{-1}v_{p+1}=-itu\frac{x_p}{\sqrt{2}}v_p.
\end{eqnarray*}
Choose $x_p=-iu\sqrt{2}t$. Then $z(p-1)=\lim_{t\to 0}bxb^{-1}$. Thus $z(p-1)\in \l_{0\overline 1}$ is the limit of the semisimple elements of $\l_{0\overline 1}$.

For a fixed $j$ with $0\leq 2j<\frac{m-2}{2}$ and for $t>0$ define $b=b(t)\in\Sg$ by
\begin{eqnarray*}
&&bv_{2j}=tv_{2j},\ bv_{2j+1}=t^{-1}v_{2j+1},\\
&&bv_{m-2j}=t^{-1}v_{m-2j},\ bv_{m-2j-1}=tv_{m-2j-1},\\
&&bv_k=v_k\ \text{for}\ k\notin\{2j, 2j+1, m-2j, m-2j-1\}.
\end{eqnarray*}
Also, let $x=x(t)\in\h_{\overline 1}$ be such that $x_j=t^2$ and $x_i=0$ for $i\ne j$. Then $lim_{t\to 0} bxb^{-1}=z(2j)$.
Hence, $z_\l=\sum_j z(2j)$ is the limit of elements of the $\L$-orbits passing through $\h_{\overline 1}$.
\bigskip


\underline{Case I.e: $\Sg=\Og_{p,p-1}\times \Sp_{2p}(\R)$, $\Og_{2p-1}(\C)\times \Sp_{2p}(\C)$, $\Ug_{p,p-1}\times \Ug_{p,p}$ or\newline 
$\Sp_{p,p-1}\times \Og^*_{4p}$ }

\medskip

(This is the only case where (\ref{basic assumption-one}) or (\ref{basic assumption-two}) occurs. 
It happens when the involution $\Bbb D\ni a\to\overline a\in \Bbb D$ is trivial.)
\begin{eqnarray*}
&&m\in 2\Bbb Z\setminus 4\Bbb Z,\ m>0;\\
&&\V=\sum_{k=1}^{m+1}\Bbb D v_k,\ v_{even}\in \V_{\overline 0},\ v_{odd}\in \V_{\overline 1};\\
&&v_k=z^kv_1\ne 0, 0\leq k\leq m, zv_{m+1}=0;\\ 
&&(v_k,v_l)=0\ \text{if}\ l\ne m+2-k,
( v_k, v_{m+2-k}) =\delta(k)(v_1,v_{m+1}),\\
&&(v_1,v_{m+1})=\delta(\frac{m+2}{2})\sgn((\cdot ,\cdot )_{\overline 0}).
\end{eqnarray*}
Here,
\begin{eqnarray*}
&&Tv_k=t_kv_{m+2-k},\ 
t_k=(-1)^k\delta(k)(v_1,v_{m+1}),\ 0\leq k\leq m.
\end{eqnarray*}
The Cartan subspace $\h_{\overline 1}$ 
consists of the linear maps $x$ defined by
\begin{eqnarray*}
&&xv_{2j-1}=x_{j}v_{2j},\ xv_{2j}=x_{j}v_{2j-1},\\ 
&&xv_{m+3-2j}=-\overline{x_{j}}v_{m+2-2j},\ xv_{m+2-2j}=\overline{x_{j}}v_{m+3-2j},\\
&&x_j\in \C = \text{the centralizer of  $i$ in}\ \Bbb D,\  0<2j<\frac{m+2}{2},\\
&&\text{if the involution $\Bbb D\ni a\to \overline a\in \Bbb D$ is trivial, then}\  xv_{\frac{m+2}{2}}=0,\\ 
&&\text{if the involution $\Bbb D\ni a\to \overline a\in \Bbb D$ is not trivial, then}\\
&&xv_{\frac{m+2}{2}}=x_{\frac{m+2}{2}}\frac{1}{\sqrt{2}}(v_{\frac{m}{2}}+iv_{\frac{m}{2}+2}),\ 
x\frac{1}{\sqrt{2}}(v_{\frac{m}{2}}+iv_{\frac{m}{2}+2})=x_{\frac{m+2}{2}}v_{\frac{m+2}{2}},\\
&&x_{\frac{m+2}{2}}\in \frac{1}{\sqrt{2}}(i-1)\R,\\
&&x\frac{1}{\sqrt{2}}(v_{\frac{m}{2}}-iv_{\frac{m}{2}+2})=0.
\end{eqnarray*}
This is the direct sum of the indecomposable Cartan subspaces which occur in \cite[Theorem 6.2(a) and (b)]{PrzebindaLocal}.

The Lie algebra $\a$ consists of the linear maps $a$ defined by
\begin{eqnarray*}
&&av_k=a_kv_{k},\ a_k\in\R,\ 1\leq k\leq m+1,\\
&&a_k=-a_{m+2-k},\ 1\leq k< \frac{m+2}{2},\ a_{\frac{m+2}{2}}=0,\\
&&a_{2j-1}=a_{2j}\ \text{for}\ 0<2j< \frac{m+2}{2},
\end{eqnarray*}

For $1\leq k<\frac{m+2}{2}$ define a linear map $z(k)$ by
\begin{eqnarray*}
&&z(k)v_k=v_{k+1},\ z(k)v_{m+1-k}=v_{m+2-k},\\
&&z(k)v_j=0,\ j\notin\{k, m+1-k\}.
\end{eqnarray*}
Then $z=\sum_{k=1}^{\frac{m+2}{2}}z(k)$, $z(k)\in\ss1$ and for $a\in \a$, $[a,z(k)]=(a_{k+1}-a_k)z(k)$. In particular $z(odd)\in \sum_{i=1}^m\l_{i\overline 1}$, $z(\frac{m+2}{2})\in \l_{0\overline 1}$ and $z(2j)\in\n_1$  for $0<2j<\frac{m+2}{2}$. 

Assume that the involution $\Bbb D\ni a\to\overline a\in \Bbb D$ is not trivial.
Set $p=\frac{m+2}{2}$, $u=\frac{-1}{\sqrt{2}}(1+i)$ and for $t\in\R^\times$ define
\begin{eqnarray*}
&&bv_{p-1}=tv_{p-1},\ bv_p=uv_p,\ bv_{p+1}=t^{-1}v_{p+1},\\
&&bv_k=v_k\ \text{if}\ k\notin\{p-1,p,p+1\}.
\end{eqnarray*}
Then $b\in \Sg$ and
\begin{eqnarray*}
&&bxb^{-1}v_{p-1}=t^{-1}u\frac{x_p}{\sqrt{2}}v_p,\ bxb^{-1}v_{p}=u^{-1}\frac{x_p}{\sqrt{2}}(tv_{p-1}+t^{-1}iv_{p+1}),\\
&&bxb^{-1}v_{p+1}=-itu\frac{x_p}{\sqrt{2}}v_p.
\end{eqnarray*}
Choose $x_p=-iu\sqrt{2}t$. Then $z(p-1)=\lim_{t\to 0} bxb^{-1}$. Thus $z(p-1)\in \l_{0\overline 1}$ is the limit of the elements of the $\L$-orbit passing through $\h_{\overline 1}$.

If  the involution $\Bbb D\ni a\to\overline a\in \Bbb D$ is trivial, then the only semisimple element in $\l_{0\overline 1}$ is $0$. 
 The restriction of $\Sg$ to the span of $v_{p-1}$, $v_p$ and $v_{p+1}$ is isomorphic to $\Og_1\times \Sp_2(\R)$, if $\Bbb D=\R$ and to $\Og_1\times \Sp_2(\C)$, if $\Bbb D=\C$. 
The complement of $0$ in  $\l_{0\overline 1}$ is a single nilpotent orbit under the action of this group.

For a fixed $j$ with $0< 2j<\frac{m+2}{2}$ and for $t>0$ define $b=b(t)\in\Sg$ by
\begin{eqnarray*}
&&bv_{2j-1}=tv_{2j-1},\ bv_{2j}=t^{-1}v_{2j},\\
&&bv_{m+3-2j}=t^{-1}v_{m+3+2j},\ bv_{m+2-2j}=tv_{m+2-2j},\\
&&bv_k=v_k\ \text{for}\ k\notin\{2j-1, 2j, m+3-2j, m+2-2j\}.
\end{eqnarray*}
Also, let $x=x(t)\in\h_{\overline 1}$ be such that $x_j=t^2$ and $x_i=0$ for $i\ne j$. Then $\lim_{t\to 0} bxb^{-1}=z(2j-1)$.
Hence, $z_\l=\sum_j z(2j-1)$ is the limit of elements of the $\L$-orbits passing through $\h_{\overline 1}$.

Notice that the same holds for the nilpotent element $z=\sum_{k=2}^{\frac{m+2}{2}}z(k)$, $z(k)\in\ss1$. 
\bigskip

\enlargethispage*{3truecm}
\underline{Case I.f: $\Sg=\Og_{p,p}\times \Sp_{2(p+1)}(\R)$ or $\Og_{2p}(\C)\times \Sp_{2(p+1)}(\C)$}

\begin{eqnarray*}
&&m\in 2\Bbb Z\setminus 4\Bbb Z,\ m\geq 0,\ \Bbb D\ne \Bbb H,\ \text{the involution}\ \Bbb D\ni a\to\overline a\in\Bbb D\ \text{is trivial};\\
&&\V=\sum_{k=0}^{m}(\Bbb D v_k+\Bbb D v'_k);\ v_{even},\ v'_{even}\in \V_{\overline 0};\ v_{odd},\ v'_{odd}\in \V_{\overline 1};\\
&&v_{j}=z^jv_0\ne 0,\ v'_{j}=z^jv'_0\ne 0,\ 0\leq j\leq m,\ zv_{m}=zv'_{m}=0;\\ 
&&( v_k, v'_{m-k}) =\delta(k),\ ( v'_k, v_{m-k})=-\delta(k),\ 0\leq k\leq m\\ 
&&\text{and all other pairings are zero};
\end{eqnarray*}
Here,
\begin{eqnarray*}
&&Tv_k=t_kv'_{m-k},\ Tv'_k=t'_kv_{m-k},\\ 
&&t_k=(-1)^k\delta(k),\ t'_k=-(-1)^k\delta(k),\ 0\leq k\leq m.
\end{eqnarray*}
The Cartan subspace $\h_{\overline 1}$ 
consists of the linear maps $x$ defined by
\begin{eqnarray*}
&&xv_{2j}=x_{j}v_{2j+1},\ xv_{2j+1}=x_{j}v_{2j},\ xv_{m}=0,\\ 
&&xv'_{m-2j}=-x_{j}v'_{m-2j-1},\ xv'_{m-2j-1}=x_{j}v'_{m-2j},\ xv'_0=0,\\
&&x_j\in\Bbb D,\  0\leq j\leq \frac{m}{2}.
\end{eqnarray*}
This is the direct sum of the indecomposable Cartan subspaces which are isomorphic to those which occur in \cite[Theorem 6.2(a)]{PrzebindaLocal}. (One has to adjust the sesquilinear forms $\tau_0$ and $\tau_1$ listed in \cite[Theorem 6.2(a)]{PrzebindaLocal} in order to get the forms $(\ ,\ )_{\overline 0}$ and $(\ ,\ )_{\overline 1}$ we are working with here.)

The Lie algebra $\a$ consists of the linear maps $a$ defined by
\begin{eqnarray*}
&&av_k=a_kv_{k},\ av'_{m-k}=-a_{k}v'_{m-k},\ a_k\in\R,\ 0\leq k<m,\\
&&av_{m}=0,av'_0=0,\\
&&a_{2j}=a_{2j+1}\ \text{for}\ 0\leq 2j< m.
\end{eqnarray*}

For $0\leq k<m$ define a linear map $z(k)$ by
\begin{eqnarray*}
&&z(k)v_k=v_{k+1},\ z(k)v'_{m-k-1}=v'_{m-k},\\
&&z(k)v_j=0,\ z(k)v'_{m-j-1}=0,\ j\ne k.
\end{eqnarray*}
Then $z=\sum_{k=0}^{m-1}z(k)$, $z(k)\in\ss1$ and for $a\in \a$, $[a,z(k)]=(a_{k+1}-a_k)z(k)$. In particular $z(even)\in \sum_{i=1}^m\l_{i\overline 1}$ and $z(odd)\in\n_1$.

For a fixed $j$ with $0\leq 2j\leq\frac{m}{2}$ and for $t>0$ define $b=b(t)\in\Sg$ by
\begin{eqnarray*}
&&bv_{2j}=tv_{2j},\ bv_{2j+1}=t^{-1}v_{2j+1},\\
&&bv_{m-2j}=t^{-1}v_{m-2j},\ bv_{m-2j-1}=tv_{m-2j-1},\\
&&bv_k=v_k\ \text{for}\ k\notin\{2j, 2j+1, m-2j, m-2j-1\}.
\end{eqnarray*}
Also, let $x=x(t)\in\h_{\overline 1}$ be such that $x_j=t^2$ and $x_i=0$ for $i\ne j$. Then $\lim_{t\to 0} bxb^{-1}=z(2j)$.
Hence, $z_\l=\sum_j z(2j)$ is the limit of elements of the $\L$-orbits passing through $\h_{\overline 1}$.
\bigskip

\enlargethispage*{3truecm}
\underline{Case I.g. $\Sg=\Og_{p,p}\times \Sp_{2p}(\R)$, $\Ug_{p,p}\times \Ug_{p,p}$ or $\Sp_{p,p}\times \Og^*_{4p}$}

\begin{eqnarray*}
&&m\in 2\Bbb Z+1;\\
&&\V=\sum_{k=0}^{m}(\Bbb D v_k+\Bbb D v'_{k+1});\ v_{even},\ v'_{even}\in \V_{\overline 0};\ v_{odd},\ v'_{odd}\in \V_{\overline 1};\\
&&v_{j}=z^jv_0\ne 0,\ v'_{j+1}=z^jv'_1\ne 0,\ 0\leq j\leq m,\ zv_{m}=zv'_{m+1}=0;\\ 
&&( v_k, v'_{m+1-k}) =\delta(k),\ ( v'_{k+1}, v_{m-k})=\delta(k+1)\delta(m),\ 0\leq k\leq m\\ 
&&\text{and all other pairings are zero};
\end{eqnarray*}
Here,
\begin{eqnarray*}
&&Tv_k=t_kv'_{m+1-k},\ Tv'_{k+1}=t'_{k+1}v_{m-k},\\ 
&&t_k=\delta(k+1),\ t'_{k+1}=\delta(m+1-k),\ 0\leq k\leq m.
\end{eqnarray*}
The Cartan subspace $\h_{\overline 1}$ 
consists of the linear maps $x$ defined by
\begin{eqnarray*}
&&xv_{2j}=x_{j}v_{2j+1},\ xv_{2j+1}=x_{j}v_{2j},\\ 
&&xv'_{m+1-2j}=-\overline{x_j}v'_{m-2j},xv'_{m-2j}=\overline{x_j}v'_{m+1-2j},\\
&&\text{if the involution}\ \Bbb D\ni a \to \overline{a}\in\Bbb D\ \text{is trivial, then}\ x_j\in\Bbb D,\\
&&\text{otherwise}\  x_j\in\C\subseteq \Bbb D,\  0\leq j\leq \frac{m}{2}.
\end{eqnarray*}
This is the direct sum of the indecomposable Cartan subspaces which are isomorphic to those which occur in \cite[Theorem 6.2(a)]{PrzebindaLocal}. (One has to adjust the sesquilinear forms $\tau_0$ and $\tau_1$ listed in \cite[Theorem 6.2(a)]{PrzebindaLocal} in order to get the forms $(\ ,\ )_{\overline 0}$ and $(\ ,\ )_{\overline 1}$ we are working with here.)

The Lie algebra $\a$ consists of the linear maps $a$ defined by
\begin{eqnarray*}
&&av_k=a_kv_{k},\ av'_{m+1-k}=-a_{k}v'_{m+1-k},\ a_k\in\R,\ 0\leq k\leq m,\\
&&a_{2j}=a_{2j+1}\ \text{for}\ 0\leq 2j< m.
\end{eqnarray*}

For $0\leq k<m$ define a linear map $z(k)$ by
\begin{eqnarray*}
&&z(k)v_k=v_{k+1},\ z(k)v'_{m-k}=v'_{m+1-k},\\
&&z(k)v_j=0,\ z(k)v'_{m+1-j}=0,\ j\ne k.
\end{eqnarray*}
Then $z=\sum_{k=0}^{m-1}z(k)$, $z(k)\in\ss1$ and for $a\in \a$, $[a,z(k)]=(a_{k+1}-a_k)z(k)$. In particular $z(even)\in \sum_{i=1}^m\l_{i\overline 1}$ and $z(odd)\in\n_1$.

For a fixed $j$ with $0\leq j\leq\frac{m}{2}$ and for $t>0$ define $b=b(t)\in\Sg$ by
\begin{eqnarray*}
&&bv_{2j}=tv_{2j},\ bv_{2j+1}=t^{-1}v_{2j+1},\\
&&bv_{m+1-2j}=t^{-1}v_{m+1-2j},\ bv_{m-2j}=tv_{m-2j},\\
&&bv_k=v_k\ \text{for}\ k\notin\{2j, 2j+1, m+1-2j, m-2j\}.
\end{eqnarray*}
Also, let $x=x(t)\in\h_{\overline 1}$ be such that $x_j=t^2$ and $x_i=0$ for $i\ne j$. Then $\lim_{t\to 0} bxb^{-1}=z(2j)$.
Hence, $z_\l=\sum_j z(2j)$ is the limit of elements of the $\L$-orbits passing through $\h_{\overline 1}$.
\bigskip


\underline{Case I.``V ": $\Sg=\Og_2\times \Sp_2(\R)$;}

\begin{eqnarray*}
&&\Bbb D=\R;\\
&&\V=\Bbb D v_1+\Bbb D v_2+\Bbb D v_3+\Bbb D v'_2;\\
&&zv_1=v_2,\ zv_2=v_3,\ zv_3=0,\ zv'_2=0;\\ 
&&(v_1,v_3)=1,\ (v_2,v_2)=(v'_2,v'_2)=-1\ \text{and all the other pairings are zero};
\end{eqnarray*}
The Cartan subspace $\h_{\overline 1}$ 
consists of the linear maps $x$ defined by

\begin{eqnarray*}
&&x\frac{1}{\sqrt{2}}(v_2+v'_2)=x_1(v_1+v_3),\\ 
&&x\frac{1}{\sqrt{2}}(-v_2+v'_2)=x_1(v_1-v_3),\\ 
&&xv_1=x_1\sqrt{2}v_2,\\ 
&&xv_3=-x_1\sqrt{2}v'_2.
\end{eqnarray*}
where $x_1\in\R$, or equivalently,
\begin{eqnarray*}
&&xv_2=x_1v_3,\ xv'_2=x_1v_1,\\
&&xv_1=x_1v_2,\ xv_3=-x_1v'_2.
\end{eqnarray*}
This is the Cartan subspace which is isomorphic to the one which occurs in \cite[Theorem 6.2(c)]{PrzebindaLocal}.

For $t>0$ define $b=b(t)\in\Sg$ by
\begin{eqnarray*}
&&bv_1=tv_1,\ bv_2=v_2,\ bv_3=t^{-1}v_3,\ bv'_2=v'_2.
\end{eqnarray*}
Also, let $x_1=t$. Then $\lim_{t\to 0} bxb^{-1}=z$.
\bigskip


\underline{Case I.``V V": $\Sg=\Og_2\times \Sp_4\supseteq (\Og_1\times \Sp_2)\times (\Og_1\times \Sp_2)$}

\begin{eqnarray*}
&&\Bbb D=\R\ \text{or}\ \C\ \text{with trivial involution};\\
&&\V=(\Bbb D v_1+\Bbb D v_2+\Bbb D v_3)+(\Bbb D v'_1+\Bbb D v'_2+\Bbb D v'_3);\\
&&zv_1=v_2,\ zv_2=v_3,\ zv_3=0,\ zv'_1=v'_2,\ zv'_2=v'_3,\ zv'_3=0;\\ 
&&(v_1,v_3)=(v'_1,v'_3)=1,\ (v_2,v_2)=(v'_2,v'_2)=-1\ \text{and all the other pairings are zero};
\end{eqnarray*}
In this case there is only one, up to conjugation, Cartan subspace $\h_{\overline 1}\subseteq \ss1$. If $\Bbb D=\R$, then $\h_{\overline 1}$ may be realized as the space consisting of linear maps $x$ defined by
\begin{eqnarray*}
&&x\frac{1}{\sqrt{2}}(v_2+v'_2)=x_1(v_1+v_3),\\ 
&&x\frac{1}{\sqrt{2}}(-v_2+v'_2)=x_1(v_1-v_3),\\ 
&&xv_1=x_1\sqrt{2}v_2,\\ 
&&xv_3=-x_1\sqrt{2}v'_2.
\end{eqnarray*}
where $x_1\in\R$, see \cite[Theorem 6.2(c)]{PrzebindaLocal}.  Equivalently,
\begin{eqnarray*}
&&xv_2=x_1v_3,\ xv'_2=x_1v_1,\\
&&xv_1=x_1v_2,\ xv_3=-x_1v'_2.
\end{eqnarray*}
If $\Bbb D=\C$, then $\h_{\overline 1}$ is isomorphic to the Cartan subspace comprised of elements which occurs in \cite[Theorem 6.2(a)]{PrzebindaLocal}.
In any case the kernel of a non-zero element of $\h_{\overline 1}$  is contained in $\V_{\overline 1}$. However any semisimple orbit in $\ss1$ passes through a Cartan subspace.
Hence, the kernel of any non-zero semisimple element of $\ss1$ is contained in $\V_{\overline 1}$. Let $x\in\ss1$ and let $x=x_s+x_n$ be Jordan decomposition of $x$, see \cite[]{PrzebindaLocal}. Since $x_n\in\ss1$ commutes with $x_s\in\ss1$ we see from the above that $x_n=0$. Thus $x$ is either nilpotent or semisimple.  
Since there are finitely many nilpotent orbits in $\ss1$, our given nilpotent $z$ is the limit of elements of the $\Sg$-orbits passing through $\h_{\overline 1}$. Similar argument applies to the previous case too. 
\bigskip


Next we consider the dual pairs of type II. As for the dual pairs of type I, $x\in\End(\V)$ belongs to $\ss1$ if and only if $(xu,v)=(u,sxv)$ for all $u,v\in\V$, where $s(v_{even}+v_{odd})=v_{even}-v_{odd}$. Also, if $\Bbb D\ne \R$, then there is a nontrivial involution on $\Bbb D$.
In any case we'll describe a basis on $\V$ consisting of vectors $v_k$. Then there is a positive definite hermitian (or symmetric if $\Bbb D=\R$) form $\eta(\ ,\ )$ on $\V$ defined by the condition $\eta(v_k,v_k)=1$. This form determines an involution $\End(\V)\ni x\to x^\dag\in \End(\V)$ defined by $\eta(xu,v)=\eta(u,x^\dag v)$. 
Then $\theta(x)=-x^\dag$ if $x\in \so$ and $\theta(x)=sx^\dag$ if $x\in \ss1$.\newline
\newline
\underline{\mbox{Case II.a:}}

\begin{eqnarray*}
&&m\in \Bbb Z,\ m>0;\\
&&\V=\sum_{k=0}^m\Bbb D v_k,\ v_{even}\in \V_{\overline 0},\ v_{odd}\in \V_{\overline 1};\\
&&v_k=z^kv_0\ne 0, 0\leq k\leq m, zv_m=0;
\end{eqnarray*}
If $m$ is odd, then the Cartan subspace $\h_{\overline 1}$ 
consists of the linear maps $x$ defined by
\begin{eqnarray*}
&&xv_{2j}=x_{j}v_{2j+1},\ xv_{2j+1}=x_{j}v_{2j},\\ 
&&\text{if}\ \Bbb D\ne \Bbb H,\ \text{then}\ \ x_j\in\Bbb D,\\
&&\text{if}\ \Bbb D=\Bbb H,\ \text{then}\ x_j\in \C\subseteq \Bbb H,\  0\leq 2j<m.
\end{eqnarray*}
If $m$ is even, then the Cartan subspace $\h_{\overline 1}$ 
consists of the linear maps $x$ defined by
\begin{eqnarray*}
&&xv_{2j+1}=x_{j}v_{2j+2},\ xv_{2j+2}=x_{j}v_{2j+1},\ xv_0=0,\\ 
&&\text{if}\ \Bbb D\ne \Bbb H,\ \text{then}\ \ x_j\in\Bbb D,\\
&&\text{if}\ \Bbb D=\Bbb H,\ \text{then}\ x_j\in \C\subseteq \Bbb H,\  0\leq 2j<m.
\end{eqnarray*}
These are the direct sums of the indecomposable Cartan subspaces which occur in \cite[Theorem 6.2(e)]{PrzebindaLocal}.

If $m$ is odd, then the Lie algebra $\a$ consists of the linear maps $a$ defined by
\begin{eqnarray*}
&&av_k=a_kv_{k},\ a_k\in\R,\ 0\leq k\leq m,\\
&&a_{2j}=a_{2j+1}\ \text{for}\ 0\leq 2j< m.
\end{eqnarray*}
If $m$ is even, then the Lie algebra $\a$ consists of the linear maps $a$ defined by
\begin{eqnarray*}
&&av_k=a_kv_{k},\ a_k\in\R,\ 1\leq k\leq m,\ av_0=0\\
&&a_{2j+1}=a_{2j+2}\ \text{for}\ 0\leq 2j< m.
\end{eqnarray*}

For $0\leq k<m$ define a linear map $z(k)$ by
\begin{eqnarray*}
&&z(k)v_k=v_{k+1},\\
&&z(k)v_j=0,\ j\ne k.
\end{eqnarray*}
Then $z=\sum_{k=0}^{m-1}z(k)$, $z(k)\in\ss1$ and for $a\in \a$, $[a,z(k)]=(a_{k+1}-a_k)z(k)$. If $m$ is odd, then $z(even)\in \sum_{i=1}^m\l_{i\overline 1}$ and $z(odd)\in\n_1$.
 If $m$ is odd, then $z(odd)\in \sum_{i=1}^m\l_{i\overline 1}$ and $z(even)\in\n_1$.

Suppose $m$ is odd. For a fixed $j$ with $0\leq 2j<m$ and for $t>0$ define $b=b(t)\in\Sg$ by
\begin{eqnarray*}
&&bv_{2j}=tv_{2j},\ bv_{2j+1}=t^{-1}v_{2j+1},\\
&&bv_k=v_k\ \text{for}\ k\notin \{2j, 2j+1\}.
\end{eqnarray*}
Also, let $x=x(t)\in\h_{\overline 1}$ be such that $x_j=t^2$ and $x_i=0$ for $i\ne j$. Then $\lim_{t\to 0}bxb^{-1}=z(2j)$.
Hence, $z_\l=\sum_j z(2j)$ is the limit of elements of the $\L$-orbits passing through $\h_{\overline 1}$.

Suppose $m$ is even. For a fixed $j$ with $0\leq 2j<m$ and for $t>0$ define $b=b(t)\in\Sg$ by
\begin{eqnarray*}
&&bv_{2j+1}=tv_{2j+1},\ bv_{2j+2}=t^{-1}v_{2j+2},\\
&&bv_k=v_k\ \text{for}\ k\notin \{2j+1, 2j+2\}.
\end{eqnarray*}
Also, let $x=x(t)\in\h_{\overline 1}$ be such that $x_j=t^2$ and $x_i=0$ for $i\ne j$. Then $\lim_{t\to 0} bxb^{-1}=z(2j+1)$.
Hence, $z_\l=\sum_j z(2j+1)$ is the limit of elements of the $\L$-orbits passing through $\h_{\overline 1}$.

\newpage

\underline{\mbox{Case II.b:}}

\begin{eqnarray*}
&&m\in \Bbb Z,\ m>0;\\
&&\V=\sum_{k=1}^{m+1}\Bbb D v_k,\ v_{even}\in \V_{\overline 0},\ v_{odd}\in \V_{\overline 1};\\
&&v_{k+1}=z^kv_1\ne 0, 0\leq k\leq m, zv_{m+1}=0;
\end{eqnarray*}
If $m$ is odd, then the Cartan subspace $\h_{\overline 1}$ 
consists of the linear maps $x$ defined by
\begin{eqnarray*}
&&xv_{2j+1}=x_{j}v_{2j+2},\ xv_{2j+2}=x_{j}v_{2j+1},\ xv_{m+1}=0,\\ 
&&\text{if}\ \Bbb D\ne \Bbb H,\ \text{then}\ \ x_j\in\Bbb D,\\
&&\text{if}\ \Bbb D=\Bbb H,\ \text{then}\ x_j\in \C\subseteq \Bbb H,\  0\leq 2j<m.
\end{eqnarray*}
If $m$ is even, then the Cartan subspace $\h_{\overline 1}$ 
consists of the linear maps $x$ defined by
\begin{eqnarray*}
&&xv_{2j+1}=x_{j}v_{2j+2},\ xv_{2j+2}=x_{j}v_{2j+1},\\ 
&&\text{if}\ \Bbb D\ne \Bbb H,\ \text{then}\ \ x_j\in\Bbb D,\\
&&\text{if}\ \Bbb D=\Bbb H,\ \text{then}\ x_j\in \C\subseteq \Bbb H,\  0\leq 2j<m.
\end{eqnarray*}
These are the direct sums of the indecomposable Cartan subspaces which occur in \cite[Theorem 6.2(e)]{PrzebindaLocal}.

If $m$ is odd, then the Lie algebra $\a$ consists of the linear maps $a$ defined by
\begin{eqnarray*}
&&av_k=a_kv_{k},\ a_k\in\R,\ 0\leq k\leq m,\ av_{m+1}=0,\\
&&a_{2j+1}=a_{2j+2}\ \text{for}\ 0\leq 2j< m.
\end{eqnarray*}
If $m$ is even, then the Lie algebra $\a$ consists of the linear maps $a$ defined by
\begin{eqnarray*}
&&av_k=a_kv_{k},\ a_k\in\R,\ 1\leq k\leq m,\\
&&a_{2j+1}=a_{2j+2}\ \text{for}\ 0\leq 2j< m.
\end{eqnarray*}

For $0\leq k<m$ define a linear map $z(k)$ by
\begin{eqnarray*}
&&z(k)v_k=v_{k+1},\\
&&z(k)v_j=0,\ j\ne k.
\end{eqnarray*}
Then $z=\sum_{k=0}^{m-1}z(k)$, $z(k)\in\ss1$ and for $a\in \a$, $[a,z(k)]=(a_{k+1}-a_k)z(k)$. Then $z(odd)\in \sum_{i=1}^m\l_{i\overline 1}$ and $z(even)\in\n_1$.

For a fixed $j$ with $0\leq 2j<m$ and for $t>0$ define $b=b(t)\in\Sg$ by
\begin{eqnarray*}
&&bv_{2j+1}=tv_{2j+1},\ bv_{2j+2}=t^{-1}v_{2j+2},\\
&&bv_k=v_k\ \text{for}\ k\notin \{2j+1, 2j+2\}.
\end{eqnarray*}
Also, let $x=x(t)\in\h_{\overline 1}$ be such that $x_j=t^2$ and $x_i=0$ for $i\ne j$. Then $\lim_{t\to 0} bxb^{-1}=z(2j+1)$.
Hence, $z_\l=\sum_j z(2j+1)$ is the limit of elements of the $\L$-orbits passing through $\h_{\overline 1}$.
\vskip 1truecm

By combining the cases Case I.a - Case II.b we see that Proposition \ref{good nilpotent elements} holds if $(z,\V)$ is indecomposable. 
Every nilpotent $(z,\V)$ is a finite direct sum of indecomposable nilpotents $(z_i,\V_i)$, \cite[Definition 3.14]{DaszKrasPrzebindaK-S2}. Suppose
$\Sg$ is not isomorphic to an ortho-symplectic pair ($\Og_{p,q}\times\Sp_{2n}(\R)$ or $\Og_{p}(\C)\times\Sp_{2n}(\C)$). Then each
$(z_i,\V_i)$ has the same property (the group $\Sg|_{\V_i}$ is not ortho-symplectic) and the Cartan subspace $\h_{\overline 1}$ can be defined as the direct sum of the Cartan subspaces constructed for each $(z_i,\V_i)$. Also the involution $\theta$ may be extended from each $\s(\V_i)$ to $\s(\V)$.
Hence, Proposition \ref{good nilpotent elements} holds for $(z,\V)$, which is the sum of the $(z_i,\V_i)$.

Suppose $\Sg$ is an ortho-symplectic pair. Let
\begin{eqnarray*}
(z,\V)=(z_1,\V_1) \oplus \dots \oplus (z_{r},\V_{r})\oplus (0,\V_0)
\end{eqnarray*}
be the decomposition into irreducibles, where each $z_i$ is non-zero. 
As we checked in Case I.c - Case I.g, the proposition holds for each $(z_i,\V_i)$. However this does not automatically imply the proposition for the sum, because the sum of the individual Cartan subspaces constructed for the $(z_i,\V_i)$ may be too small (could be zero). This problem does not occur if $dim_\Bbb D\,\V_{i\overline 0}$ is even for each $i$. Thus we may assume that $dim_\Bbb D\,\V_{i\overline 0}$ is odd for each $i$.

If $dim_\Bbb D\,\V_{i\overline 0}>dim_\Bbb D\,\V_{i\overline 1}$, then we are in Case I.a, and if $dim_\Bbb D\,\V_{i\overline 0}=dim_\Bbb D\,\V_{i\overline 1}$, then we are in Case I.g, the proposition holds for the sum of such $(z_i,\V_i)$ with the Cartan subspace equal to the sum of the individual Cartan subspaces.
 
Thus we may assume that $dim_\Bbb D\,\V_{i\overline 0}<dim_\Bbb D\,\V_{i\overline 1}$ for each $i$. Thus each $(z_i,\V_i)$ is as in Case I.e. Here we combine Case I.e with either Case I.``V V" or Case I.``V ." to construct the Cartan subspace and the involution $\theta$ for the sum.
\section*{Appendix C}
\label{appendix:limit-ss}
\setcounter{thh}{0}
\renewcommand{\thethh}{C.\fontindex{thh}}
\setcounter{equation}{0}
\renewcommand{\theequation}{C.\fontindex{equation}}
\renewcommand{\theequation}{C.\fontindex{equation}}
In this appendix we conclude the proof of Theorem \ref{density}. Let $(\Sg,\s)$ be a supergroup associated with a dual pair 
$(\Ug_n,\Ug_n)$ or $(\GL_n(\DD),\GL_n(\DD))$ with $\DD=\R$, $\C$ or $\Ha$. Let $x=x_s+x_n$ be the Jordan decomposition of an 
element $x \in \ss1$. We suppose that $x_s\neq 0$ and $x_n\neq 0$. Moreover, we may assume that $x_s$ belongs to a
isotypic Cartan subalgebra built up from indecomposable blocks as in \cite[Theorems 5.2(b) and 5.3]{PrzebindaLocal}. 

\medskip 

\noindent\underline{$\Sg=\Ug_n \times \Ug_n$:} \\
$\V=\V_{\overline 0} \times \V_{\overline 1}$ with $\V_\alpha=\sum_{k=1}^n \C v_{\alpha,k}$ \quad ($\alpha \in \Zb/2\Zb)$. \\
$x_s=x_s(a)=\begin{pmatrix} 0 & a\Id \\ a\Id & 0\end{pmatrix}$ with $a \in \C\setminus \{0\}$ and $\Id$ the $n\times n$ identity matrix. \\
$x_n=\begin{pmatrix} 0 & w \\ w* & 0\end{pmatrix}$ with $w\in \M_n(\C)$ and $w^*=i\bar{w}^t$.\\
Since $x_n x_s=x_s x_n$, the matrix $w$ must satisfy $w^*=w$. Hence $x_n^2=\begin{pmatrix} w^2 & 0 \\ 0 & w^2\end{pmatrix} \in \so=\u_n\times \u_n$.
Since $x_n^2$ is nilpotent, so must be $w^2$. But every matrix in $\u_n$ is diagonalizable. Hence $w^2=0$. Taking traces, we obtain 
$i{\rm tr}(\bar{w}^tw)={\rm tr}(w^2)=0$, which implies $w=0$. This shows that in this case every element in $\ss1$ is either nilpotent or semisimple.

\medskip
\noindent\underline{$\Sg=\GL_n(\C) \times \GL_n(\C)$:} \\
$\V$ and $x_s=x_s(a)$ are as in the previous case. $x_n=\begin{pmatrix} 0 & w \\ w' & 0\end{pmatrix}$ with $w, w'\in \M_n(\C)$.\\
From $x_n x_s=x_s x_n$ we obtain that $w'=w$. Thus $x=x_s+x_n=\begin{pmatrix} 0 & a\Id+w \\ a\Id+w & 0\end{pmatrix}$.
By the density of the semisimple matrices in $\M_n(\C)$, we can find $y_m$ semisimple so that $\lim_{m\to \infty} y_m=a\Id+w$. Hence
$x=\lim_{m\to \infty}  \begin{pmatrix} 0 & y_m \\ y_m & 0\end{pmatrix}$. Notice that the latter matrix is semisimple. Indeed, if 
$y\in \M_n(\C)$ is semisimple and $gyg^{-1}=d$ is a diagonal matrix, then $\begin{pmatrix} g & 0\\0  & g\end{pmatrix}
\begin{pmatrix} 0 & y \\ y & 0\end{pmatrix}\begin{pmatrix} g^{-1} & 0\\0  & g^{-1}\end{pmatrix} =\begin{pmatrix} 0 & d\\ d & 0\end{pmatrix}$
is also diagonalizable.

\medskip
\noindent\underline{$\Sg=\GL_n(\R) \times \GL_n(\R)$:} \\
There are three cases to consider.
\begin{enumerate}
\renewcommand{\theenumi}{\Roman{enumi}}
\item The analogue of the case of $\GL_n(\C) \times \GL_n(\C)$, with $\C$ replaced by $\R$. One proceeds as above.
\item 
$\V=\V_{\overline 0} \times \V_{\overline 1}$ with $\V_\alpha=\sum_{k=1}^n \R v_{\alpha,k}$ \qquad ($\alpha \in \Zb/2\Zb)$. \\
$x_s=x_s(a)=\begin{pmatrix} 0 & a\Id \\ -a\Id & 0\end{pmatrix}$ with $a \in \R\setminus \{0\}$. \\
$x_n=\begin{pmatrix} 0 & w \\ w' & 0\end{pmatrix}$ with $w, w'\in \M_n(\R)$.\\
From $x_n x_s=x_s x_n$ we obtain $w'=-w$. Hence $x=x_s+x_n=\begin{pmatrix} 0 & a\Id+w \\ -a\Id-w & 0\end{pmatrix}$.
Let $y_m \in \M_n(\R)$ be semisimple and so that $\lim_{m\to \infty} y_m=a\Id+w$. Then 
$x=\lim_{m\to \infty}  \begin{pmatrix} 0 & y_m \\ -y_m & 0\end{pmatrix}$. Notice that the latter matrix is semisimple. Indeed, over $\C$, if 
$gyg^{-1}=d$ is a diagonal matrix, then $\begin{pmatrix} g & 0\\0  & -g\end{pmatrix}
\begin{pmatrix} 0 & y \\ -y & 0\end{pmatrix}\begin{pmatrix} g^{-1} & 0\\0  & (-g)^{-1}\end{pmatrix} =\begin{pmatrix} 0 & d\\ -d & 0\end{pmatrix}$
is also diagonalizable.
\item 
$\V=\V_{\overline 0} \times \V_{\overline 1}$ with $\V_\alpha=\sum_{k=1}^n (\R v_{\alpha,k}+\R v'_{\alpha,k})$ \qquad ($\alpha \in \Zb/2\Zb)$. \\
$x_s=x_s(a)=\begin{pmatrix} 0 & A \\ A & 0\end{pmatrix}$ where $A={\rm diag}(a,\dots,a)$ is a block diagonal matrix with equal $2\times 2$  
diagonal blocks $a=\begin{pmatrix} \beta & \gamma \\ -\gamma & \beta\end{pmatrix}$, $\beta,\gamma \in\R$. We can assume $\beta\neq 0$ and $\gamma\neq 0$, otherwise we are reduced to the previous cases. \\
$x_n=\begin{pmatrix} 0 & w \\ w' & 0\end{pmatrix}$ with $w, w'\in \M_n(\R)$.\\
From $x_n x_s=x_s x_n$ we obtain $w'A=Aw$ and $Aw'=wA$. This implies that $w'=A^{-1}wA$ and $A^2w=w A^2$.
Notice that $$a^2=\begin{pmatrix} \beta^2-\gamma^2 & 2\beta\gamma \\ -2\beta\gamma & \beta^2-\gamma^2\end{pmatrix}=\Id+2\beta\gamma j$$
where $j=\begin{pmatrix} 0 & 1\\ -1 & 0\end{pmatrix}\,$. 
Write $w=(w_{rs})$ where each $w_{rs}$ is a $2\times 2$ block. 
As $A^2w=w A^2$, we have that $Jw=wJ$ where $J={\rm diag}(j,\dots,j)$. Hence $Jw_{rs}=w_{rs}J$, i.e. $w_{rs}=\begin{pmatrix} \beta_{rs} & \gamma_{rs} \\ -\gamma_{rs} & \beta_{rs}\end{pmatrix}$. Therefore $w$ commutes with $A$, i.e. $w'=w$.
The conclusion follows by the same argument as before.
\end{enumerate}

\medskip
\underline{$\Sg=\GL_n(\H) \times \GL_n(\H)$:} \\
There is only one case, and it is as (III) for $\GL_n(\R) \times \GL_n(\R)$.

\section*{Appendix D}
\label{appendix:IN}
\setcounter{thh}{0}
\renewcommand{\thethh}{D.\fontindex{thh}}
\setcounter{equation}{0}
\renewcommand{\theequation}{D.\fontindex{equation}}
\renewcommand{\theequation}{D.\fontindex{equation}}
Recall the function $I_N$, (\ref{3.4}). Here we verify the following elementary lemma. 
\begin{lem}\label{A.4}
For $N>0$,
\begin{equation}\label{A.4.a}
\int_0^\infty (1+(a^2+a^{-2})t^2)^{-N}\,da/a\leq \frac{2}{N}I_N(t)\qquad (t>0).
\end{equation}
Furthermore, 
\begin{equation}\label{A.4.b}
\int_0^\infty (1+a^2b^2+a^{-2}c^2)t^2)^{-N}\,da/a\leq \frac{2}{N}I_N(\sqrt{bc}) \qquad (b,c>0).
\end{equation}
\end{lem}
\begin{prf}
Since the measure $da/a$ is invariant under the substitution, $a\to at$, $t>0$, (\ref{A.4.b}) follows from (\ref{A.4.a}). Since $a^2+\overline a^{2}\geq 2>1$, the left hand side of (\ref{A.4.a}) is dominated by
\begin{eqnarray}\label{A.5}
&&(1+t^2)^{-N/2}\int_0^\infty (1+(a^2+a^{-2})t^2)^{-N/2}\,da/a\\
&\leq&(1+t^2)^{-N/2}2\int_1^\infty (1+(at)^2)^{-N/2}\,da/a\nn\\
&=&(1+t^2)^{-N/2}2\int_t^\infty (1+a^2)^{-N/2}\,da/a.\nn
\end{eqnarray}
If $t\geq 1$ then the last expression in (\ref{A.5}) is less or equal to
\begin{eqnarray*}
t^{-N}2\int_1^\infty (1+a^2)^{-N/2}\,da/a\leq t^{-N}2\int_1^\infty a^{-N-1}\,da=\frac{2}{N} t^{-N}.
\end{eqnarray*}
If $0<t< 1$ then the last expression in (\ref{A.5}) is dominated by
\begin{eqnarray*}
&&2\int_t^\infty (1+a^2)^{-N/2}\,da/a\\
&\leq& 2\int_1^\infty (1+a^2)^{-N/2}\,da/a
+2\int_t^1\,da/a \leq \frac{2}{N}-2\,\ln(t)=\frac{2}{N}(1-N\,\ln(t)).
\end{eqnarray*}
\end{prf}


\end{document}